\documentclass[11pt, reqno]{amsart}
\usepackage{latexsym}
\usepackage{bbm}
\usepackage{amsthm}
\usepackage{eepic}
\usepackage{graphicx}
\usepackage{amssymb}
\usepackage{amsfonts}
\usepackage{amsmath}
\usepackage{fourier}
\usepackage{enumerate}
\usepackage{color}

\def\UseSection{
        \numberwithin{equation}{section}
    \theoremstyle{plain}
        \newtheorem{theorem}    {Theorem}[section]
        \DefineTheorems 
}

\def\DefineTheorems{
    
    \newtheorem{lemma}      [theorem] {Lemma}
    
    \newtheorem{prop}       [theorem] {Proposition}
    
    \newtheorem{cor}        [theorem] {Corollary}

    \theoremstyle{definition}
    \newtheorem{defn}       [theorem] {Definition}

    \newtheorem{rk}     [theorem] {Remark}

    \theoremstyle{definition}

}

\newcommand{\bt}   {\begin{theorem}}
\newcommand{\et}   {\end  {theorem}}
\newcommand{\bl}   {\begin{lemma}}
\newcommand{\el}   {\end  {lemma}}
\newcommand{\bp}   {\begin{prop}}
\newcommand{\ep}   {\end  {prop}}
\newcommand{\bc}   {\begin{cor}}
\newcommand{\ec}   {\end  {cor}}
\newcommand{\bd}   {\begin{defn}}
\newcommand{\ed}   {\end  {defn}}

\newcommand{\ba}   {\begin{array}}
\newcommand{\ea}   {\end  {array}}
\newcommand{\be}   {\begin{enumerate}}
\newcommand{\ee}   {\end  {enumerate}}
\newcommand{\bi}   {\begin{itemize}}
\newcommand{\ei}   {\end  {itemize}}

\def\eq#1\en{\begin{equation}#1\end{equation}}
\def\eqsplit#1\ensplit{
    \begin{equation}\begin{split}#1\end{split}\end{equation}
    }
\def\eqalign#1\enalign{
    \begin{align}#1\end{align}
    }
\def\eqmul#1\enmul{
    \begin{multline}#1\end{multline}
    }
\newcommand{\eqarrstar} {\begin{eqnarray*}}
\newcommand{\enarrstar} {\end{eqnarray*}}
\newcommand{\eqarray}   {\begin{eqnarray}}
\newcommand{\enarray}   {\end{eqnarray}}
\newcommand{\nnb}   {\nonumber \\}

\newcommand{\lbeq}[1]  {\label{e:#1}}
\newcommand{\refeq}[1] {\eqref{e:#1}}    

\makeatletter
\newcommand{\labelcounter}[2]{{%
    \stepcounter{#1}
    \protected@write\@auxout{}%
    {\string\newlabel{#2}{{\csname the#1\endcsname}{\thepage}}}%
    {\ref{#2}}
    }}
\makeatother
\newcommand{\sss}   { \scriptscriptstyle }

\newcommand{\Ebold} {{\mathbb E}}

\newcommand{\Pbold} {{\mathbb P}}

\newcommand{\Zbold} {{\mathbb Z}}

\newcommand{\Acal}   {\mathcal{A}}
\newcommand{\Bcal}   {\mathcal{B}}
\newcommand{\Ccal}   {\mathcal{C}}
\newcommand{\Dcal}   {\mathcal{D}}
\newcommand{\Ecal}   {\mathcal{E}}

\newcommand{\Gcal}   {\mathcal{G}}

\newcommand{\Mcal}   {\mathcal{M}}

\newcommand{\Pcal}   {\mathcal{P}}

\newcommand{\Rcal}   {\mathcal{R}}
\newcommand{\Scal}   {\mathcal{S}}
\newcommand{\Tcal}   {\mathcal{T}}
\newcommand{\Ucal}   {\mathcal{U}}
\newcommand{\Vcal}   {\mathcal{V}}

\newcommand{\Xcal}   {\mathcal{X}}
\newcommand{\Ycal}   {\mathcal{Y}}
\newcommand{\Zcal}   {\mathcal{Z}}

\newcommand{\Zd}    {{ {\Zbold}^d }}

\newcommand{\spose}[1] {{\hbox to 0pt{#1\hss}} }
\newcommand{\ltapprox} {\mathrel{\spose{\lower 3pt\hbox{$\mathchar"218$}}
 \raise 2.0pt\hbox{$\mathchar"13C$}}}
\newcommand{\gtapprox} {\mathrel{\spose{\lower 3pt\hbox{$\mathchar"218$}}
 \raise 2.0pt\hbox{$\mathchar"13E$}}}

\newcommand{\nin}  {{ \not\in }}

\UseSection  
\newcommand{\R}     {\mathbb{R}}
\newcommand{\N}     {\mathbb{N}}
\renewcommand{\P}   {\mathbb{P}}

\newcommand{\E}     {\mathbb{E}}
\newcommand{\Z}     {\mathbb{Z}}

\newcommand\1{\mathbbm{1}}

\newcommand{\ua}{\nearrow}

\newcommand{\Qp}{\mathbb{Q}_{p}}
\newcommand{\Qiic}{\mathbb{Q}_{\mathsf{ \sss IIC}}}
\newcommand{\Piic}{\mathbb{P}_{\mathsf{ \sss IIC}}}
\newcommand{\Eiic}{\mathbb{E}_{\mathsf{ \sss IIC}}}
\newcommand{\Pp}{\mathbb{P}_p}
\newcommand{\Ppc}{\mathbb{P}_{p_c}}
\newcommand{\events}{\mathfrak{F}}
\newcommand{\conn}{\leftrightarrow}

\newcommand{\Conn}{\Longleftrightarrow}

\newcommand{\edges}{E}

\newcommand{\vep}{\varepsilon}

\newcommand{\egr}\Exp

\newcommand{\Exp}{{{\mathbb E}_p}}

\def\1{{\mathchoice {1\mskip-4mu\mathrm l}      
{1\mskip-4mu\mathrm l}
{1\mskip-4.5mu\mathrm l} {1\mskip-5mu\mathrm l}}}

\newcommand{\tbar}{\bar{\triangle}}

\newcommand{\indi}{\mathbbm{1}}
\newcommand{\ole}{\overline{e}}
\newcommand{\ule}{\underline{e}}

\newcommand{\ulb}{\underline{b}}
\newcommand{\olb}{\overline{b}}
\newcommand{\limp}{\lim_{p \ua p_{c}}}
\newcommand{\taupc}{\tau_{p_c}}
\newcommand{\Be}{\Ecal_r}
\newcommand{\Ee}{\Ecal_r}
\newcommand{\twa}{(2 \wedge \alpha)}

\newcommand{\QR}{Q_{r}\cap \Ccal}
\newcommand{\intFourier}[2]{\int\limits_{#1} #2 \frac{\dd k}{(2\pi)^d}}
\newcommand{\intF}[1]{\int #1 \frac{\dd k}{(2 \pi)^d}}

\newcommand{\Riic}{\mathbb{R}_{\mathsf{ \sss IIC}}}
\newcommand{\Nbb}{N_{\sss \mathsf{Bb}}}
\newcommand{\Q}{\mathbb{Q}}
\newcommand{\dqr}{Q_{r}^{c}}
\newcommand{\tCcal}{\tilde{\Ccal}}
\newcommand{\Epc}{\mathbb{E}_{p_c}}
\newcommand{\Aconn}{\stackrel{A}{\longleftrightarrow}}
\newcommand{\ttaupc}{\widetilde{\tau}_{p_c}}

\newcommand{\taup}{\tau_p}

\newcommand{\htau}{\hat{\tau}}
\newcommand{\httau}{\widehat{\widetilde{\tau}}}

\newcommand{\htaupc}{\hat \tau_{p_c}}

\newcommand{\e}{\mathrm{e}}
\newcommand{\Ck}{\hat{C}(k)}
\newcommand{\Cmu}{\hat{C}(k-\vec{u})}
\newcommand{\Cpu}{\hat{C}(k+\vec{u})}
\newcommand{\Cmv}{\hat{C}(k-\vec{v})}
\newcommand{\Cpv}{\hat{C}(k+\vec{v})}
\newcommand{\Cml}{\hat{C}(k-q)}
\newcommand{\Cpl}{\hat{C}(k+q)}

\newcommand{\Anij}{\Rcal^{\sss (n)}_{\sss (i,j)}(\vec{u}, \vec{v})}
\newcommand{\Anmij}{\Rcal^{\sss (n,m)}_{\sss (i,j)}(\vec{u}, \vec{v})}
\newcommand{\Gio}{\Gcal^{\iota}}
\newcommand{\Fiox}{F^{\iota} (\vec{z},\vec{t},\vec{s},\vec{l},x)}
\newcommand{\Fio}{F^{\iota}}
\newcommand{\hFio}{\widehat{\Fio}}
\newcommand{\uDc}{\underline{\Dcal}}
\newcommand{\bDc}{\overline{\Dcal}}

\newcommand{\qsum}[2]{\sum_{x \in Q_{#1}} \sum_{y \in Q_{#2}}}
\newcommand{\tqr}{\Lambda_r}
\newcommand{\dd}{\operatorname{d}^d \hspace{-2pt}}
\newcommand{\dnd}{\operatorname{d} \hspace{-2pt}}

\newcommand{\Bb}{\mathsf{Bb}}
\newcommand{\Bbp}{\mathsf{Bb}_{\mathsf{piv}}}

\newcommand{\Sp}{\mathsf{S}_{\mathsf{piv}}}

\newcommand{\smi}{\mathsf{S}_m^{\sss \infty}}
\newcommand{\smr}{\mathsf{S}_m^{\sss (r)}}

\newcommand{\smx}{\mathsf{S}_m^x}

\newcommand{\sni}{\mathsf{Z}_m^{\sss \infty}}
\newcommand{\snr}{\mathsf{Z}_m^{\sss (r)}}
\newcommand{\snR}{\mathsf{Z}_m^{\sss (R)}}
\newcommand{\snx}{\mathsf{Z}_m^x}

\newcommand{\iic}{\mathsf{IIC}}

\newcommand{\liR}{\Lambda^{\sss \infty}_{\sss (R)}}

\newcommand{\lrR}{\Lambda^{\sss (r)}_{\sss (R)}}
\newcommand{\lrx}{\Lambda^{\sss (r)}_{\sss x}}
\newcommand{\bliR}{(\Lambda^{\sss \infty}_{\sss (R)})^c}
\newcommand{\blir}{(\Lambda^{\sss \infty}_{\sss (r)})^c}
\newcommand{\blrR}{(\Lambda^{\sss (r)}_{\sss (R)})^c}

\newcommand{\Picl}{\Pi^{{\rm \sss classical}}}

\begin{document}

\title{High-dimensional incipient infinite clusters revisited}
\author{Markus Heydenreich}
\address{Mathematisch Instituut, Universiteit Leiden, P.O.~Box 9512, 2300 RA Leiden, The~Netherlands;
 Centrum voor Wiskunde \& Informatica, P.O.~Box 94079, 1090~GB Amsterdam, The~Netherlands.}
\email{markus@math.leidenuniv.nl}
\author{Remco van der Hofstad}
\address{Department of Mathematics and
	Computer Science, Eindhoven University of Technology, P.O.\ Box 513,
	5600 MB Eindhoven, The Netherlands.}
\email{rhofstad@win.tue.nl, w.j.t.hulshof@tue.nl}
\author{Tim Hulshof}

\date{\today}

\begin{abstract} The incipient infinite cluster (IIC) measure is the percolation measure at criticality conditioned on the cluster of the origin to be infinite. Using the lace expansion, we construct the IIC measure for high-dimensional percolation models in three different ways, extending previous work by the second author and J\'arai. We show that each construction yields the same measure, indicating that the IIC is a robust object. Furthermore, our constructions apply to spread-out versions of both finite-range and long-range percolation models. We also obtain estimates on structural properties of the IIC, such as the volume of the intersection between the IIC and Euclidean balls.
\end{abstract}

\maketitle
{\small
\vspace{1cm}

}

\vspace{1cm}
\hrule
\vspace{1em}


\section{Introduction and main results}
It is a widely believed conjecture for bond percolation on the hypercubic lattice $\Zd$ with $d \ge 2$ that there are no infinite clusters at the critical point. This conjecture has been verified for $d=2$ \cite{Harr60},\cite{Kest80}, and $d \ge 19$  \cite{BarAiz91}, \cite{HarSla90a}. Verifying the conjecture for the intermediate values of $d$, especially the values $d=3$ to $d=6$, is arguably one of the most challenging problems in probability today.

While there are no infinite clusters at the critical point, it is nevertheless reasonable to believe that one can construct infinite clusters at the critical point through suitable conditioning and limiting schemes. These clusters are known as \emph{incipient infinite clusters} (IICs). The first construction of an IIC was carried out by Kesten \cite{Kest86a} in two dimensions. He proposed two different limiting schemes, proved their existence, and showed that both constructions have the same limit. This two-dimensional work was later extended by J\'arai \cite{Jara03b}, who proved that various other natural constructions for the IIC yield the same limiting measure.

Van der Hofstad and J\'arai \cite{HofJar04} constructed the incipient infinite cluster in high dimensions, and earlier, van der Hofstad, den Hollander and Slade constructed the IIC for high-dimensional \emph{oriented percolation} \cite{HofHolSla02}. Both constructions were achieved by making use of the \emph{lace expansion}. The lace expansion for percolation was developed by Hara and Slade \cite{HarSla90a} to treat high-dimensional percolation rigorously.

The main aim of this paper is to expand on the results by van der Hofstad and J\'arai in the following ways:
\begin{enumerate}
    \item We extend all known constructions of the IIC in high dimensions to models of \emph{long-range spread-out percolation}, so they can be dealt with under the same formalism (modulo certain assumptions).
    \item We prove a new construction of the IIC that uses the asymptotics of the \emph{one-arm probability}. Under certain assumptions we can show that this construction yields the \emph{same} limiting measure as other known constructions. This is the high-dimensional equivalent of the first construction of the IIC as given by Kesten \cite{Kest86a}. It is the first construction that has been shown to work for both two- and high-dimensional models. We also show that this measure can be obtained by taking a subsequential limit instead.
    \item We prove structural properties of the IIC, such as bounds on the volume of the intersection of the IIC with Euclidean balls centered at the origin, and the density of pivotal edges for the backbone of the IIC in such balls. These bounds will form the cornerstone for the analysis of the extrinsic properties of random walk on the IIC as performed in a sequel to this paper, \cite{HeyHofHul12b}.
    \item We introduce several new techniques for bounding probabilities and expectations in terms of the asymptotics of the two-point function in Fourier space.
    \item We prove a lower bound on the extrinsic one-arm probability for long-range spread-out percolation.
\end{enumerate}

We now start by formally introducing the model.

\subsubsection*{Bond percolation on $\Zd$.}
We consider the graph $\Zd$ as a complete graph, i.e., the set of edges (or bonds) is $\edges=\left\{\{x,y\}\mid x,y\in\Zd\right\}$.
We study bond percolation on this graph: we make the edges of the graph \emph{open} in a random way and study the resulting subgraph of open edges.
For every $x,y \in \Zd$, let the edge $\{x,y\}$ be open independently with probability $pD(x, y)$, where $D$ is such that $\sum_{y \in \Zd} D(x, y) = 1$. Thus $p$ is the average number of open edges per vertex.
In this paper, the function $D(\:\cdot\:,\:\cdot\:)$ is considered to be invariant under lattice symmetries and rotations by $90^\circ$, for instance, $D(u,v)=D(0,v-u)$. We often abbreviate $D(x)= D(0,x)$. We assume that $p \in [0, \|D\|_{\infty}^{-1}]$, so that $pD(x,y)\le 1$ for all $x,y \in \Zd$.

We consider the following three important examples.
The first example is the well-studied case of \emph{nearest-neighbor} percolation, where an edge $\{x,y\}$ is open with probability $q\in[0,1]$ whenever $|x-y|=1$, and closed otherwise. Here $|x|$ denotes the Euclidean norm of $x \in \Zd$. In terms of the above general setting, this corresponds to letting $D(x)=(2d)^{-1} \indi_{\{\lvert x \rvert =1\}}$ and $p = 2dq$.

The second example is \emph{finite-range spread-out} percolation, where the edge $\{x,y\}$ is open with probability
\begin{equation}\label{e:defFinRangePerco}
    p D(x, y) = \frac{p}{(2 L +1)^d-1}\indi_{\{0 < \|x - y\|_{\infty} \le L\}}.
\end{equation}
As is common in high-dimensional percolation, we introduce the \emph{spread-out} parameter $L$ for technical reasons. Consider it fixed at a large (integer) value for the remainder of the text.

The third example is \emph{long-range spread-out} percolation, where the edge $\{x,y\}$ is open with probability
\begin{equation}\label{e:defLongRangePerco}
    p D(x, y) =\frac{p}{\max\left\{\frac{|x-y|}L, 1\right\}^{d+\alpha}}\;.
\end{equation}

Our results hold under more general definitions for long-range and finite-range spread-out percolation: the conditions considered in \cite{HeyHofSak08} suffice.

The exponent $\alpha$ can be any positive real number, though the most interesting results are obtained for $\alpha\in(0,2]$. In such cases, the spatial variance of $D$ is infinite: $\sum_x |x|^2 D(x)= \infty$. Throughout the rest of this paper we consider $\alpha \in (0,2)\cup(2,\infty)$, that is, we consider all allowed values except $\alpha =2$. When $\alpha=2$ we get logarithmic corrections to many of the bounds, and although these do not cause complications to any of the proofs, specifying them explicitly steps would make our results cumbersome to read.

In all cases, $p$ is the parameter of the model, and it is well known that percolation undergoes a phase transition at the critical threshold
\begin{equation}
    p_c=\sup\{p\mid \chi(p)<\infty\},
\end{equation}
where
\begin{equation}
    \chi(p)=\sum_{x\in\Zd}\P_p(0\conn x)
\end{equation}
is the `expected cluster size' (or `susceptibility'), $\P_p$ is the product measure belonging to parameter $p$, and $\{x\conn y\}$ denotes the event that the vertices $x$ and $y$ are connected by a path of open edges. Note that our definition of $p_c$ differs from the standard definition
\begin{equation}
    p_c = \inf \{p: \theta(p)>0\}
\end{equation}
where $\theta(p) = \P_p(|\Ccal(0)| = \infty)$ and $\Ccal(0)$ is the connected component of the origin. However, both definitions have been proved to be equivalent in our context, cf.\ \cite{AizBar87}, \cite{Mens86}.

\subsubsection*{Mean-field behavior in high dimensions.}
Understanding percolation at the critical point $p_c$ is in general a difficult (and in many cases unsolved) problem. In the high-dimensional case some significant advances have been made.
In the context of percolation, `high-dimensional' has the rather precise meaning that the triangle diagram
\begin{equation}\label{e:TriangleDiagram}
    \triangle_{p} (0)\equiv \sum_{x,y\in\Zd}\Pp(0\conn x)\,\Pp(x\conn y)\,\Pp(y\conn 0)
\end{equation}
is finite whenever $p \le p_c$. We call this the \emph{triangle condition}.
For the nearest-neighbor model this is believed to be true whenever $d>6$, though it has been proved to hold only if $d\ge 19$ by Hara and Slade \cite{HarSla94}. For finite-range spread-out models Hara and Slade \cite{HarSla90a} were able to prove that it holds for $d>6$ if the spread-out parameter $L$ is chosen large enough.
In addition, the triangle condition is known to hold for a large class of spread-out models, and in particular for the long-range percolation model \refeq{defLongRangePerco} if $d>3 \twa$ and $L$ sufficiently large, as shown in \cite{HeyHofSak08}. Here and throughout the rest of the paper, we write $\twa$ as a shorthand for $\min\{2,\alpha\}$ when considering long-range spread-out percolation with parameter $\alpha$. Our results also hold for nearest-neighbor and finite-range spread-out percolation, with every instance of $\twa$ replaced by $2$. Instead of mentioning this every time, we simply consider $\alpha=\infty$ for these models throughout the rest of the text.
All of these proofs use the lace expansion, a method invented by Brydges and Spencer to study weakly self-avoiding walk \cite{BrySpe85} and applied to percolation by Hara and Slade \cite{HarSla90a}.

A stronger version of the triangle condition is aptly called the \emph{strong triangle condition}; it is given by
\begin{equation}\label{strongtrianglecondition}
    \triangle_{p_c} (0)  = 1+ O(\beta)
\end{equation}
with $\beta = K L^{-d}$ in the case of spread-out finite- and long-range percolation, and $\beta= K/d$ in the case of nearest-neighbor percolation. Here $K$ is a constant depending only on $d$ and $\alpha$ in the spread-out case, and is independent of $d$ in the nearest-neighbor case. It was proved in \cite{HeyHofSak08} that the strong triangle condition holds for a broad range of models, including long-range percolation. In fact, with the exceptions of \cite{Kozm10} and \cite{Scho01}, for any model for which the triangle condition has been proved, actually the strong triangle condition was obtained.

Under the triangle condition (i.e., if $\triangle_{p_c} (0) <\infty$), various critical exponents exist and take on the same value as for percolation on an infinite tree, see e.g.\ Aizenman and Newman \cite{AizNew84} and Barsky and Aizenman \cite{BarAiz91}.
Based on an analogy with the Ising model, these values are called `mean-field values'.

Here and throughout the paper, $f =o(g)$ denotes that $\lim_{n \to \infty} f(n)/g(n) = 0$ (or some other appropriate limit), $f \simeq g$ denotes that $f = c g(1+o(1))$ for some constant $c$ and $f\asymp g$ denotes that both $f\le Cg$ and $f\ge cg$ hold asymptotically for some constants $c,C >0$.
We define the \emph{two-point function}
\begin{equation}
    \P_p (x \conn y) = \tau_p (x-y).
\end{equation}
For nearest-neighbor percolation in dimension $d\ge 19$ and for finite-range spread-out percolation in dimension $d>6$, Hara \cite{Hara08} and Hara, van der Hofstad and Slade \cite{HarHofSla03}, respectively, prove the two-point function estimate
\begin{equation}\label{e:TauXStrongAsymp}
    \taupc(x-y) \simeq |x-y|^{2-d},
\end{equation}
which implies
\begin{equation}\label{e:TauXAsymp}
    \taupc(x-y) \asymp|x-y|^{2-d}.
\end{equation}

The asymptotics \refeq{TauXStrongAsymp} and \refeq{TauXAsymp} are \emph{not} true for the long-range model with $\alpha < 2$ as this would imply that $\sum_{|x| \le r} \tau(x) \asymp r^2$. However, we prove later on that $\sum_{|x| \le r} \tau(x) \asymp r^{\twa}$, so the connectivity function cannot possibly scale as $|x|^{2-d}$ when $\alpha \in (0,2)$.

\subsubsection*{The incipient infinite cluster.}
In two or high dimensions it is known that $\theta(p_c)=0$ \cite{Kest80}, \cite{BarAiz91}, \cite{HarSla90a}. However, for high-dimensional models it was proved that in a box of width $n$ around the origin, with high probability, there are several clusters whose diameter is also of order $n$ \cite{Aize97}. In other words, clusters of all sizes can be found near the origin at the critical point.

Motivated by physics literature (e.g.\ \cite{AleOrb82},\cite{LeySta83}), where random walk on critical percolation clusters were studied, Kesten \cite{Kest86a} proposed to consider the critical percolation cluster conditioned to be infinite. The object he constructed is known as Kesten's incipient infinite cluster (IIC). Since then, several different but equivalent ways of constructing IIC measures have been found, both in two and in high dimensions, cf.\ \cite{HofHolSla02}, \cite{HofJar04}, and \cite{Jara03b}.

As an example, consider the following two constructions of the IIC:
Denote by $\events_0$ the algebra of \emph{cylinder events} (i.e., events that are determined by finitely many bonds), and by $\events$ the $\sigma$-algebra of \emph{events} (i.e., the $\sigma$-algebra generated by $\events_0$).
The first construction is
\begin{equation}\label{def:Piic}
    \Piic(F)\equiv \lim_{|x| \to \infty} \P_x (F) \equiv \lim_{|x|\to\infty}\Ppc(F\mid 0\conn x),
    \qquad F\in\events_0,
\end{equation}
whenever the limit exists.
The second construction is
\begin{equation}\label{def:Qiic}
    \Qiic(F)
    \equiv \lim_{p\ua p_c}\Qp(F)
    \equiv \lim_{p\ua p_c}\frac{1}{\chi(p)}\sum_{x\in\Zd}\P_p(F\cap\{0\conn x\}),
    \qquad F\in\events_0,
\end{equation}
whenever the limit exists.
Here $\Piic$ and $\Qiic$ are understood as limits in the space of probability measures on $\{0,1\}^{\edges}$ in the weak topology. It is a priori unclear whether these limits exist, we shall elaborate on that issue later on.

We call $\Qp$ the \emph{susceptibility measure} because of the appearance of the susceptibility $\chi(p)$. It will play an important role in our analysis.

Van der Hofstad and J\'arai \cite{HofJar04} proved that, subject to \refeq{TauXAsymp}, the measures $\Piic$ and $\Qiic$ exist and are equivalent. We conjecture that this is the case in all dimensions.
However, the proof depends crucially on the aforementioned bounds on the connectivity function, restricting its use to models where such convergence is known. In the case of long-range percolation we have no useful bounds on the connectivity function $\Ppc (x \conn y)$. This means that we cannot use such a relation to bound the IIC measure for high-dimensional percolation.  The following theorem circumvents this problem by making use of the (weaker) `strong triangle condition' instead of bounds on the connectivity function.
\bt[Existence IIC measure under strong triangle condition] \label{QiicTh1}
     If the strong triangle condition (\ref{strongtrianglecondition}) holds for some $\beta$ sufficiently small, then the limit (\ref{def:Qiic}) exists for any cylinder event $F$. Consequently, $\Qiic$ can be extended to the $\sigma$-algebra of events $\sigma(\events_0) = \events$.
\et
Although this is only a minor improvement, it will turn out to be a very useful one.
In Section \ref{sec:OA} we give an outline of the adaptations that need to be made to the proof in \cite{HofJar04} to prove Theorem \ref{QiicTh1}.

We show that there exist two more constructions that both give the same IIC measure as in Theorem \ref{QiicTh1}. These constructions are based on assumptions that we make about the properties of critical percolation. These properties are not proved, but are in the spirit of some results from \cite{CheSak11} and \cite{KozNac11}.
The first assumption that we make is that the following connectivity function bounds hold for long-range percolation:
\begin{equation}\label{e:CheSakConn}
    \Ppc (x \conn y) \simeq |x - y|^{\twa -d}.
\end{equation}
A conditional proof of this relation has been circulated by Chen and Sakai \cite{CheSak12}.

To state the second assumption we need a few definitions. The vertex set $Q_r$ is defined to be the Euclidean ball of radius $r$ around the origin, that is,
\begin{equation}\label{def:QR}
    Q_r = \{x \in \Zd : |x| \le r\}.
\end{equation}

It is generally conjectured that at criticality, the probability of having a path from $0$ to $\dqr$ (the outside of a ball of radius $r$) asymptotically behaves as a power of $r$,
\begin{equation}
    \Ppc(0 \conn \dqr) \asymp r^{-1/\rho}
\end{equation}
where $\rho$ is the \emph{one-arm exponent} (cf.\ \cite[Section 9.1]{Grim99}). Based on work of Smirnov \cite{Smir01}, Lawler, Schramm, and Werner \cite{LawSchWer02b} proved that $\rho = 48/5$ for site percolation on the two-dimensional triangular lattice. They also conjectured that this is the value of the exponent for any planar lattice.

Kozma and Nachmias \cite{KozNac11} proved the following one-arm exponent for high-dimensional percolation when $\taupc (x) \asymp |x|^{2-d}$:
\begin{equation}\label{e:KozNacArm}
    \Ppc (0 \conn \dqr) \asymp r^{-2}.
\end{equation}
As mentioned before, the condition on the $x$-space asymptotics of $\taupc$ have been proved for nearest-neighbor percolation and finite-range spread-out percolation, but not for long-range spread-out percolation.
We will assume that the one-arm exponent also exists for long-range percolation, but since we do not know its value in this case, we will write $1/\rho$. Our conjecture is that for long-range percolation the correct value for $\rho$ is $2/(4 \wedge \alpha)$. Although Theorem \ref{th:OAprob} below establishes that this is a valid lower bound, we will not assume this. Instead we will use the weaker assumption that $\rho$ is well defined and $\rho \in [1/\twa,\infty)$. Furthermore, in point (iii) of the theorem below, we also assume that the asymptotics are stronger than upper and lower bound, that is, the relation is ``$\simeq$'' instead of ``$\asymp$''.
Note that we only make use of these assumptions in the statement and proof of Theorem \ref{QiicTh2}, the statement and proof of Theorems \ref{QiicTh1}, \ref{th:OAprob} and \ref{ExpectationBounds} do not require these assumptions.

\bt[Conditional IIC measure existence] \label{QiicTh2}
     Under the assumptions of Theorem \ref{QiicTh1},
     \begin{enumerate}
        \item for finite-range percolation models, and for long-range percolation models under the assumption \refeq{CheSakConn}, the limit
            \begin{equation}
                \Piic (F) \equiv \lim_{|x| \to \infty} \P_x (F) \equiv \lim_{|x| \to \infty} \Ppc (F \mid 0 \conn x)
            \end{equation}
            exists for any cylinder event $F$;
        \item for finite-range percolation models, and for long-range percolation models under the assumption that there exists $\rho \in [1/\twa,\infty)$ such that
            \begin{equation}\label{e:LRPArm}
                \Ppc (0 \conn \dqr) \asymp r^{-1/\rho},
            \end{equation}
            there exists an increasing subsequence $(r_n)$ such that along this subsequence, the limit
            \begin{equation}
                \Riic (F) \equiv \lim_{n \to \infty} \R_{r_n} (F) \equiv \lim_{n \to \infty} \Ppc (F \mid 0 \conn Q_{r_n}^{c})
            \end{equation}
            exists for any cylinder event $F$;
        \item for finite-range percolation models, and for long-range percolation models under the assumption there exists $\rho \in [1/\twa,\infty)$ such that
            \begin{equation}\label{e:LRPArmStrong}
                \Ppc(0 \conn \dqr) \simeq r^{-1/\rho},
            \end{equation}
            the limit
            \begin{equation}
                \Riic (F) \equiv \lim_{r \to \infty} \R_r (F) \equiv \lim_{r \to \infty} \Ppc (F \mid 0 \conn \dqr)
            \end{equation}
            exists for any cylinder event $F$;
        \item when the measures $\Qiic$, $\Piic$ and $\Riic$ exist, they are equal, i.e., $\Qiic = \Piic = \Riic$.
     \end{enumerate}
\et
\begin{rk} Theorem \ref{QiicTh2}(iv) is the crux of the above theorem. Theorem \ref{QiicTh2}(i) has already been proved for finite-range models by van der Hofstad and J\'arai. We repeat it here for completeness. Theorem \ref{QiicTh2}(ii) and (iii) yield versions of the IIC as in Kesten's first construction \cite{Kest86a}.
\end{rk}

In Theorem \ref{QiicTh2}(ii) and (iii) we assume that for long-range percolation the one-arm critical exponent $\rho$ exists and $\rho \in [1/\twa, \infty)$. We prove that if $\rho$ exists, then $\rho \ge 2/(4 \wedge \alpha)$:

\begin{theorem}[A lower bound on the one-arm probability for long-range percolation]\label{th:OAprob}
    When $d > d_c$, there exists $c>0$ such that for critical long-range spread-out percolation with parameter $\alpha$,
    \begin{equation}
        \Ppc(0 \conn \dqr) \ge \frac{c}{r^{(4 \wedge \alpha)/2}}.
    \end{equation}
\end{theorem}
We prove this theorem in Section \ref{sec:OAprob}. The heuristics of the proof are simple: if the cluster reaches distance $r$, then either the cluster contains many vertices, or the cluster contains an edge that is very long (of order $r$). To bound the probability that the cluster is large, we use a simple second moment estimate. This contributes the dominant term to the lower bound when $\alpha \ge 4$, as the probability of finding a long edge is negligible in this regime. When $\alpha \le 4$ however, this is not the case anymore, and the dominant contribution will be due to the existence of long edges. To establish this, we show that the existence of long edges is only weakly dependent on the size of the cluster, and vice versa.

We conjecture that $2/(4 \wedge \alpha)$ is indeed the correct value for $\rho$. Supporting evidence for this comes from Janson and Marckert's analysis of the one-dimensional discrete snake with long-range step distribution \cite{JanMar05}. Indeed, their results indicate that the maximal displacement of critical branching random walk exceeds $R$ is proportional to $R^{2/(4 \wedge \alpha)}$. Since branching random walk can be considered as mean-field model for high-dimensional percolation, one expects that the maximal displacement of branching random walk behaves similarly to the one-arm probability of high-dimensional percolation.
It would be of interest to show that $\rho = 2/(4 \wedge \alpha)$ indeed holds for percolation in high dimensions. For instance, \refeq{CheSakConn} and an adaptation of \cite{KozNac11} might be used to prove that $\rho = 1/2$ when $\alpha > 4$. However, another approach seems necessary for smaller values of $\alpha$.

\subsubsection*{Euclidean distance.}
Using properties of $\Qiic$ allows us to estimate the expected volume of the intersection between Euclidean balls and the cluster at the origin.

Let $\Eiic$ be the expectation with respect to $\Qiic$, and let $\iic = \iic(\omega)$ be the (infinite) connected component of $0$. Let $\Nbb(r)$ be the number of edges in the \emph{backbone} of $\iic$ at Euclidean distance at most $r$ from $0$, that is, all `directed' edges $b=(\ulb, \olb)$ with $\ulb \in Q_r \cap \iic$ such that $\{0 \conn \ulb\}$ and $\{\olb \conn \infty\}$ occur disjointly and $b$ is open.
\bt[Cluster and backbone volume bounds]\label{ExpectationBounds}
     Under the assumptions of Theorem \ref{QiicTh1},
        \begin{eqnarray}
            \label{e:pcVolBound} \Epc [|Q_r  \cap \Ccal(0)|] &\asymp& r^{\twa};\\
            \label{e:IICVolBound}\Eiic [|Q_r  \cap \iic|] &\asymp& r^{2\twa};\\
            \label{e:BBVolBound} \Eiic [\Nbb(r)] &\asymp& r^{\twa}.
        \end{eqnarray}
\et

Let $B_r(0; \Gcal)$ be the \emph{graph-metric ball} of radius $r$ around $0$, where the graph-metric $d_{\Gcal} (x,y)$, for all $x,y \in \Zd$, is given by the number of edges on a shortest path between $x$ and $y$ in the graph $\Gcal$. The graph-metric is also referred to as the \emph{intrinsic distance}, because it is defined by the graph structure rather than by spatial properties of the graph.
Theorem \ref{ExpectationBounds} can be contrasted with \cite[Theorems 1.3, 1.4]{KozNac09} where $\Epc[|B_r(0; \Ccal(0))|]$ is proved to be of order $r$, regardless of the range of the model (i.e., the value of $\alpha$ does not influence the asymptotics).

This paper is organized as follows:
\begin{enumerate}
    \item In Section \ref{sec:LE} we perform a lace expansion for the measure $\Riic$.
    \item In Section \ref{sec:OA} we use this lace expansion to prove Theorem \ref{QiicTh2}, subject to Proposition \ref{prop:TermBounds}. Theorem \ref{QiicTh2}(ii) and (iii) are proved in full detail, whereas we only present an indication of the proof of Theorem \ref{QiicTh2}(i). We also give an indication of the proof of Theorem \ref{QiicTh1}.
    \item In Section \ref{sec:propIIC} we prove Theorem \ref{ExpectationBounds} using Fourier space techniques. We also prove a useful lemma that establishes a way of `reversing the limit' for $\Qiic$. Both results are important ingredients in the analyses of \cite{HeyHofHul12b} and \cite{HeyHofHulMie12}.
    \item In Section \ref{sec:OAprob} we prove Theorem \ref{th:OAprob}.
    \item In Sections \ref{app:Prop} and \ref{sec:proofprop252} we prove Proposition \ref{prop:TermBounds}.
\end{enumerate}

\section{The lace expansion}\label{sec:LE}
Lace expansions for percolation have been presented in numerous papers, cf.\  \cite{BorChaHofSlaSpe05b}, \cite{HarSla90a}, \cite{Slad06}. In particular, van der Hofstad and J\'arai \cite{HofJar04} performed it with limiting schemes for the IIC in mind. Our approach is quite similar to theirs, and given that our three limiting schemes require only slightly different lace expansions, we refer the reader to the expansions for $\P_x$ and $\Q_p$ in \cite{HofJar04} and focus mainly on the lace expansion of $\R_r$. This expansion is the most involved of the three, and it actually contains almost all of the elements that are required for the expansion of the other two measures. At the end of Section \ref{sec:OA} we indicate how the other two lace expansions are done. In the sections that follow we show how the limiting behavior of the terms in the expansion can be used to show that all three measures converge to the IIC measure.

Before we start the expansion we restate an important lemma that is at the heart of every lace expansion, namely the \emph{Factorization Lemma} (Lemma \ref{lem-cut1} below).

\pagebreak
\subsection{The Factorization Lemma}
\label{sec-factlem}
Parts of this subsection are taken verbatim from \cite[Section 2]{HofHolSla07b}, where also the proof of Lemma \ref{lem-cut1} appears. We start with a few definitions.

\begin{defn}\label{def:def}
\label{def-onin2}
\begin{enumerate}
\item[(i)] For any pair $x,y \in \Zd$, we write $\{x,y\}$ to signify the \emph{undirected edge} between $x$ and $y$, and we write $(x,y)$ to signify the \emph{directed edge} from $x$ to $y$. When dealing with directed edges $b=(\ulb,\olb)$, we call $\ulb$ the `bottom' vertex, and $\olb$ the `top' vertex. We define $\Ee = \{(\ulb,\olb) : \ulb \in Q_r, \olb \in \Zd\}$, the set of directed edges with the bottom vertex inside $Q_r$ and the top vertex in $\Zd$.
\item[(ii)] Let $\omega$ be an edge configuration and $b$ an (open or closed) edge. Let $\omega^{b}$ be the same edge configuration with the status of the edge $b$ changed. We say an edge $b$ is a \emph{pivotal edge} for the configuration $\omega$ and the event $E$, if $\omega \in E$ and $\omega^{b} \,\nin\, E$, or if $\omega \,\nin\, E$ and $\omega^{b} \in E$. An edge $b$ that is pivotal for a configuration $\omega$ and a connection event $\{A \conn B\}$ will always be assumed to be directed, i.e., $b = (\ulb, \olb)$, in such a way that $\omega, \omega^{b} \in \{A \conn \ulb\}\cap\{\olb \conn B\}$. When we say that an edge is pivotal for an event this should be taken to mean that it is pivotal for that event in some fixed but unspecified configuration.
\item[(iii)] Given a (deterministic or random) set of vertices $A$ and an edge configuration $\omega$, we define $\omega_A$, the restriction of $\omega$ to $A$, to be
    \eq
    \omega_A(\{x, y\})
   =
   \left\{
   \begin{array}{lll}
   &\omega(\{x, y\})   &\text{if  }x, y \in A,\\
   &0           &\text{otherwise},
   \end{array}
   \right.
   \en
for every $x, y$ such that $\{x, y\}$ is an edge. In other words, $\omega_A$ is obtained from $\omega$ by making every edge that does not have both endpoints in $A$ closed.

\item[(iv)] Given a (deterministic or random) set of vertices $A$ and an event $E$, we say that \emph{$E$ occurs on $A$}, and write $\{E$ on $A\}$, if $\omega_A\in E$. In other words, $\{E$ on $A\}$ means that $E$ occurs on the (possibly modified) configuration in which every edge that does not have both endpoints in $A$ is made closed.  We adopt the convention that $\{x \conn x$ on $A\}$ occurs if and only if $x\in A$.

    Similarly, we say that \emph{$E$ occurs off $A$}, and write $\{E$ off $A\}$, if $\{E$ on $A^c\}$, where $A^c$ is the complement of $A$.

    We say that \emph{$E$ occurs through $A$}, and write $\{E$ through $A\}$ for the event that $E$ occurs, but $E$ does \emph{not} occur if all the edges with at least one endpoint in $A$ are made closed, that is, $\{E$ through $A\}= E \setminus \{E$ off $A\}$. For a two-point event $\{x \conn y$ through $A\}$ we write $\{x \Aconn y\}$.

\item[(v)] Given a (deterministic or random) set of vertices $A$, we define the \emph{restricted percolation measure} for any event $E$:
    \begin{equation}
        \Pp^A(E) = \Pp (E \text{ off }A).
    \end{equation}
    Given two vertices, $x$ and $y$, we define the \emph{restricted two-point function:}
    \begin{equation}
        \taup^A (x,y) = \Pp(\{x \conn y\} \text{ off }A) = \Pp^A (x \conn y).
    \end{equation}

\item[(vi)] Given an edge configuration and a set $A \subseteq \Z^d$, we define $\Ccal(A)$ to be the set of vertices to which $A$ is connected, i.e., $\Ccal(A)=\{y \in \Z^d: A \conn y\}$. Given an edge configuration and an edge $b$, we define the restricted cluster $\tilde{\Ccal}^{b}(A)$ to be the set of vertices $y \in \Ccal(A)$ to which $A$ is connected in the (possibly modified) configuration in which $b$ is made closed. When $A = \{x\}$ for some $x \in \Zd$, as will often occur, we write $\Ccal(\{x\})=\Ccal(x)$.

\end{enumerate}
\end{defn}

The statement of the Factorization Lemma is in terms of {\it two} independent percolation configurations, whose laws are indicated by subscripts $0$ and $1$.
We use the same subscripts for random variables, to indicate which law describes their distribution. Thus, the law of $\tCcal^{(u, v)}_{\sss 0}(y)$ is described by $\Pbold_{\sss 0}$, with corresponding expectation $\Ebold_{\sss 0}$.

\begin{lemma}[Factorization Lemma, \cite{HofHolSla07b}]
    \label{lem-cut1}
    Fix $p \in [0,\|D\|_\infty^{-1}]$, a directed edge $(u, v)$, a vertex $y$, and events $E, F$. Assume that $p$ is such that
    $\theta(p)=0$. Then,
    \eqalign
    \lbeq{lemcut(i)}
        & \Ebold \Big(
        \indi_{\{E \text{\emph{ on }} \tCcal^{(u, v)}(y),\;
        F \text{\emph{ off }} \tCcal^{(u, v)}(y)\}}\Big)
        =
        \Ebold_{\sss 0} \Big(\indi_{\{E \text{\emph{ on }}
        \tCcal^{(u, v)}_{\sss 0}(y)\}}
        \Ebold_{\sss 1}\big(\indi_{\{F \text{\emph{ off }} \tCcal^{(u, v)}_{\sss 0}(y)\}}\big)\Big).
     \enalign
     Moreover, when $E\subseteq \{u\,\in\,\tCcal^{(u, v)}(y), v\,\nin\, \tCcal^{(u, v)}(y)\}$, the event on the left-hand side of \refeq{lemcut(i)} is independent of the occupation status of $(u, v)$.
\end{lemma}

\subsection{The lace expansion of the one-arm IIC measure}\label{sec:LaceExpansionOA}
In this section we give the lace expansion for the measure $\Riic$ as defined in Theorem \ref{QiicTh2}.
The measure is defined for cylinder events and two-point events.
The aim is to show that for some increasing subsequence $(r_n)$, the measure
\begin{equation}\label{e:OneArmIIC}
    \Riic (F) = \lim_{n \to \infty} \R_{r_n} (F) = \lim_{n \to \infty} \Ppc (F \mid 0 \conn Q_{r_n}^c)=\lim_{r \to \infty} \frac{\Ppc (F, 0 \conn Q_{r_n}^c)}{\Ppc(0 \conn Q_{r_n}^c)}
\end{equation}
converges. We assume $F \in \mathfrak{F}_0$ to be determined by the edges in $Q_m$, for some $1 \le m \le r_1$.

Repeatedly using the inclusion-exclusion principle, we will chip away at the event $\{F, 0 \conn \dqr\}$, separating out increasingly improbable events, until we end up with a complicated but manageable expression for the right-hand side of \refeq{OneArmIIC}. In Section \ref{sec:OA} we show that the limit $\Riic$ exists and equals $\Qiic$.

When the event $\{0 \conn \dqr\}$ occurs, this implies that $\{Q_{m} \conn \dqr\}$ also occurs for any $m \le r$. Now there are two cases:
The first case is that there are no pivotal edges for $\{Q_{m} \conn \dqr\}$. This implies that both $\{0 \conn \dqr\}$ and $\{Q_{m} \Conn \dqr\}$ occur.
The second case is that there is a pivotal edge for $\{Q_{m} \conn \dqr\}$. In this case, let $(u,v)$ denote the first pivotal edge for $Q_{m} \conn \dqr$. Since $0 \conn \dqr$, the edge $(u,v)$ is also pivotal for $0 \conn \dqr$. We can therefore write
\begin{equation}\label{e:FirstStep}
    \begin{split}
        \Ppc (F, 0 \conn \dqr)  &= \Ppc(F, 0 \conn \dqr, Q_{m} \Conn \dqr) \\
                                &\quad + \sum_{(u,v)\in \Ee} \Ppc (F, 0 \conn u, Q_{m} \Conn u,\,(u,v) \text{ open and pivotal for }Q_{m} \conn \dqr)\\
                                &= \sum_{(u,v)\in \Ee} \Ppc(\{F \cap \{0 \conn u\}\cap \{Q_{m} \Conn u\}\cap\{Q_{m} \conn \dqr\}^c \text{ on }\tCcal^{(u,v)}(Q_{m})\}\\
                                &\qquad\quad \cap \{(u,v) \text{ open}\}\cap\{v \conn \dqr \text{ off }\tCcal^{(u,v)}(Q_{m})\})+ \Ppc(F, 0 \conn \dqr, Q_{m} \Conn \dqr).
    \end{split}
\end{equation}
In the last step we used a standard partition of an event involving a fixed pivotal edge into a part that occurs \emph{before} the edge (i.e., on $\tCcal^{(u,v)}(Q_{m})$) and a part occurring \emph{after} the edge (i.e., off $\tCcal^{(u,v)}(Q_{m})$). The extra event $\{Q_{m} \conn \dqr\}^c$ that occurs on $\tCcal^{(u,v)}(Q_{m})$ on the right-hand side of \refeq{FirstStep} is there to ensure that the edge $(u,v)$ is still pivotal after the partition.

We define
\begin{equation}\label{e:RhoZero}
    \xi^{(0)} (r; F) = \Ppc (F \cap\{ 0 \conn \dqr, Q_{m} \Conn \dqr\})
\end{equation}
and
\begin{equation}\label{e:GammaZero}
    \gamma^{(0)} (r;F) = \sum_{(u,v)\in \Ee} p_{uv} \E_0[\indi_{\{F \cap \{0 \conn u, Q_{m} \Conn u, Q_{m} \conn \dqr\} \text{ on }\tCcal^{(u,v)}(Q_{m}) \}} \P_{1}^{\tCcal^{(u,v)}(Q_{m})} (v \conn \dqr)]
\end{equation}
where
\begin{equation}
    p_{uv} = p_c D(v-u).
\end{equation}

Using that for any event $E$, $\indi_{E^c} = 1 - \indi_{E}$, and applying this and the Factorization Lemma to the right-hand side of \refeq{FirstStep} yields
\begin{equation}
    \Ppc(F, 0 \conn \dqr) = \xi^{(0)}(r; F) - \gamma^{(0)}(r;F) + \sum_{(u,v)\in \Ee} p_{uv} \Epc [ \indi_{F \cap\{ 0 \conn u, Q_{m} \Conn u\}} \Ppc^{\tCcal^{(u,v)}(Q_{m})}(v \conn \dqr)].
\end{equation}

For $x \in \Zd$, define
\begin{equation}
    \label{e:PiZero}  \pi^{(0)} (x,r ; F) = \Ppc(F \cap\{ 0 \conn x, Q_{m} \Conn x\}).
\end{equation}
Although $\pi^{(0)}(x,r;F)$ is independent of $r$, the higher order terms $\pi^{(n)}(x,r;F)$ do depend on $r$, so we write the redundant argument $r$ here for compatibility later on.
For $v \in \Zd$ define
\begin{eqnarray}
    \label{e:PsiZero} \psi^{(0)} (v,r; F) &= &\sum_{u \in Q_r} p_{uv} \pi^{(0)} (u,r; F);\\
    \label{e:Rzero}   R^{(0)} (r;F)& =& \sum_{(u,v)\in \Ee} p_{uv} \Epc[\indi_{F \cap\{ 0 \conn u, Q_{m} \Conn u\}}(\Ppc(v \conn \dqr) - \Ppc^{\tCcal^{(u,v)}(Q_{m})}(v \conn \dqr))].
\end{eqnarray}
Then,
\begin{equation}
    \sum_{(u,v)\in \Ee} p_{uv} \Ppc(F \cap\{ 0 \conn u, Q_{m} \Conn u\}) \Ppc (v \conn \dqr) = \sum_{v \in \Zd} \psi^{(0)} (v,r;F) \Ppc(v \conn \dqr).
\end{equation}
and using the identity
\begin{equation}
    \Ppc^{A}(v \conn \dqr)  = \Ppc(v \conn \dqr)-[\Ppc(v \conn \dqr) - \Ppc^{A}(v \conn \dqr)]
\end{equation}
with $A = \tCcal^{(u,v)}(Q_{m})$, we can write
\begin{equation}
    \Ppc(F, 0\conn \dqr) = \xi^{(0)} (r;F) - \gamma^{(0)}(r;F) + \sum_{v \in \Zd} \psi^{(0)} (v,r;F)\Ppc(v \conn \dqr) - R^{(0)}(r;F).
\end{equation}
The aim of the following derivation is to expand $R^{(0)}(r;F)$ further.

For $A \subseteq \Zd$, define the events
\begin{equation}\label{e:Eprime}
        E'(v,x;A)   = \left\{    \begin{array}{l} v \Aconn x \text{ and there is no pivotal edge }(u_1, v_1)\\
                                    \text{ for the connection }v \conn x \text{ such that } v \Aconn u_1
                                \end{array}\right\},
\end{equation}
and
\begin{equation}\label{e:Edprime}
        E''(v,r;A)   = \left\{    \begin{array}{l} v \Aconn \dqr \text{ and there is no pivotal edge }(u_1, v_1)\\
                                    \text{ for the connection }v \conn \dqr \text{ such that } v \Aconn u_1
                                \end{array}\right\}.
\end{equation}

Let $A \dot{\cup}B$ denote the \emph{disjoint union} of $A$ and $B$: the union of the events $A$ and $B$ that have no elements in common (i.e., $A \cap B = \varnothing$).
We can partition the event $\{v \Aconn \dqr\}$ into a disjoint union of events $E'$ and $E''$:
\bl\label{lem:partition}
     For any $v \in \Zd$, $A \subseteq \Zd$ and $r \in \N$:
    \begin{equation}
        \{v \Aconn \dqr\} = E''(v,r;A) \dot{\cup}\, \dot{\bigcup_{b=(\ulb,\olb)\in \Ee}} [E'(v,\ulb;A) \cap \{b \text{ is open and pivotal for }v \conn \dqr\}].
    \end{equation}
\el
\proof We decompose the event $\{v \Aconn \dqr\}$ according to whether or not there is an open pivotal edge $b=(\ulb,\olb)$ such that $\{v \Aconn \ulb\}$ and $b \in \Ee$ is the first such edge along the path from $v$ to $\dqr$ that has this property. When such an edge does not exist, the event $E''(v,r;A)$ occurs. If an edge $b$ with these properties does exist, then, since it is the first edge that is pivotal for $\{v \conn \dqr\}$ and $\{v \Aconn \ulb\}$ occurs, there can be no other edge $b'\in \Ee$ that is open and pivotal for $\{v \conn \ulb\}$ such that $\{v \Aconn \ulb'\}$. Therefore, $E'(v,\ulb;A)$ holds. \qed

By the Factorization Lemma, for any $v \in Q_r$, $r\in \N$, $b=(\ulb,\olb) \in \Ee$ and $A \subseteq \Zd$,
    \begin{equation}\label{e:Edprimepart}
        \Epc[\indi_{\{E'(v,\ulb;A)\cap\{b \text{ open \& piv. for }v \conn \dqr\}\}}] = p_b \Epc[\indi_{\{E'(v,\ulb;A) \cap \{v \conn \dqr\}^c \text{ on }\tCcal^b (v)\}} \Ppc^{\tCcal^{b}(v)}(\olb \conn \dqr)]
    \end{equation}
where $p_b = p_{\ulb\olb}$.

Using $\{v \Aconn \dqr\} = \{v \conn \dqr\} \setminus \{v \conn \dqr \text{ off } A\}$, Definition \ref{def:def} and Lemma \ref{lem:partition} and \refeq{Edprimepart}, we can write
\begin{equation}\label{e:Iteration} \begin{split}
    \Ppc(v \Aconn \dqr) =& \Ppc(v \conn \dqr) - \Ppc^{A} (v \conn \dqr)\\
        =& \Ppc(E'' (v,x;A)) + \sum_{b\in \Ee} p_{b} \Epc[ \indi_{E'(v,\ulb;A)}\, \Ppc^{\tCcal^{b}(v)}(\olb  \conn \dqr)]\\
        &-\sum_{b\in \Ee} p_{b} \Epc[ \indi_{\{E'(v,\ulb;A) \cap \{v \conn \dqr\} \text{ on }\tCcal^{b}(v)\}}\, \Ppc^{\tCcal^{b}(v)}(\olb \conn \dqr)].
    \end{split}
\end{equation}

Given edges $(u_0, v_0)$, $(u_1, v_1), \dotso$, we denote
\begin{equation}
    \tCcal_0 = \tCcal^{(u_0, v_0)}(Q_{m}) \quad\text{ and }\quad \tCcal_j = \tCcal^{(u_j, v_j)}(v_{j-1}) \quad\text{ for } j\ge1
\end{equation}
and write
\begin{equation}
    \indi_j = \indi_{E'(v_{j-1}, u_j ; \tCcal_{j-1})}, \quad \text{ for }j\ge1.
\end{equation}
Inserting \refeq{Iteration} with $v=0$ and $A=\tCcal_0$ into \refeq{Rzero}, yields
\begin{equation}\label{e:RzeroB}
    \begin{split}
        R^{(0)} (r;F) = & \sum_{(u_0,v_0)\in \Ee} p_{u_0 v_0} \E_0[\indi_{F \cap \{ 0\conn u_0, Q_{m} \Conn u_0\}} \E_1 [\indi_{E'' (v_0, r; \tCcal_0)}]]\\
                        & + \sum_{(u_0,v_0)\in \Be} p_{u_0 v_0} \sum_{(u_1, v_1)\in \Ee} p_{u_1 v_1} \E_0 [\indi_{F \cap \{ 0 \conn u_0, Q_{m} \Conn u_0\}} \E_1[\indi_1 \Ppc^{\tCcal_1} (v_1 \conn \dqr)]]\\
                        &-\sum_{(u_0,v_0)\in \Be} p_{u_0 v_0} \sum_{(u_1, v_1)\in \Ee} p_{u_1 v_1} \E_0 [\indi_{F \cap \{ 0 \conn u_0, Q_{m} \Conn u_0\}}\\
                        & \qquad \qquad \qquad\qquad \quad \times \E_1[\indi_{\{E'(v_0, u_1; \tCcal_0)\cap\{v_0 \conn \dqr\} \text{ on }\tCcal_1\}} \Ppc^{\tCcal_1} (v_1 \conn \dqr)]].
    \end{split}
\end{equation}
We define the first term on the right-hand side as $\xi^{(1)}(r;F)$ and the last term as $\gamma^{(1)}(r;F)$. We define
\begin{equation}
    \psi^{(1)}(v,r;F)= \sum_{u \in Q_r} p_{u v} \sum_{(u_0 v_0)\in \Be} p_{u_0 v_0} \E_0 \left[\indi_{F \cap \{0 \conn u, Q_m \Conn u_0\}} \E_1\left[\indi_{E'(v_0, u; \tCcal_0)}\right]\right]
\end{equation}
and we define $R^{(1)}(r;F)$ such that
\begin{multline}
    \sum_{(u_0,v_0)\in \Be} p_{u_0 v_0} \sum_{(u_1 v_1)\in \Ee} p_{u_1 v_1} \E_0 [\indi_{F \cap \{ 0 \conn u_0, Q_{m} \Conn u_0\}} \E_1[\indi_1 \Ppc^{\tCcal_1} (v_1 \conn \dqr)]] \\
        = \sum_{v \in Q_r} \psi^{(1)}(v,r;F)\Ppc(v \conn \dqr) - R^{(1)}(r;F)
\end{multline}
where we used that $\Ppc^{\tCcal_1}(v \conn \dqr) =\Ppc(v \conn \dqr)-\Ppc(v \stackrel{\tCcal_1}{\longleftrightarrow} \dqr)$.

Hence, we can write $R^{(0)}(r;F)$ as
\begin{equation}
    R^{(0)}(r;F) = \xi^{(1)}(r;F) - \gamma^{(1)}(r;F) + \sum_{v \in \Zd} \psi^{(1)}(v,r;F) \Ppc(v \conn \dqr) - R^{(1)}(r;F).
\end{equation}
From here we continue to extract terms $\xi^{(2)}(r;F)$, $\gamma^{(2)}(r;F)$, $\sum_{v \in Q_r} \psi^{(2)}(v,r;F) \Ppc(v \conn \dqr)$ and $R^{(2)}(r;F)$ from $R^{(1)}(r;F)$, and so forth. We end up with the following:

\begin{prop}[The lace expansion]
For $N\ge0$ and $p \le p_c$,
\begin{equation}\label{e:Expansion1}
    \begin{split}
        \Ppc(F , 0 \conn \dqr)  &= \sum_{n=0}^N (-1)^n \xi^{(n)} (r;F) - \sum_{n=0}^N (-1)^n \gamma^{(n)}(r;F)\\
                                &\quad +\sum_{n=0}^N (-1)^n \sum_{v \in \Zd} \psi^{(n)} (v,r;F) \Ppc(v \conn \dqr) + (-1)^{N+1} R^{(N)}(r;F).
    \end{split}
\end{equation}
Here, $\xi^{(0)}(r;F)$ is given by \refeq{RhoZero}, $\gamma^{(0)}(r;F)$ is given by \refeq{GammaZero}, $\pi^{(0)}(x,r;F)$ is given by \refeq{PiZero}, and for $n \ge1$,
\begin{eqnarray}
    \label{e:Expansion2} \xi^{(n)}(r;F) &=& \sum_{(u_0, v_0)\in \Be} p_{u_0 v_0} \dotsm \sum_{(u_{n-1}, v_{n-1})\in \Ee} p_{u_{n-1} v_{n-1}} \E_0 \bigg[\indi_{F \cap \{ 0 \conn u_0, Q_{m} \Conn u_0\}}\\
    \nonumber           & &\quad \times \E_1\Big[\indi_1 \E_2\big[ \indi_2 \dotsm \E_{n-1}[\indi_{n-1} \E_n[ \indi_{E''(v_{n-1}, r;\tCcal_{n-1})}]]\dotsm \big]\Big]\bigg];\\
   \label{e:Expansion6} \gamma^{(n)}(r;F) &=& \sum_{(u_0, v_0)\in \Be} p_{u_0 v_0} \dotsm \sum_{(u_{n}, v_{n})\in \Ee} p_{u_{n} v_{n}}\E_0 \bigg[\indi_{F \cap \{ 0 \conn u_0, Q_{m} \Conn u_0\}}\\
    \nonumber           & &\quad \times \E_1\Big[\indi_1 \E_2\big[ \indi_2 \dotsm \E_{n-1}[\indi_{\{E'(v_{n-1}, u_n; \tCcal_{n-1})\cap\{v_{n-1} \conn \dqr\} \text{\emph{ on }}\tCcal_n\}} \P_{n}^{\tCcal_n} (v_n \conn \dqr)]\dotsm \big]\Big]\bigg];\\
   \label{e:Expansion3} \pi^{(n)}(x,r;F) &=& \sum_{(u_0 v_0)\in\Be} p_{u_0, v_0} \dotsm \sum_{(u_{n-1}, v_{n-1})\in \Be} p_{u_{n-1} v_{n-1}} \E_0 \bigg[\indi_{F \cap \{ 0 \conn u_0, Q_{m} \Conn u_0\}}\\
     \nonumber          & &\quad \times \E_1\Big[\indi_1 \E_2 \big[ \indi_2 \dotsm \E_{n-1}[\indi_{n-1} \E_n[ \indi_{E'(v_{n-1}, x;\tCcal_{n-1})}]]\dotsm\big]\Big]\bigg].
\end{eqnarray}
Also, for $n \ge 0$,
\begin{equation}\label{e:Expansion4}
    \psi^{(n)}(v,r;F) = \sum_{u \in Q_r} p_{uv} \pi^{(n)} (u,r;F),
\end{equation}
and
\begin{equation}\label{e:Expansion5}
    \begin{split}
        R^{(N)}(r;F)    &= \sum_{(u_0, v_0)\in \Be} p_{u_0 v_0} \dotsm \sum_{(u_N, v_N)\in \Ee} p_{u_N v_N} \E_0 \bigg[ \indi_{F \cap \{ 0\conn u_0, Q_{m} \Conn u_0\}}\\
                        &\quad \times \E_1 \Big[\indi_1 \E_2 \big[\indi_2 \dotsm \E_N [\indi_N (\Ppc(v_N \conn \dqr) - \Ppc^{\tCcal_N}(v_N \conn \dqr))]\dotsm \big]\Big]\bigg].
    \end{split}
\end{equation}
\end{prop}

\subsection{Bounds on the expansion terms}
In order to prove Theorem \ref{QiicTh2}, we have to give bounds on the terms of the expansion. We do this using the following proposition:
\begin{prop}[Bounds on expansion terms]\label{prop:TermBounds}
Under the same assumptions as Theorems \ref{QiicTh1} and \ref{QiicTh2}(ii), the following holds:
\be
    \item   For $L \ge L_0$, any $r\in \N$ and any $F \in \events_0$, there is a constant $K=K(F,L,d,\alpha)$ and an $\vep >0$ such that
            \begin{equation}\label{e:XiBound}
                \sum_{n=0}^\infty \xi^{(n)} (r;F) \le \frac{K}{r^{1/\rho+ \vep}} \qquad \text{ and }\qquad \sum_{n=0}^\infty \gamma^{(n)} (r;F) \le \frac{K}{r^{1/\rho +\vep}};
            \end{equation}
    \item   For some $\delta > 0$, any $r \in \N$ and any $L \ge L_0$ and any $F \in \events_0$ there is a $K'=K'(F,L,d,\alpha, \delta)$ such that
            \begin{equation}
                \sum_{x \in \Zd} \sum_{n=0}^\infty \lvert x \rvert^{\twa +\delta} \pi^{(n)}(x,r;F) \le K';
            \end{equation}
    \item   For any $r \in \N$
            \begin{equation}
                \lim_{N \to \infty} R^{(N)} (r;F) = 0.
            \end{equation}
\ee
\end{prop}
We prove this proposition in Sections \ref{app:Prop} and \ref{sec:proofprop252}.

\begin{rk}
By setting $r=\infty$ and $F=\Omega$ in \refeq{Expansion3}, we can define
\begin{equation}
    \Picl(x) \equiv \sum_{n=0}^{\infty} (-1)^n \pi^{(n)}(x, \infty; \Omega).
\end{equation}
It is a well known result (see e.g.\ \cite[Chapter 10]{Slad06}) that the `classical' inclusion-exclusion lace expansion for the percolation two-point function  yields the convolution equation
\begin{equation}
    \taupc(x) = \Picl (x) + (p_c D * \taupc * \Picl)(x).
\end{equation}
With minor modifications to the proof of Proposition \ref{prop:TermBounds}(ii) one can show that for some $\delta'>0$ and $L \ge L_0$, there is a $K''= K''(L,d,\alpha,\delta')$ such that
\begin{equation}
    \sum_{x \in \Zd} \lvert x \rvert^{\twa+\delta'} \pi^{(n)}(x,\infty;\Omega) \le K''.
\end{equation}
\end{rk}

\section{Existence of the IIC in various constructions}\label{sec:OA}

\subsection{Existence of the one-arm IIC measure}
\noindent
\emph{Proof of Theorem \ref{QiicTh2}(ii) and (iii) subject to Proposition \ref{prop:TermBounds}.}
Define
\begin{eqnarray}
    \label{e:Xidef}\Xi(r;F) &=& \sum_{n=0}^\infty (-1)^n \xi^{(n)}(r;F);\\
    \label{e:Gammadef}\Gamma(r;F) &=& \sum_{n=0}^{\infty} (-1)^n \gamma^{(n)}(r;F);\\
    \label{e:Pidef}\Pi(x,r;F) &=& \sum_{n=0}^\infty (-1)^n \pi^{(n)}(x,r;F);\\
    \label{e:Psidef}\Psi(x,r;F)&=& \sum_{n=0}^\infty (-1)^n \psi^{(n)}(x,r;F).
\end{eqnarray}
By Proposition \ref{prop:TermBounds} and \refeq{Expansion4} these sums converge. Therefore we may take the limit $N \to \infty$ in \refeq{Expansion1} to obtain
\begin{equation}\label{e:ExpansionAfterLim}
    \Ppc(F, 0 \conn \dqr) = \Xi(r;F) -\Gamma(r;F) + \sum_{y \in \Zd} \Psi(y,r;F)\Ppc(y \conn \dqr).
\end{equation}
Dividing \refeq{ExpansionAfterLim} by $\Ppc(0 \conn \dqr)$ gives
\begin{equation}\label{e:Riicprelimit}
    \R_r (F) = \frac{\Xi (r;F) - \Gamma(r;F)}{\Ppc(0 \conn \dqr)} + \sum_{y \in \Zd} \Psi(y,r;F) \frac{\Ppc(y \conn \dqr)}{\Ppc(0 \conn \dqr)}.
\end{equation}
The aim is now to show that $\lim_{r \to \infty} \R_r (F) = \sum_{y \in \Zd} \Psi(y ;F)$, and that $\Psi(y;F) = \lim_{r \to \infty} \Psi(y,r;F)$ exists.

By \refeq{LRPArm} and Proposition \ref{prop:TermBounds}(i),
\begin{equation}
    \lim_{r \to \infty} \frac{\Xi(r;F) - \Gamma(r;F)}{\Ppc(0 \conn \dqr)} =0.
\end{equation}
We are left to deal with the second term of \refeq{ExpansionAfterLim}. By \refeq{Expansion4},
\begin{equation}
        \sum_{y \in \Zd} \Psi(y,r;F) \frac{\Ppc(y \conn \dqr)}{\Ppc(0 \conn \dqr)} \le \sum_{y\in \Zd} \sum_{x \in Q_r} p_{xy} \lvert \Pi (x,r;F) \rvert \frac{\Ppc(y \conn \dqr)}{\Ppc (0 \conn \dqr)}.
\end{equation}
Note that for any edge $(x,y) \in \Ee$,
\begin{equation}
    p_{xy} \Ppc(y \conn \dqr) \le \Ppc(x \conn \dqr),
\end{equation}
and also note that $\sum_{y\in \Zd} D(y-x)=1$. Therefore,
\begin{equation}
    \sum_{y\in \Zd} \sum_{x \in Q_r} p_{xy} \lvert \Pi (x,r;F) \rvert \frac{\Ppc(y \conn \dqr)}{\Ppc (0 \conn \dqr)} \le p_c \sum_{x \in Q_r} \lvert \Pi (x,r;F) \rvert \frac{\Ppc(x \conn \dqr)}{\Ppc (0 \conn \dqr)}.
\end{equation}
To evaluate the right-hand side, we split up the sum over $x \in Q_r$ into two parts. For $a\in (0,1)$ we evaluate separately the contributions to the sum from $|x| \le r^a$ and $|x|>r^a$. We start with the latter. We first prove
\begin{equation}\label{e:Pioutside}
    \sum_{\substack{ x \in Q_r \\ |x| > r^a}} \lvert \Pi (x,r;F)\rvert \frac{\Ppc(x \conn \dqr)}{\Ppc(0 \conn \dqr)} = o(1),
\end{equation}
so the dominant contributions to the sum arise from $x \in Q_{r^a}$. Splitting once more gives
\begin{equation}\label{e:convbigr}
    \sum_{\substack{ x \in Q_r \\ |x| > r^a}} \lvert \Pi (x,r;F)\rvert \frac{\Ppc(x \conn \dqr)}{\Ppc(0 \conn \dqr)} \le \sum_{r^a < |x| < \frac{r}{4}}  \lvert \Pi (x,r;F)\rvert \frac{\Ppc(x \conn \dqr)}{\Ppc(0 \conn \dqr)} + \sum_{\frac{r}{4} \le |x| \le r} \frac{\lvert \Pi (x,r;F)\rvert}{\Ppc(0 \conn \dqr)}.
    \end{equation}
Bounding the first term on the right-hand side, we use that for all $x$ with $|x| < r/4$,
\begin{equation}
    \Ppc(x \conn \dqr) \le \Ppc(0 \conn Q_{r/2}^{c}) \le C r^{-1/\rho},
\end{equation}
so from \refeq{LRPArm} it follows that
\begin{equation}
    \frac{\Ppc(x \conn \dqr)}{\Ppc(0 \conn \dqr)} \le C.
\end{equation}
Therefore,
\begin{equation}
    \sum_{r^a < |x| < \frac{r}{4}}  \lvert \Pi (x,r;F)\rvert \frac{\Ppc(x \conn \dqr)}{\Ppc(0 \conn \dqr)} \le C \sum_{r^a < |x| < \frac{r}{4}}  \lvert \Pi (x,r;F)\rvert.
\end{equation}

For all $x$ such that $r^a < \lvert x \rvert$, we have $\lvert x \rvert / r^a >1$, so by Proposition \ref{prop:TermBounds}(ii)
\begin{equation}
     C \sum_{r^a < |x| < \frac{r}{4}}  \lvert \Pi(x,r;F)\rvert \le \frac{C}{r^{a(\twa + \delta)}}  \sum_{r^a < |x| < \frac{r}{4}} \lvert x \rvert^{\twa + \delta} \lvert \Pi (x,r;F)\rvert \le \frac{C}{r^{a(\twa + \delta)}} = o(1).
\end{equation}
Hence, the first term on the right-hand side of \refeq{convbigr} is $o(1)$.

To bound the second term on the right-hand side of \refeq{convbigr} we also use Proposition \ref{prop:TermBounds}(ii): when $\lvert x \rvert > r/4$, we have by \refeq{LRPArm} that $\Ppc(0 \conn \dqr) \ge c r^{-1/\rho} \ge c |4 x|^{-\twa}$. Furthermore, $(4 |x|)^\delta / r^\delta \ge 1$, so it follows that
\begin{equation}
    \sum_{\frac{r}{4} \le |x| \le r} \frac{\lvert\Pi (x,r;F)\rvert}{\Ppc(0 \conn \dqr)} \le \sum_{\frac{r}{4} \le |x| \le r} C |x|^{\twa} \lvert\Pi (x,r;F)\rvert
        \le \frac{4^\delta C}{r^\delta} \sum_{x \in \Zd} |x|^{\twa+\delta} \lvert\Pi (x,r;F)\rvert = O(r^{-\delta}) = o(1).
\end{equation}
This proves \refeq{Pioutside}.

In Theorem \ref{QiicTh2}(iii) we assumed that $\Ppc(0 \conn \dqr) \simeq r^{-1/\rho}$, which implies by monotonicity of the one-arm probability in $r$, that the ratio of the one-arm probabilities converges to 1 whenever $|x|$ is sufficiently small, i.e.,
\begin{equation}
    \lim_{r \to \infty} \frac{\Ppc(x \conn \dqr)}{\Ppc(0 \conn \dqr)} \le \lim_{r \to \infty} \frac{r^{1/\rho}}{(r-r^a)^{1/\rho}}(1+o(1))  =1.
\end{equation}
Furthermore, $\pi^{(n)}(x,r;F)$ is monotonically increasing as $r$ increases.
Hence, taking the limit $r \to \infty$ in \refeq{Riicprelimit}, it follows that
\begin{equation}\label{e:FinalLimit}
    \Riic (F) = \lim_{r \to \infty} \R_r (F) = \sum_{x \in \Zd} \Psi(x; F) = p_c \sum_{x \in \Zd} \Pi (x;F),
\end{equation}
exists by monotone convergence. Here $\Pi(x;F)$ is the function $\Pi(x,r;F)$ with all the summations over edges extended to the set $\Zd \times \Zd$, and a similar definition for $\Psi(x;F)$. The last step follows from $\sum_v D(v-u)=1$, and we are done. Note that the right-hand side of \refeq{FinalLimit} is the exact same expression for any of the other known limiting schemes for construction of the IIC, so $\Riic$ is in fact the same measure as $\Piic$ and $\Qiic$.

\medskip
To prove Theorem \ref{QiicTh2}(ii), a more involved analysis of the limit ratio is required.
The important contributions come from the vertices near the origin, i.e.,
$|x| \le r^a$. We show that for such $x$, the ratio of the probabilities converges to 1 along some subsequence of $(r)$.
\bl[Convergence of the ratio]\label{lem:RatioConvergence}
    Under the assumptions of Theorem \ref{QiicTh2}(ii) there exists a sequence $(r_n)$ with $r_n \to \infty$ as $n \to \infty$, such that for any $a \in (0,1)$
    \begin{equation}\label{e:RatioLimit}
        \limsup_{n \to \infty} \max_{\lvert x \rvert \le r^a} \left\lvert \frac{\Ppc(x \conn Q_{r_n}^c)}{\Ppc(0 \conn Q_{r_n}^c)} -1 \right\rvert =0.
    \end{equation}
\el
\proof Let $\Acal$ be the set of accumulation points of
\begin{equation}
    \{r^{1/\rho} \Ppc (0 \conn \dqr) | r \in \R\}.
\end{equation}
The set $\Acal$ is closed. In the finite-range setting, Kozma and Nachmias proved that $\Acal$ is bounded and positive \cite{KozNac11}, and for long-range percolation this is our assumption. Hence, it is compact and contains a positive minimum:
\begin{equation}
    A = \min \Acal \in (0,\infty).
\end{equation}
Since $A$ is an accumulation point, there exists a subsequence $(\tilde{r}_{n})_{n \in \N}$ such that
\begin{equation}
    \lim_{n \to \infty} \tilde{r}_{n}^{1/\rho} \Ppc (0 \conn Q_{\rho_{n}}^{c}) = A.
\end{equation}
Choose the sequence $(r_n)$ such that $r_n - r_{n}^a = \tilde{r}_{n}$.

Take $N \in \N$ such that $|x| \le r_{n}^a$ for all $n \ge N$.
By translation invariance of the measure $\Ppc$ and monotonicity of the event $\{0 \conn \dqr\}$ as $r$ increases, we have for $x \in Q_{r^a}$ the bounds
\begin{equation}
    \Ppc(0 \conn  Q_{r+r^a}^{c}) \le \Ppc(x \conn \dqr) \le \Ppc(0 \conn Q_{r-r^a}^{c}).
\end{equation}
This implies
\begin{equation}
    \max_{\lvert x \rvert \le r_n^a} \left\lvert \frac{\Ppc(x \conn Q_{r_n}^c)}{\Ppc(0 \conn Q_{r_n}^c)} -1 \right\rvert \le \max \left(1- \frac{\Ppc(0 \conn Q_{r_n+r_n^a}^{c})}{\Ppc(0 \conn Q_{r_n}^c)}, \frac{\Ppc(0 \conn Q_{r_n -r_n^a}^c)}{\Ppc(0 \conn Q_{r_n}^c)} -1\right)
\end{equation}
so that it suffices to show that there exists a (single) subsequence $(r_n)_{n \in \N}$ for which
\begin{equation}\label{e:Pratio}
    \lim_{n \to \infty} \frac{\Ppc(0 \conn Q_{r_n+r_n^a}^{c})}{\Ppc(0 \conn Q_{r_n}^c)} = 1 \quad\text{ and }\quad\lim_{n \to \infty} \frac{\Ppc(0 \conn Q_{r_n}^c)}{\Ppc(0 \conn Q_{r_n -r_n^a}^c)} =1.
\end{equation}
Since $\Ppc(0 \conn Q_r^c)$ is monotonically decreasing in $r$, both ratios are at most equal to $1$. Monotonicity also implies that $\Ppc(0 \conn Q_{r_n + r_{n}^a}^c) \le \Ppc(0 \conn Q_{r_n - r_n^a}^c)$, so \refeq{Pratio} follows once we show that
\begin{equation}
    \liminf_{n\to\infty} \frac{\Ppc(0 \conn Q_{r_n+r_n^a}^{c})}{\Ppc(0 \conn Q_{r_n -r_n^a}^c)} =1.
\end{equation}
It is obvious that the left-hand side is at most $1$, so we will focus on proving that it also is at least $1$.

We give a proof by contradiction. Suppose that there exists an $0< \vep <1$ such that, for the sequence $(r_n)_{n \in \N}$,
\begin{equation}
    \Ppc(0 \conn Q_{r_n + r_{n}^a}^{c} \mid 0 \conn Q_{r_n - r_{n}^a}^{c}) \le 1- \vep.
\end{equation}
Then also
\begin{equation}
    \frac{(r_n -r_{n}^a)^{1/\rho}}{(r_n + r_{n}^a)^{1/\rho}} \frac{(r_n + r_{n}^a)^{1/\rho}\, \Ppc(0 \conn Q_{r_n + r_{n}^a}^{c})}{(r_n -r_{n}^a)^{1/\rho} \, \Ppc(0 \conn Q_{r_n - r_{n}^a}^{c})} \le 1-\vep.
\end{equation}
There exists an $N' \in \N$, such that for $n \ge N'$,
\begin{equation}
     \frac{(r_{n}-r_{n}^a)^{1/\rho}}{(r_n + r_{n}^a)^{1/\rho}} \ge \sqrt{1- \vep},
\end{equation}
so it follows that
\begin{equation}\label{e:inequality}
     \sqrt{1-\vep}\, (r_n+ r_{n}^a)^{1/\rho}\, \Ppc (0 \conn Q_{r_n + r_{n}^a}^{c}) \le (1-\vep)\, (r_n - r_{n}^a)^{1/\rho}\, \Ppc( 0 \conn  Q_{r_n - r_{n}^a}^{c}).
\end{equation}

Note that our choice of $r_n$ implies $r_n + r_{n}^a \le \tilde{r}_{n} + 3 \tilde{r}_{n}^a$ as soon as $r_n \ge 3^{1/(1-a)}$.
Taking $\liminf$ on both sides of \refeq{inequality} and using the fact that $A$ is the minimum of $\Acal$, we obtain
\begin{equation}
    \begin{split}
        \sqrt{1-\vep}A  & \le \sqrt{1-\vep} \liminf_{n} (\tilde{r}_{n} + 3 \tilde{r}_{n}^a)^{1/\rho} \Ppc(0 \conn  Q_{\tilde{r}_{n} + 3 \tilde{r}_{n}^a}^{c})\\
                        & \le (1-\vep) \liminf_n \tilde{r}_{n}^{1/\rho} \Ppc (0 \conn Q_{\tilde{r}_{n}}^{c}) = (1-\vep) A,
    \end{split}
\end{equation}
which yields a contradiction. This proves \refeq{Pratio} and hence the claim of the lemma follows. \qed

Applying Lemma \ref{lem:RatioConvergence}, it follows that
\begin{equation}\label{e:FinalLimitRn}
    \Riic (F) = \lim_{n \to \infty} \R_{r_n} (F) = \sum_{x \in \Zd} \Psi(x ; F) = p_c \sum_{x \in \Zd} \Pi (x;F),
\end{equation}
along a subsequence $(r_n)$, proving Theorem \ref{QiicTh2}(ii). We observe that this limit is the same as the one we obtained in \refeq{FinalLimit}, and furthermore, following the proofs outlined in the next subsection, one can easily verify that this is also the limit for $\lim_{x \to \infty} \P_x (F)$ and $\limp \Q_p (F)$, so $\Qiic = \Piic = \Riic$, and Theorem \ref{QiicTh2}(iv) follows. \qed

\subsection{Existence of the IIC susceptibility and two-point constructions}
The IIC susceptibility construction is quite similar to the one given in \cite{HofJar04}, so the proof of Theorem \ref{QiicTh1} can be given along the same lines as that of the original construction (i.e., \cite[Theorem 1.2]{HofJar04}).  The single modification of the argument is that the $x$-space bounds on the two-point function used in \cite{HofJar04} to bound the lace expansion diagrams are replaced by bounds on the triangle diagram in Fourier space to achieve the same effect. The main reasoning, however, remains unchanged. Hence, we will not give the proof.

The proof of Theorem \ref{QiicTh2}(i), in turn, is very similar to \cite[Theorem 1.1]{HofJar04}. The only modification that we need to make here is replacing every instance of nearest-neighbor and finite-range connectivity function by the long-range connectivity function, and subsequently, we must apply our assumed bound \refeq{CheSakConn} to every instance $\taupc(x)$ instead of the bound \refeq{TauXStrongAsymp} that is used in \cite{HofJar04}.

Going through the steps of the proof with this replacement, it is not hard to see that the only thing we need to show is that for $d > 3 \twa$ and some constant $C'$ that depends only on $L,d$ and $\alpha$,
\begin{equation}
    \lvert \Pi (x;F)\rvert \le \frac{C}{(\lvert x \rvert +1)^{2(d -\twa)}}
\end{equation}
to complete the proof.
That this bound indeed holds follows immediately from \refeq{CheSakConn} and \cite[Proposition 1.8(c)]{HarHofSla03}.

\section{Volume estimates of the IIC: proof of Theorem \ref{ExpectationBounds}}\label{sec:propIIC}
In this section we calculate upper and lower bounds for the expectations of the volume of critical percolation clusters and IICs inside Euclidean balls. We also calculate bounds on the expected volume of the IIC backbone.

\subsection{The BK-inequality and bounds on two triangle diagrams}
An important tool in the coming analysis is the \emph{van den Berg-Kesten inequality} (BK for short) \cite{BerKes85},\cite{Grim99}. We call an event $A$ \emph{increasing} if for any two configurations $\omega$ and $\omega'$ such that $\omega \preceq \omega'$ (that is, any edge that is open in $\omega$ is also open in $\omega'$), $\omega \in A$ implies $\omega' \in A$. Hence, by a standard coupling argument, if $A$ is increasing, then $\Pp(A) \le \P_{p'}(A)$ whenever $p < p'$. For two increasing events $A$ and $B$ we write  $A \circ B$  to indicate the disjoint occurrence of $A$ and $B$. This is the event that the set of edges can be split in two parts, say $K$ and $K^c$,  such that $A$ occurs on $K$ (i.e., $\omega_K \in A$) and $B$ occurs on $K^c$ (i.e., $\omega_{K^c} \in B$).
The BK-inequality states
\begin{equation}
    \P_p(A\circ B)\le \P_p(A)\,\P_p(B).
\end{equation}

We use the BK-inequality to bound the probability of complicated events, such as those in the various lace expansion terms, by a product of the probabilities of its constituent disjoint events.

It is often convenient to reduce events to triplets of disjointly occurring path event. Taking the probability of these triplets results in the so-called \emph{triangle diagrams}
\begin{equation}\label{e:triangle}
         \triangle_p (x)= (\tau_p * \tau_p * \tau_p)(x), \quad \tbar_p  = \sup_{x \in \Z^d} \triangle_p (x),
\end{equation}
\begin{equation}\label{e:triangleT}
        T_p (x) = (\tau_p * \tau_p * D *\tau_p)(x), \quad T_p = \sup_{x \in \Z^d} T_p (x).
\end{equation}
We can give bounds on the two triangle diagrams for the percolation models as described in the introduction, and in the more general setting of \cite{HeyHofSak08} in terms of the parameter $\beta$ (as defined below (\ref{strongtrianglecondition})):
\bl[An upper bound on the open triangle]\label{orderbounds}
    When $D(\;\cdot\;)$ satisfies the assumptions stated in \cite{HeyHofSak08} for all $p \le p_c$, then $\tbar_p\le 1 + O(\beta)$ and $T_p \le O(\beta)$.
\el
Variants of this lemma have been proved numerous times in the lace expansion literature so we  not prove it here (see for instance the proof of \cite[Lemma 5.5]{BorChaHofSlaSpe05b}).

\subsection{Bounds on the expected volume of critical clusters in a ball: proof of \refeq{pcVolBound}}\label{ChapUBVol}
The aim of this section is to bound the volume of a critical percolation configuration inside a Euclidean ball, that is, we prove bounds on
\begin{equation}
    \E_{p_c} [|\QR (0)|] = \sum_{x \in Q_r} \taupc (x) = \sum_{x \in \Zd} \taupc (x) \indi_{Q_r}(-x)= (\tau_{p_c}* \indi_{Q_r})(0).
\end{equation}
We start with the upper bound. We use the Fourier-space bounds on $\taupc$ proved in \cite{HeyHofSak08}. The Fourier transform of the indicator function may be negative, making it difficult to use. An alternative for the indicator function is a function $g_r (x)$ that satisfies the following criteria:
\begin{itemize}
    \item[(i)] $g_r (x) \ge \indi_{Q_r} (x)$ for all $x \in Q_r$;
    \item[(ii)] $\hat{g}_r (k) \ge 0 $ for all $k \in [-\pi, \pi]^d$;
    \item[(iii)] $\sum_{x \in \Zd} \taupc (x) g_r (-x) \le C r^{\twa}$.
\end{itemize}
When all three criteria are satisfied the result is the desired upper bound $\E_{p_c} [|\QR(0)|] \le C r^{\twa}$.

Let $\tqr = \{x \in \Zd : \|x\|_{\infty} \le r\}$.
As a choice of $g_r (x)$ we propose
\begin{equation}\label{e:grxdef}
    g_r (x)= (2 r+1)^d (p_r * p_r)(x) \qquad \text{ with } \qquad p_r (x) = \frac{ \indi_{\tqr}(x)}{(2\lfloor r \rfloor +1)^{d}},
\end{equation}
so $p_r$ is in fact a probability distribution. It is easy to verify criterion (i), and (ii) follows from the fact that $\hat{g}_r (k) = (2 r+1)^d \hat{p}_r (k)^2$. All that is left is to check criterion (iii):
We start by mentioning the bounds
\begin{equation}\label{e:htaubd}
    0 \le \hat{\tau}_p (k) \le \frac{1+O(\beta)}{1- \hat{D}(k)}.
\end{equation}
The lower bound was established in \cite{AizNew84} for all percolation two-point functions. In \cite{HeyHofSak08}, the upper bound for long-range percolation in dimension $d > 3\twa$ is proved and in \cite{HarSla90a} the same bound is proved for finite-range spread-out models in dimensions $d>6$ and for nearest-neighbor percolation in dimension $d \ge 19$.
Furthermore,
\begin{equation}\label{e:lbDhat}
    1 - \hat{D}(k)  \ge \left\{
                \begin{array}{ll}
                    c_1 L^{\twa} \lvert k \rvert^{\twa} &\text{for } \|k\|_\infty \le L^{-1};\\
                    c_2 &\text{for } \|k\|_\infty > L^{-1}
                \end{array} \right.
\end{equation}
have been established in \cite{HeyHofSak08}.
The Fourier transform of $(p_r * p_r) (x)$ is
\begin{equation}
    \widehat{(p_r * p_r)}(k) = \frac{ \widehat{\indi_{\tqr}}(k)^2 }{(2 \lfloor r \rfloor +1)^{2d}}= \frac{1}{(2 \lfloor r \rfloor +1)^{2d}} \prod_{i=1}^{d} \left(\frac{\sin([2 \lfloor r \rfloor +1] k_i /2)}{\sin(k_i / 2)}\right)^2,
\end{equation}
that is, it is the $d$-dimensional Dirichlet kernel squared.
Since $p_r$ is a probability distribution, its Fourier transform has a maximum value of 1.
Using the above bounds,
\begin{equation}\begin{split}
    (\taupc * g_r)(0) &= (2 r+1)^d \intFourier{[-\pi,\pi]^d}{\hat{\tau}_{p_c} (k) \hat{p}_r (k)^2}\\
    &\le C_A r^d \intFourier{|k| \le 1/r}{\frac{1}{|k|^{\twa}}} + C_B r^{d+\twa} \intFourier{|k| \ge 1/r}{\hat{p}_r(k)^2}.
    \end{split}
\end{equation}
Here, $C_A$ and $C_B$ are both positive constants that depend only on $\alpha,\: L$ and $d$. We bound the two terms separately.
For the first term on the right-hand side we have
\begin{equation}\label{e:pcA}
    C_A r^d \intFourier{|k| \le 1/r}{\frac{1}{|k|^{\twa}}} = C'_A r^d \int\limits_{0}^{1/r} k^{d-1-\twa} \dnd k \le C r^{\twa}.
\end{equation}
To bound the second term on the right-hand side we can extend the integration over $k$ to $[-\pi, \pi]^d$ and obtain
\begin{equation}\label{e:pcB}
    \begin{split}
     C_B r^{d+\twa} \intFourier{|k| \ge 1/r}{\hat{p}_r(k)^2} &\le C_B r^{d+\twa} \intFourier{[-\pi, \pi]^d}{\hat{p}_r (k)^2}\\
      & = \frac{C_B r^{d+\twa}}{(2 \lfloor r \rfloor +1)^{2d}} (\indi_{\tqr} * \indi_{\tqr})(0) = \frac{C_B r^{d+\twa}}{(2 \lfloor r \rfloor +1)^{d}} \le C r^{\twa}.
    \end{split}
\end{equation}
Combining the bounds for \refeq{pcA} and \refeq{pcB} yields the desired upper bound.

In much the same way we can determine a lower bound for $\E_{p_c} [|\QR (0)|]$.
Now we use a function $h_{r}(x)$, satisfying
\begin{itemize}
    \item[(i)] $h_{r} (x) \le \indi_{Q_r} (x)$ for all $x \in Q_r$;
    \item[(ii)] $\hat{h}_{r} (k) \ge 0 $ for all $k \in [-\pi, \pi]^d$;
    \item[(iii)] $\sum_{x \in \Zd} \tau_{p_c} (x) h_{r} (-x) \ge c r^{\twa}$.
\end{itemize}
With these criteria we get a lower bound $\E_{p_c} [|\QR(0)|] \ge c r^{\twa}$.

We choose
\begin{equation}\label{e:hrxdef}
    h_r (x) =  \frac{r^d}{d^{d/2}} (p_{q} * p_{q})(x) \qquad\text{ where }\qquad q = \left\lfloor \frac{r}{2 \sqrt{d}} \right\rfloor.
\end{equation}
That criterion (i) holds follows from the fact that $h_r(x)$ is monotonically decreasing in $\lvert x \rvert$,
\begin{equation}
    h_r(0) =  \frac{r^d}{d^{d/2}} \frac{\lvert \Lambda_q \rvert}{(2q+1)^{2d}} \le 1
\end{equation}
by the choice of $q$ and $h_r(x)=0$ for all $x$ such that $\lvert x \rvert > r$.
Just as in the case of the upper bound, criterion (ii) is easily verified.
To verify the last criterion, we use a bound on $1-\hat{D}(k)$ from \cite[Lemma 1.1]{Heyd11}:
for $\vep > 0$ small enough,
\begin{equation}\label{e:Dhatlb}
   1- \hat{D} (k) \le w_{\alpha} |k|^{\twa}, \qquad\text{ when } |k| \le \vep,
\end{equation}
where $0 <w_{\alpha} \le O(L^{\twa})$ is a constant.
Then
\begin{equation}
        (\taupc * h_{r})(0) = \intFourier{[-\pi,\pi]^d}{\hat{\tau}_{p_c}(k) \hat{h}_{r} (k)}
                            = \frac{r^d}{d^{d/2}(2q+1)^{2d}} \intFourier{[-\pi,\pi]^d}{\hat{\tau}_{p_c}(k) \prod_{i=1}^{d} \left(\frac{\sin([2q+1]k_i /2)}{\sin(k_i/ 2)}\right)^2}.
\end{equation}
Since $x^2/2 \le \sin(x)^2 \le x^2$ for $\lvert x \rvert \le 1$, we bound, for $(2q+1) \lvert k_i \rvert \le 1$,
\begin{equation}
        \frac{1}{(2q+1)^{2d}}\prod_{i=1}^{d} \left(\frac{\sin([2q+1]k_i /2)}{\sin(k_i /2)}\right)^2 \ge \frac{1}{(2q+1)^{2d}}\prod_{i=1}^{d}  \frac{ (\frac{2q+1}{2})^2 k_{i}^2}{k_{i}^2/8} \ge 2^{-d},
\end{equation}
Using this bound, we get
\begin{equation}\label{e:pcLowB}
        \intFourier{[-\pi,\pi]^d}{\hat{\tau}_{p_c}(k) \hat{h}_{r} (k)} \ge \frac{c r^d}{d^{d/2} 2^d} \intFourier{|k| \le \vep/r}{\frac{1}{|k|^{\twa}}}
        \ge \frac{c r^{d}}{d^{d/2} 2^d r^{d-\twa}} \ge c(d) r^{\twa}
\end{equation}
for $\vep >0$ sufficiently small and $c(d)$ a constant depending on $d$. This completes the proof of \refeq{pcVolBound}.
\qed

\subsection{Bounds on the expectation of the backbone volume: proof of \refeq{BBVolBound}}
Backbone edges are those edges that have a path from $0$ to one end of the edge, disjointly from a path from the other end of the edge to infinity. Therefore
\begin{equation}\label{e:bbBall}
    \Eiic[\Nbb (r)] = \sum_{b \in \Be} \Eiic [\indi_{\{b \text{ is a backbone edge}\}}] =  \sum_{b \in \Be} \Eiic [\indi_{\{0 \conn \ulb\}\circ\{b \text{ open}\}\circ\{\olb \conn \infty\}}].
\end{equation}
Backbone events are by definition \emph{not} cylinder events, and hence it is a priori unclear whether the limiting scheme that yields $\Qiic$ can be reversed.
The aim of this section is to show that we can.

We call an open edge $b=\{x,y\}\in \Zd$ \emph{backbone-pivotal}
when every infinite self-avoiding walk in the IIC starting at the origin
makes use of this edge.
The backbone-pivotal edges can be ordered as $(b_i)_{i=1}^{\infty}$, in such a way that
every infinite self-avoiding walk starting at the origin passes through $b_i$ before passing through $b_{i+1}$.
Also, we can think of the backbone-pivotal edges as being \emph{directed} edges $b=(x,y)$, where the direction is such that $\{0\conn x\}$ makes use of different edges than $\{y\conn \infty\}$.
For a directed edge $b=(x,y)$, we let $\ulb=x$ denote its \emph{bottom}, and $\olb=y$ its \emph{top}.
Writing $b_m$ for the $m$th backbone-pivotal edge, we define
\begin{equation}
     \smi\equiv \tilde\Ccal^{b_m}(0) \setminus \tCcal^{b_{m-1}}(0)
\end{equation}
to be the subgraph of the $m$th ``backbone sausage'' (where, by convention $\tCcal^{b_{0}}(0) = \varnothing$).

If $0$ is connected with $\dqr$, and there are precisely $n$ open pivotal edges for this connection, we can again impose an ordering on the open pivotal edges $(b_i)_{i=1}^n$ in such a way that any self-avoiding path from 0 to $\dqr$ passes through $b_i$ before passing through $b_{i+1}$.
If $n\ge m$, we let $\smr\equiv \tilde\Ccal^{b_m}(0) \setminus \tCcal^{b_{m-1}}(0)$ and we let $\smr=\varnothing$ whenever $0\nleftrightarrow \dqr$ or $n<m$.

In the same way, we let $\smx\equiv \tilde\Ccal^{b_m}(0)\setminus \tCcal^{b_{m-1}}(0)$ where $b_m$ now is the $m$th open pivotal edge for $\{0\conn x\}$, and $\smx=\varnothing$ if no $m$th pivotal bond exist for the connection $\{0\conn x\}$.

We are interested in events that take place on the first $m$ backbone sausages. To this end, define
\begin{equation}
    \sni \equiv \bigcup_{i=1}^m \mathsf{S}_{i}^{\sss \infty}, \qquad \snr \equiv \bigcup_{i=1}^m \mathsf{S}_{i}^{\sss (r)}, \qquad \text{ and } \qquad \snx \equiv \bigcup_{i=1}^m \mathsf{S}_{i}^x.
\end{equation}
Even though events occurring on $\sni$ are not necessarily cylinder events, it is still possible to reverse the IIC-limit for such events, as the following lemma demonstrates.
\bl[Backbone limit reversal lemma] \label{lem:BackboneReversal}
For any event $E$ and any $m \in \N$,
\begin{equation}
    \Qiic\left(E \text{\emph{ on }}\sni\right) = \limp \frac{1}{\chi(p)} \sum_{x \in \Zd} \Pp\left(\{E \text{\emph{ on }} \snx\}\cap\{0\conn x\}\right).
\end{equation}
\el
\proof
Fix $m$ throughout the proof. We prove the lemma via comparison of $\sni$, $\snr$ and $\snx$.
To this end, we define the events
\begin{eqnarray}
    \liR &\equiv& \left\{\omega \colon \sni = \snR\right\};\\
    \lrR &\equiv& \left\{\omega \colon \snr = \snR \text{ and there are at least $m$ pivotals for }\{0 \conn Q_r^c\}\right\};\\
    \lrx &\equiv& \left\{\omega \colon \snr = \snx \text{ and there are at least $m$ pivotals for }\{0 \conn Q_r^c\}\right\}.
\end{eqnarray}
We show that it is improbable that these sets are different (when compared within the same configuration, near the origin), so that, upon taking a suitable limit, replacing one with the other is justified.

We start by observing that for any $R$,
\begin{equation}\label{e:splitF}
    \left\{E \text{ on } \sni\right\} = \left(\left\{E \text{ on } \sni \right\} \cap \liR \right)\dot \cup\left(\left\{E \text{ on } \sni\right\} \cap \bliR \right) \equiv F_m^1 (R) \dot \cup F_m^2 (R).
\end{equation}
At the end of the proof we take the limit $R \to \infty$.
In this limit, the event $\bliR$ has probability $0$ under $\Qiic$ as it implies that there exists a path from one of the first $m$ sausages to $Q_R^c$ disjoint of the backbone, which in the limit $R \to \infty$ implies the existence of \emph{two} disjoint connections to $\infty$; an event that does not occur $\Qiic$-almost surely.

For $F_m^1 (R)$, the occurrence of $\liR$ implies $\{E \text{ on }\sni\} = \{E \text{ on }\snR\}$.
Furthermore, for any $r$ such that $0 < r < R$ we can write
\begin{equation}\label{e:splitG}
    F_m^1 (R) = \left(\left\{E \text{ on } \snR \right\} \cap \liR \cap \lrR \right) \dot \cup\left(\left\{E \text{ on } \snR\right\} \cap \liR \cap \blrR \right) \equiv G_m^1 (R,r) \dot \cup G_m^2 (R,r).
\end{equation}
In the double limit where first $R\to\infty$ and then $r\to\infty$, the probability of $G_m^2 (R,r)$ vanishes as
\begin{equation}\label{e:G2bound}
	\lim_{r \to \infty} \lim_{R \to \infty} \Qiic(G_m^2 (R,r))
	\le\lim_{r \to \infty} \lim_{R \to \infty} \Qiic(\blrR) = \lim_{r \to \infty} \Qiic(\blir) =0,
\end{equation}
since otherwise again there are two disjoint paths to $\infty$.

We can rewrite $G_m^1 (R,r)$ as follows:
\begin{equation}\label{e:splitH}
    G_m^1 (R,r) = \left(\left\{E \text{ on } \snR\right\} \cap \lrR \right) \setminus \left(\left\{E \text{ on } \snR\right\} \cap \lrR \cap \bliR\right) \equiv H_m^1 (R,r) \setminus H_m^2 (R,r).
\end{equation}
Since $H_m^2 (R,r) \subseteq \bliR$ we again have that $\Qiic(H_m^2 (R,r)) \to 0$ as $R \to \infty$.

Now, $H_m^1 (R,r)$ is a cylinder event, so that \eqref{def:Qiic} applies,
\begin{equation}
    \Qiic(H_m^1 (R,r)) = \limp \frac{1}{\chi(p)} \sum_{x \in Q_R^c} \Pp\big(H_m^1 (R,r)\cap \{0 \conn x\}\big),
\end{equation}
(where the sum over $x\in Q_R$ vanishes in the $p\nearrow p_c$ limit).

The crucial observation is that for $r<R$ and $x\in Q_R^c$ we have $\lrR=\lrx$, so that
\begin{equation}
	\left\{E \text{ on } \snR\right\} \cap \lrR
	=\left\{E \text{ on } \snx\right\} \cap \lrR.
\end{equation}
Consequently,
\begin{equation}\label{e:splitM}
    \begin{split}
        H_m^1 (R,r) \cap \{0 \conn x\} &= \left\{E \text{ on } \snR\right\} \cap \lrR \cap \{0 \conn x\}
        	 = \left\{E \text{ on }\snx\right\}\cap \lrR \cap \{0 \conn x\} \\
                &= \left(\left\{E \text{ on } \snx\right\} \cap \{0 \conn x\}\right) \setminus \left(\left\{E \text{ on } \snx\right\} \cap \blrR \cap \{0 \conn x\}\right) \equiv M_m^1 (x) \setminus M_m^2 (R,r,x).
    \end{split}
\end{equation}
For $M_m^2 (R,r,x)$ we note that $\blrR$ is a cylinder event, so that \eqref{def:Qiic} implies
\begin{equation}\label{e:M2bound}
    \begin{split}
        \lim_{r \to \infty}\lim_{R \to \infty}  \limp \frac{1}{\chi(p)} \sum_{x \in Q_R^c} \Pp(M_m^2 (R,r,x))
        &\le \lim_{r \to \infty} \lim_{R \to \infty} \limp \frac{1}{\chi(p)} \sum_{x \in Q_R^c} \Pp(\blrR, 0 \conn x)\\
        &\le \lim_{r \to \infty} \lim_{R \to \infty} \Qiic(\blrR) = \lim_{r \to \infty} \Qiic(\blir) =0.
    \end{split}
\end{equation}

Combining \refeq{splitF}--\refeq{splitM},
\begin{eqnarray}
    \Qiic\left(E  \text{ on }  \sni\right) &= &
    \Qiic(F_m^2 (R)) + \Qiic(G_m^2 (R,r)) - \Qiic(H_m^2 (R,r))\nnb
     &&{}+ \limp \frac{1}{\chi(p)} \sum_{x \in Q_R^c}
     \big(\Pp(M_m^1 (x))  - \Pp(M_m^2 (R,r,x))\big).
\end{eqnarray}
Now we add $0=\limp {\chi(p)^{-1}} \sum_{x \in Q_R} \Pp(M_m^1(x))$
to the right hand side, so that the term involving $M^1_m(x)$ is independent of $r$ and $R$.
Then we let $R\to\infty$, so that $\Qiic(F_m^2 (R))$ and $\Qiic(H_m^2 (R,r))$ vanish.
Subsequently, we let also $r\to\infty$, so that also the terms involving $G_m^2 (R,r)$ and $M_m^2 (R,r,x)$ disappear, by \eqref{e:G2bound} and \eqref{e:M2bound}.
The result is
\begin{equation}
    \Qiic\left(E \text{ on } \sni\right) = \limp \frac{1}{\chi(p)} \sum_{x\in \Zd} \Pp(M_m^1 (x)) = \limp \frac{1}{\chi(p)} \sum_{x \in \Zd} \Pp\left(\left\{E \text{ on } \snx\right\}\cap\{0\conn x\}\right),
\end{equation}
completing the proof.
\qed

\medskip
Let $\Bb(\omega)$ denote the backbone edge set of a configuration $\omega$, and let $\mathsf{S}[A,B](\omega)$ denote the set of open edges between the sets $A$ and $B$,
that is, $\{u,v\}\in\mathsf{S}[A,B](\omega)$ whenever $\{u,v\}$ is open and $\{a\conn u\}\circ\{v\conn b\}$ for some $a\in A,b\in B$.
Similarly, write $\Bbp(\omega)$ for the set of backbone pivotal edges, and $\Sp [A,B](\omega)$ for the set of open pivotal edges for the event that there exists a connection between the sets $A$ and $B$.

In this paper we use two specific cases of the above lemma.
\bc[Backbone limit reversal lemma for sets of edges]\label{BPRlem}
Let $\{b_i\}_{i=1}^n$ be a fixed and finite set of edges. Then,
\begin{enumerate}
    \item   \begin{equation}
                \Qiic\left(\{b_i\}_{i=1}^{n} \subseteq \Bb\right) = \limp \frac{1}{\chi(p)} \sum_{x \in \Zd} \Pp\left(\{b_i\}_{i=1}^{n} \subseteq \mathsf{S}[0,x]\right);
            \end{equation}
    \item   \begin{equation}
                \Qiic\left(\{b\}_{i=1}^n  \subseteq \Bbp\right) = \limp \frac{1}{\chi(p)} \sum_{x \in \Zd} \Pp \left(\{b_i\}_{i=1}^n \subseteq \Sp[0,x]\right).
            \end{equation}
\end{enumerate}
\ec
\proof The proof for both cases follows by the same argument, so we only prove it for (i). Define
\begin{eqnarray}
    A_{m} \equiv \left\{\{b_i\}_{i=1}^{n} \subseteq \sni\right\},
    &&
    A_{\infty} \equiv \left\{\{b_i\}_{i=1}^{n} \subseteq \Bb\right\} =  \bigcup_{m=1}^{\infty} A_{m};\\
    B_{m}(x) \equiv \left\{\{b_i\}_{i=1}^{n} \subseteq \snx\right\}\cap \{0 \conn x\};
    &&
    B_{\infty}(x) \equiv \left\{\{b_i\}_{i=1}^{n} \subseteq \mathsf{S}[0,x]\right\}\cap\{0 \conn x\} =  \bigcup_{m=1}^{\infty} B_{m}.
\end{eqnarray}
Since $A_m\subseteq A_{m+1}$ for $m\ge1$, we partition $A_\infty$ as $A_\infty=A_1\cup\,\bigcup_{m\ge1}(A_{m+1}\setminus A_m)$, where the union is over disjoint subsets.
We can write a similar partition for $B_\infty(x)$, for every $x\in\Zd$.
Next we apply Lemma \ref{lem:BackboneReversal} to each $A_m$-term, followed by dominated convergence to deduce
\begin{eqnarray}
	\Qiic(A_\infty)
	&=&\Qiic(A_1)+\sum_{m=1}^{\infty} \left(\Qiic(A_{m+1})-\Qiic(A_m)\right)\nnb
	 &=&\limp\frac{1}{\chi(p)}\sum_{x\in\Zd}\Pp(B_1(x))+\sum_{m=1}^{\infty}\left(\limp\frac{1}{\chi(p)}\sum_{x\in\Zd}\Pp(B_{m+1}(x))-\limp\frac{1}{\chi(p)}\sum_{x\in\Zd}\Pp(B_m(x))\right)\nnb
	&=&\limp\frac{1}{\chi(p)}\sum_{x\in\Zd}\left(\Pp(B_1(x))+\sum_{m=1}^{\infty} \left(\Pp(B_{m+1}(x))-\Pp(B_{m}(x))\right)\right).
\end{eqnarray}
Since clearly $\Pp(B_m (x)) \to \Pp(B_\infty (x))$ as $m \to \infty$, the telescoping sum on the right-hand side is equal to $\Pp(B_{\infty}(x)) - \Pp(B_1(x))$, so
\begin{equation}
	\Qiic(A_{\infty}) = \limp\frac{1}{\chi(p)}\sum_{x\in\Zd}\Pp(B_\infty(x)),
\end{equation}
as we set out to prove.
\qed

\subsubsection{Upper bound on the expectation of the backbone volume}
Applying Corollary \ref{BPRlem}(i) to \refeq{bbBall}, we obtain
\begin{equation}\label{e:bbDisjoint}
    \begin{split}
        \Eiic[\Nbb (r)]= & \sum_{b \in \Be} \Qiic(b \in \Bb) = \limp \frac{1}{\chi(p)} \sum_{b \in \Be} \sum_{x \in \Zd} \Pp(b \text{ open and pivotal for }0 \conn x)\\
         \le & \limp \frac{1}{\chi(p)} \sum_{b \in \Be} \sum_{x \in \Zd} \E_p [\indi_{\{0 \conn \ulb\}\circ\{b \text{ open}\}\circ\{\olb \conn x \}}].
    \end{split}
\end{equation}
Applying the BK-inequality to \refeq{bbDisjoint} gives
\begin{equation}
     \Eiic[\Nbb (r)] \le \limp \frac{1}{\chi(p)} \sum_{b \in \Be} \sum_{x \in \Zd} \taup (\ulb) pD(b) \taup(x-\olb).
\end{equation}
Summing over $x$ and then $\olb$ and bounding $b \in \Be$ by $\ulb \in Q_r$, we obtain a factor $\chi(p)$ and a factor $p$, respectively, after which we take the limit $p \ua p_c$:
\begin{equation}\label{e:BBsum}
    \Eiic[\Nbb (r)] \le p_c \sum_{\ulb \in Q_r} \taupc(\ulb).
\end{equation}
The upper bound in \refeq{pcVolBound}, which is proved in the previous section, completes the proof of upper bound in \refeq{BBVolBound}. \qed

\subsubsection{Lower bound on the expectation of the backbone volume}
To establish a lower bound on $\Eiic[\Nbb (r)]$ we count only the backbone-pivotal edges.
Recall the definition of $h_r$ given in \refeq{hrxdef}. We bound
\begin{equation}\label{e:bbBall1}
        \Eiic[\Nbb(r)]  = \sum_{b \in \Be} \Qiic(b \in \Bb) \ge \sum_{b \in \Be} \Qiic(b \in \Bbp) \ge \sum_{b \in \Zd \times \Zd} h_r (\ulb) \Qiic(b \in \Bbp).
\end{equation}
That the second inequality is necessary is not immediately obvious, but it will turn out to be crucial for obtaining a good bound in the case of nearest-neighbor percolation. Now we apply Corollary \ref{BPRlem}(ii) to obtain
\begin{equation}\label{e:bbBall2}
    \sum_{b \in \Zd \times \Zd} h_r (\ulb) \Qiic(b \in \Bbp) = \limp \frac{1}{\chi(p)} \sum_{x, \ulb,\olb \in \Zd \times \Zd} h_r (\ulb) \Pp(b \in \Sp[0,x]).
\end{equation}
By the definition of $\Sp[0,x]$ we have
\begin{equation}
    \{b \in \Sp[0,x]\} = \{0 \conn \ulb \text{ on }\tCcal^b (0)\} \circ \{b \text{ open}\} \circ \{\olb \conn x \text{ off } \tCcal^b (0)\},
\end{equation}
so we can apply the Factorization Lemma to the right-hand side of \refeq{bbBall2} to obtain
\begin{equation}\label{e:Nbbrbd}
    \begin{split}
        \Eiic[\Nbb (r)] &\ge \limp \frac{1}{\chi(p)}  \sum_{x,\ulb,\olb \in \Zd} p D(b) h_r (\ulb) \E_0[\indi_{\{0 \conn \ulb \text{ on }\tCcal^{b} (0)\}} \E_1[\indi_{\{\olb \conn x \text{ off }\tCcal^{b}(0)\}}]]\\
                        &= \limp \frac{1}{\chi(p)} \sum_{x,\ulb,\olb \in \Zd} p D(b) h_r (\ulb) \E_0[\indi_{\{0 \conn \ulb\}} \E_1[\indi_{\{\olb \conn x \text{ off }\tCcal^{b}(0)\}}]].
    \end{split}
\end{equation}
In the second line we omitted the condition ``on $\tCcal^{b}(0)$'' because
\begin{equation}
    \{0 \conn \ulb \text{ on } \tCcal^{b}(0)\} = \{0 \conn \ulb\} \setminus \{0 \conn \ulb \text{ through }\Zd \setminus \tCcal^b(0)\},
\end{equation}
but $\{0 \conn \ulb \text{ through }\Zd \setminus \tCcal^b(0)\}$ implies that $b$ is pivotal for the connection $\{0 \conn \ulb\}$, which means that the event $\{0 \conn \olb\}$ has to occur. However, the indicator $\indi_{\{\olb \conn x \text{ off }\tCcal^{b}(0)\}}$ is always $0$ for such events, so the change from $\{0 \conn \ulb \text{ on } \tCcal^{b}(0)\}$ to $\{0 \conn \ulb\}$ has no effect on the expectation.

We write $\tCcal^b_0 (0)$ to remind us that the cluster is random w.r.t. $\E_0$, but fixed w.r.t. $\E_1$. We bound the expectations in \refeq{Nbbrbd} from the inside out:
\begin{equation}
         \Eiic[\Nbb (r)] \ge \limp \frac{1}{\chi(p)} \sum_{x,\ulb,\olb \in \Zd} p D(b) h_r (\ulb) \Big( \taup (\ulb)\taup(x - \olb)- \E_0\Big[ \indi_{\{0 \conn \ulb\}} \P_p \bigg(\olb \stackrel{\tCcal^b_0 (0)}{\longleftrightarrow} x \bigg)\Big]\Big) \equiv N_1 - N_2.
\end{equation}
For the inequality we used the identity $\{E $ off $A\}=E \setminus \{E$ through $A \}$. We bound $N_1$ and $N_2$ separately.

Consider $N_1$ first:
\begin{equation}
    \begin{split}
        N_1     &= \limp \frac{1}{\chi(p)} \sum_{x, \ulb, \olb \in \Zd} h_r (\ulb) \taup (\ulb) p D(b) \taup(x - \olb)\\
                &= p_c \sum_{\ulb \in \Zd} h_r (\ulb) \taupc (\ulb) =p_c \intFourier{[-\pi,\pi]^d}{\htaupc(k) \hat{h}_r(k)}.
    \end{split}
\end{equation}
To obtain the second equality we summed over $x$ and then $\olb$, as we did for the upper bound.

The bound on $N_2$ is harder:
\begin{equation}\label{e:N2intermediate}
        N_2 = \limp \frac{1}{\chi(p)} \sum_{x, \ulb, \olb \in \Zd} p D(b) h_r (\ulb) \E_0\Big[ \indi_{\{0 \conn \ulb\}} \P_p \bigg(\olb \stackrel{\tCcal^b_0 (0)}{\longleftrightarrow} x \bigg)\Big]
\end{equation}
and note that here we need an upper bound.
The dependence in the second connectivity function effectively implies that there is a path from some vertex along the path $\olb \conn x$ to another vertex on the path $0 \conn \ulb$, and that this path does not use the edge $b$.
Consider a fixed set of vertices $A \subset \Zd$. Then,
\begin{equation}
    \{\olb \conn x \text{ through }A\} \subseteq \bigcup_{a \in A} \{\olb \conn a\}\circ\{a \conn x\}.
\end{equation}
Therefore,
\begin{equation}\label{e:genEvent}
    \begin{split}
        \Pp(\olb \Aconn x)  &\le \Pp\left(\bigcup_{a \in A} \{\olb \conn a\}\circ\{a \conn x\}\right) \le \sum_{a \in \Zd} \indi_{\{a \in A\}} \Pp(\{\olb \conn a\}\circ\{a \conn x\})\\
                            &\le \sum_{a \in \Zd} \indi_{\{a \in A\}} \taup(a-\olb)\taup(x-a).
    \end{split}
\end{equation}

Since the set $\tCcal^b_0 (0)$ is fixed with respect to the expectation $\E_1$ we may apply \refeq{genEvent} to the expectation on the right-hand side of \refeq{N2intermediate} with $A = \Ccal_0 (0) \supset \tCcal^b_0 (0)$:
\begin{equation}\label{e:N2m}
    \begin{split}
        N_2 &\le  \limp \frac{1}{\chi(p)} \sum_{x, \ulb, \olb, a \in \Zd} p D(b) h_r (\ulb) \E_0\left[ \indi_{\{0 \conn \ulb\}} \indi_{\{0 \conn a\}}\taup(a-\olb)\taup(x-a)\right]\\
            &= \limp \frac{1}{\chi(p)} \sum_{x, \ulb, \olb, a \in \Zd}  p D(b) h_r (\ulb) \Pp(0 \conn \ulb, 0 \conn a)\taup(a-\olb)\taup(x-a).
    \end{split}
\end{equation}
We use the tree-graph bound \cite{AizNew84}:
\begin{equation}\label{e:TreeGraph}
    \Pp(0 \conn \ulb, 0 \conn a) \le \sum_{z \in \Zd} \taup(z)\taup(\ulb-z)\taup(a-z)
\end{equation}
and insert the above inequality into \refeq{N2m} to obtain
\begin{equation}
    \begin{split}
        N_2 &\le \limp \frac{1}{\chi(p)} \sum_{x,\ulb,\olb,a,z \in \Zd}  h_r (\ulb) \taup(z) \taup(\ulb -z) p D(b) \taup (a-\olb) \taup(a -z) \taup(x-a)\\
            &= p_c \sum_{\ulb,\olb,a,z \in \Zd} h_r(\ulb) \taupc (z) \taupc (\ulb -z)  D(b) \taupc(a - \olb) \taupc(a -z).
    \end{split}
\end{equation}
Define
\begin{equation}
        T'_{p_c} (x)   = \taupc(x) (D* \taupc * \taupc)(x).
\end{equation}
An upper bound on its Fourier transform is
\begin{equation}\label{e:Tppcbd}
        |\hat{T}'_{p_c} (k)| = \left\lvert\sum_{x \in \Zd} \e^{i k\cdot x} T'_{p_c} (x)\right\rvert \le |\hat{T}'_{p_c}(0)| \le \sum_{v,w,y\in \Zd}  D(v) \taupc(w-v) \taupc(y-w) \taupc(y) = T_{p_c}(0) \le C \beta,
\end{equation}
with $T_p (x)$ as given by $\refeq{triangleT}$. The bound on $T_{p_c}(0)$ follows from Lemma \ref{orderbounds}.

With this definition we can write
\begin{equation}
    N_2 \le p_c \sum_{\ulb, \olb,a,z \in \Zd} h_r(\ulb) T'_{p_c}(\ulb-z) \taupc(z) = p_c (\taupc * T'_{p_c} * h_r )(0).
\end{equation}
We can bound $N_2$ by expressing the right-hand side in terms of its Fourier transform:
\begin{equation}\label{e:N2bd}
        N_2 \le p_c (\taupc * T'_{p_c} * h_r)(0) = p_c \intFourier{[-\pi, \pi]^d}{\htaupc (k) \hat{T'}_{p_c} (k) \hat{h}_r (k)} \le  C \beta \intFourier{[-\pi, \pi]^d}{\htaupc (k) \hat{h}_r (k)},
\end{equation}
where the second inequality follows from \refeq{Tppcbd}.

With bounds on both $N_1$ and $N_2$ we can conclude that, when $\beta$ is small enough,
\begin{equation}
    \Eiic[\Nbb (r)] \ge N_1 - N_2 \ge  p_c (1- C \beta)\intFourier{[-\pi, \pi]^d}{\htaupc(k) \hat{h}_r (k)} \ge c'(d) r^{\twa} .
\end{equation}
for some constant $c'(d)$ that only depends on $d$. The final inequality follows from \refeq{pcLowB}.
This concludes the proof of the lower bound. The upper and lower bound combined complete the proof of \refeq{BBVolBound}. \qed
\subsection{Bounds on the expected IIC volume in a ball: proof of \refeq{IICVolBound} }
Define the IIC connectivity function
\begin{equation}\label{ConnIIC}
        \rho(y) \equiv \Qiic (0 \conn y)
               = \limp \frac{1}{\chi(p)} \sum_{x \in \Zd} \Pp (0 \conn y, 0 \conn x).
\end{equation}
Since the event $\{0 \conn y\}$ is not a cylinder event, it is not immediately obvious that we can write it as a limit. However, in \cite{HofJar04} it is proved that this is allowed.

Using the techniques of the previous paragraphs, we can easily find an upper bound. The lower bound requires more work.

\subsubsection{IIC volume expectation upper bound}
We start by bounding (\ref{ConnIIC}) using the tree-graph bound \refeq{TreeGraph}:
\begin{equation}
    \rho(y) \le \limp \frac{1}{\chi(p)} \sum_{x,z \in \Zd} \tau_p (z) \tau_p (x-z) \tau_p (y-z).
\end{equation}
Keeping $z$ fixed and summing over $x$ we get a factor $\chi(p)$. Then, with the divergence of the susceptibility canceled, we can take the limit $p \ua p_{c}$:
\begin{equation}
        \rho(y) \le \sum_{z \in \Zd} \tau_{p_c} (z) \tau_{p_c} (y-z)
                = (\tau_{p_c} * \tau_{p_c})(y).
\end{equation}
The expected volume of the IIC in a Euclidean ball is given by
\begin{equation}
    \Eiic [|Q_r \cap \iic|] = \sum_{y \in Q_r} \rho(y).
\end{equation}
Using the same techniques as in the proof of \refeq{pcVolBound}, we obtain
\begin{equation}
    \begin{split}
        \Eiic[|Q_r \cap \iic|]   &\le  \sum_{x, y \in \Zd} \tau_{p_c} (x) \tau_{p_c} (y-x) g_r (-y)
                        \le  (2 r)^d \int\limits_{[-\pi, \pi]^d} \frac{\hat{p}_r (k)^{2}}{[1 - \hat{D}(k)]^2} \frac{\dd k}{(2 \pi)^d}\\
                        &\le C_A r^d \intFourier{|k| \le 1/r}{\frac{1}{|k|^{2\twa}}}+ C_B r^{d+\twa} \intFourier{|k|\ge 1/r}{\hat{p}_r (k)^2} \le C r^{2\twa}.
    \end{split}
\end{equation}
\qed

\subsubsection{IIC volume expectation lower bound}
This bound is the most demanding one, as we are required to use the Factorization Lemma twice.
We bound (\ref{ConnIIC}) from below by
\begin{equation}
        \rho(y) \ge \limp \frac{1}{\chi(p)} \sum_{x, \ulb,\olb \in \Zd} \P_p\left( \begin{array}{l} 0\conn y, 0\conn x, b=(\ulb,\olb) \text{ is the first edge that is}\\ \text{open and pivotal for }0\conn x \text{ but not for }0\conn y
    \end{array} \right)
 \end{equation}
Observe that
\begin{equation}
    \begin{split}
        \{0\conn y, 0\conn x, &\, b \text{ is the first edge that is open and pivotal for }0\conn x \text{ but not for }0\conn y\} \\
            &= \{b \text{ open}\} \circ \{\{0 \conn \ulb\}\circ\{ \ulb \conn y \} \text{ on } \tilde{\Ccal}^b(0)\} \circ \{\olb \conn x \text{ on } \Zd \setminus \tilde{\Ccal}^b (0)\}.
    \end{split}
\end{equation}
Applying the Factorization Lemma gives
\begin{equation}\label{e:FacLemI}
    \rho(y) \ge \limp \frac{1}{\chi(p)} \sum_{x,\ulb,\olb \in \Zd} p D(b) \E_0[\indi_{\{0 \conn \ulb\}\circ\{\ulb \conn y\}}\E_1[\indi_{\{\olb \conn x \text{ off }\tCcal^b_0 (0)\}}]],
\end{equation}
where we left out the condition ``on $\tCcal^{b}_0 (0)$'' again, for the same reason that we were allowed to leave it out in \refeq{Nbbrbd}.
For a fixed set of vertices $A$,
\begin{equation}\label{e:offtothrough}
    \Pp(x \conn y \text{ off }A ) = \taup(y-x) - \E_p[\indi_{\{x \conn y \text{ through }A\}}].
\end{equation}
Since $\tCcal^b_0 (0)$ is fixed with respect to $\E_1$ we may apply this identity to \refeq{FacLemI} and sum over $y \in Q_r$ to obtain
\begin{equation}
    \begin{split}
        \sum_{y \in Q_r} \rho(y) \ge & \limp \frac{1}{\chi(p)} \sum_{y \in Q_r} \sum_{x,\ulb,\olb \in \Zd}p D(b) \E_0 [ \indi_{\{0 \conn \ulb\} \circ \{\ulb \conn y\}} (\taup(x - \olb) - \E_1 [\indi_{\{\olb \conn x \text{ through }\tCcal^b_0 (0)\}}])]\\
        \ge & \limp \frac{1}{\chi(p)} \sum_{y \in \Zd} \sum_{x,\ulb,\olb \in \Zd}p D(b)h_r(y) \E_0 [ \indi_{\{0 \conn \ulb\} \circ \{\ulb \conn y\}} (\taup(x - \olb) - \E_1 [\indi_{\{\olb \conn x \text{ through }\tCcal^b_0 (0)\}}])]\\
         \equiv & S_1 - S_2.
    \end{split}
\end{equation}
In the second inequality we have again replaced the sum over $y \in Q_r$ by a sum over $y \in \Zd$ and inserted a factor $h_r (y)$. This is a necessary step for obtaining a good bound on $S_1$.

We first give an upper bound on $S_2$, and then establish a lower bound on $S_1$. As mentioned, the set $\tCcal^b_0 (0)$ is fixed w.r.t. $\E_1$. Hence, we obtain an upper bound on $S_2$ by making use of \refeq{genEvent} with $A = \Ccal_0 (0) \supset \tCcal^b_0 (0)$:
\begin{equation}
    S_2 \le \limp \frac{1}{\chi(p)} \sum_{x,y,\ulb,\olb,a \in \Zd} p D(b) h_r(y) \Pp(\{0 \conn \ulb\}\circ\{\ulb \conn y\}) \taup(a) \taup(a - \olb) \taup(x -a).
\end{equation}
We can now sum over $x$ to obtain a factor $\chi(p)$ and subsequently take the limit $p \ua p_c$.

By the BK inequality we can also bound
\begin{equation}
    \Pp(\{0 \conn \ulb\}\circ\{\ulb \conn y\}) \le \taup(\ulb)\taup(y-\ulb).
\end{equation}
Applying this bound we obtain
\begin{equation}
    \begin{split}
        S_2 \le& p_c  \sum_{y,\ulb,\olb, a \in \Zd} D(b) h_r(y) \taupc(\ulb)\taupc(y-\ulb)\taupc(a) \taupc(a - \olb)\\
            = & \sum_{y \in \Zd} h_r(y) p_c (\taupc * T'_{p_c} )(y) = p_c (\taupc * T'_{p_c} * h_r)(0),
    \end{split}
\end{equation}
where the last inequality follows from the symmetries of $h_r$.
We end up with the same bound as on $N_2$ in the previous section. Hence, by \refeq{N2bd}, \refeq{pcA} and \refeq{pcB}, we obtain
\begin{equation}\label{e:S2bd}
    S_2 \le C \beta r^{\twa}.
\end{equation}

We now establish an upper bound on $S_1$. Immediately we can sum over $x$ and $\olb$ to obtain factors $\chi(p)$ and $p$ and take the limit $p \ua p_c$:
\begin{equation}
    S_1 =   \sum_{y, \ulb,\olb \in \Zd} p_c^2 h_r(y) \Ppc(\{0 \conn \ulb \}\circ\{\ulb \conn y\}).
\end{equation}
Observe that
\begin{equation}
    \{0 \conn \ulb\}\circ \{\ulb \conn y\} \supseteq \bigcup_{e: \ule=\ulb} \{e \text{ is open and pivotal for }0 \conn y\}.
\end{equation}
and
\begin{equation}
    \{e \text{ is open and pivotal for }0 \conn y\} = \{0 \conn \ule \text{ on }\tCcal^e_0 (0)\} \cap \{e \text{ open}\} \cap \{\ole \conn y \text{ off } \tCcal^e (0)\}.
\end{equation}
Making this replacement, applying the Factorization Lemma again, and applying \refeq{offtothrough} we obtain the lower bound
\begin{equation}
        S_1 \ge \sum_{y, \ulb \in \Zd} \sum_{e: \ule = \ulb} p_c^2 h_r(y) D(e) \E_0[ \indi_{\{0 \conn \ule\}}
               (\taupc(y -\ole) - \E_1[ \indi_{\{\ole \conn y \text{ through }\tCcal^e_0 (0)\}}])] \equiv S_{1,1} -S_{1,2}
\end{equation}
where we again left out the condition ``on $\tCcal^{e}_0 (0)$'' for the same reason that we were allowed to leave it out in \refeq{Nbbrbd}.

Writing $S_{1,1}$ in terms of its Fourier transform, we obtain
\begin{equation}
    \begin{split}
        S_{1,1} &=\sum_{y, \ulb \in \Zd} \sum_{e: \ule = \ulb} p_{c}^2 D(e) h_r(y) \taupc(\ule)\taupc(y-\ole) =p_{c}^2 (\taupc * D * \taupc * h_r)(0)\\
                &= p_c^2 \intFourier{[-\pi,\pi]^d}{\htaupc(k)^2 \hat{D}(k) \hat{h}_r (k)}.
    \end{split}
\end{equation}
Rewriting the right-hand side gives
\begin{equation}\label{e:S1bound}
    \begin{split}
        S_{1,1} &= p_{c}^2  \intFourier{[-\pi, \pi]^d}{ \htaupc(k)^2 (1-[1-\hat{D}(k)]) \hat{h}_{r}(k)}\\
        &=p_{c}^2  \intFourier{[-\pi, \pi]^d}{ \htaupc(k)^2 \hat{h}_{r}(k)} - p_{c}^2  \intFourier{[-\pi, \pi]^d}{ \htaupc(k)^2 [1-\hat{D}(k)] \hat{h}_{r}(k)}.
    \end{split}
\end{equation}
For the second integral we obtain an upper bound:
\begin{equation}
        p_{c}^2  \intFourier{[-\pi, \pi]^d}{ \htaupc(k)^2 [1-\hat{D}(k)] \hat{h}_{r}(k)} \le C  \intFourier{[-\pi, \pi]^d}{ \frac{\hat{h}_{r}(k) [1-\hat{D}(k)]}{[1-\hat{D}(k)]^2}}.
\end{equation}
We split up the integral and bound:
\begin{equation}
    \begin{split}
        C  \intFourier{[-\pi, \pi]^d}{ \frac{\hat{h}_{r}(k)}{[1-\hat{D}(k)]}} &\le C_A \intFourier{\lvert k \rvert \le 1/r}{ \frac{1}{\lvert k \rvert^{\twa}}} + C_B r^{\twa} \intFourier{[-\pi,\pi]^d}{\hat{h}_r(k)}\\
        & \le C_A r^{\twa} + C_B r^{\twa}.
    \end{split}
\end{equation}
Here the first integral has been bounded in the same way as \refeq{pcA} and the second one in the same way as \refeq{pcB}.
Combining both bounds, we can conclude that
\begin{equation}
    S_{1,1} \ge p_c^2 \intFourier{[-\pi, \pi]^d}{ \htaupc(k)^2 \hat{h}_{r}(k)} - C r^{\twa}.
\end{equation}

For $S_{1,2}$ we need an upper bound.
Using \refeq{genEvent} and the tree-graph bound \refeq{TreeGraph}, we get
\begin{equation}
    \begin{split}
        S_{1,2} &=\sum_{y, \ulb \in \Zd} \sum_{e: \ule=\ulb} p_{c}^2 h_r(y)D(e) \E_0[\indi_{\{0 \conn \ule\}} \E_1 [\indi_{\{\ole \conn y \text{ through }\tCcal^e (0)\}}]]\\
                &\le \sum_{y, \ulb,v,v' \in \Zd} \sum_{e: \ule=\ulb} p_{c}^2 h_r (y) D(e) \taupc(v)\taupc(\ule-v)\taupc(v' -\ole)\taupc(v'-v) \taupc(y-v').
    \end{split}
\end{equation}
Observe that
\begin{equation}
    \taupc(v'-v)\left[\sum_{\ulb \in \Zd}\sum_{e:\ule=\ulb}\taupc(\ule-v) D(\ole-\ule) \taupc(v'-\ole)\right] \le T'_{p_c} (v'-v).
\end{equation}
This implies
\begin{equation}
    S_{1,2} \le p_{c}^2 \sum_{y, v,v' \in \Zd} h_r(y) \taupc(v)T'_{p_c} (v'-v)\taupc(y-v')= p_{c}^2(\taupc * T'_{p_c} * \taupc * h_r)(0).
\end{equation}
We rewrite the right-hand side in terms of its Fourier transform and apply \refeq{Tppcbd}:
\begin{equation}\label{e:S2bound}
    S_{1,2} \le  p_{c}^2 \intFourier{[-\pi, \pi]^d}{\htaupc(k)^2 \hat{T'}_{p_c} (k) \hat{h}_r (k)} \le p_c^2 C' \beta \intFourier{[-\pi, \pi]^d}{\htaupc(k)^2  \hat{h}_r (k)}.
\end{equation}

Finally, combining the bounds \refeq{S1bound}, \refeq{S2bound} and \refeq{S2bd}, we obtain, for $\beta$ small enough,
\begin{equation}
        \Eiic[|Q_r  \cap \Ccal(0)|] \ge S_{1,1}-S_{1,2}-S_2
                                    \ge p_c^2 (1- C' \beta) \intFourier{[-\pi, \pi]^d}{\hat{\tau}_{p_c}(k)^2 \hat{h}_r (k)} - C r^{\twa} \ge c''(d) r^{2\twa}
\end{equation}
for some constant $c''(d)$. The last inequality follows from a similar bound as \refeq{pcLowB}.
This completes the proof of Theorem \ref{ExpectationBounds}. \qed

\section{A lower bound on the long-range one-arm probability: proof of Theorem \ref{th:OAprob}}\label{sec:OAprob}

In this section we restrict ourselves to models of long-range spread-out percolation only.
\proof[Proof of Theorem \ref{th:OAprob}]
We start by proving $\Ppc(0 \conn \dqr) \ge c/r^{\alpha/2}$.
Let $\Ccal_r (0)$ be the $r$\emph{-truncated cluster} of $0$, that is, the percolation cluster of $0$ generated by using the edge probability
\[
    D_r (x) = D(x) \indi_{\{|x| \le r\}}
\]
instead of $D(x)$. Note that there exists $\zeta >0$ such that
\begin{equation}\label{e:totalDr}
    \sum_{x \in \Zd} D_r(x) = \sum_{x \in Q_r} D(x) \le 1 - \zeta  r^{-\alpha}
\end{equation}
when $r$ is sufficiently large.

One way for a path from $0$ to reach $\dqr$ is if $\Ccal_r (0)$ is at least of size $k$ (we will fix the value of $k$ later), and at least one of the vertices, say $v$, of $\Ccal_r(0)$ is an endpoint of an open edge $e$ that has length at least $2r$. Then either $v \in \dqr$ and so there exists a path, or, perhaps more likely, $v \in Q_r$, which would imply that the other endpoint of $e$ is in $\dqr$.
Hence, we have
\begin{equation}
    \Ppc(0 \conn \dqr) \ge \Ppc( \lvert \Ccal_r (0) \rvert \ge k, \exists e=(\ule,\ole) \text{ such that } \ule \in \Ccal_r (0), \lvert e \rvert > 2r).
\end{equation}
Edge probabilities are translation invariant and independent, and there are at least $k$ vertices in $\Ccal_r(0)$, so we have a lower bound on the right-hand side,
\begin{multline}
    \Ppc(\lvert \Ccal_r (0) \rvert \ge k)\left(1- \Ppc(\nexists\, \ole \in Q_{2r}^c \text{ such that } \{0, \ole\} \text{ is open})^k\right)\\ = \Ppc(\lvert \Ccal_r (0) \rvert \ge k)\left(1- \left(p_c \sum_{\ole \in Q_{2r}} D(\ole)\right)^k\right) \ge \Ppc(\lvert \Ccal_r (0) \rvert \ge k)\left(1- \left(1 - \frac{\zeta}{(2 r)^{\alpha}}\right)^k\right)
\end{multline}
where we used \refeq{totalDr} and the fact that $p_c \ge 1$ in the last step. When $k < (2 r)^{\alpha}/\zeta$, we may bound
\begin{equation}\label{e:polybd}
    \left(1 - \frac{\zeta }{(2 r)^{\alpha}} \right)^k \le 1 - \frac{\zeta  k}{2(2 r)^{\alpha}}.
\end{equation}
Thus,
\begin{equation}\label{e:firststep}
    \Ppc(0 \conn \dqr) \ge \frac{\zeta  k }{(2 r)^{\alpha}} \Ppc( \lvert \Ccal_r (0) \rvert \ge k) \qquad \text{ for all }k < \frac{(2 r)^{\alpha}}{\zeta}.
\end{equation}
We are left to prove a lower bound on $\Ppc(\lvert \Ccal_r (0) \rvert \ge k)$.

Combining results of \cite{BarAiz91} and \cite{HeyHofSak08}, we may conclude that there exists a constants $C_1 \ge c_1 >0$ such that
\begin{equation}\label{e:BAsize}
   \frac{c_1}{\sqrt{s}} \le  \Ppc(\lvert \Ccal(0) \rvert \ge s) \le \frac{C_1}{\sqrt{s}}
\end{equation}
holds for long-range percolation when $d > d_c$. Furthermore, we have
\begin{equation}
    \Ppc(\lvert \Ccal_r(0) \rvert \ge k) \ge \Ppc(\lvert \Ccal(0) \rvert \ge k) - \Ppc( \lvert \Ccal(0) \rvert \ge k, \lvert \Ccal_r (0) \rvert < k).
\end{equation}
To bound the first term on the right-hand side we use \refeq{BAsize}. For the second term we need an upper bound.  Given that $\lvert \Ccal_r(0)\rvert <k$, for $\rvert\Ccal(0)\lvert \ge k$ to hold, there needs to exist at least one open edge that is longer than $r$ with at least one endpoint in $\Ccal_r (0)$. Thus,
\begin{equation}
   \Ppc(\lvert \Ccal(0) \rvert \ge k, \lvert \Ccal_r (0) \rvert < k) \le \Ppc(\lvert \Ccal_r (0) \rvert < k,  \exists  e=\{\ule, \ole\} \text{ with } \ule \in \Ccal_r (0) \text{ s.t. } \lvert e \rvert > r , e \text{ open}).
\end{equation}
The probability of having such an edge only depends on $\lvert \Ccal_r(0)\rvert$, the number of possible allowed endpoints for this edge. Hence, we may condition on the size of $\Ccal_r(0)$ and use translation invariance and independence of edges for an upper bound:
\begin{equation}\label{e:secondstep}
    \begin{split}
    \Ppc\big(\exists  e=\{\ule, \ole\} \text{ with } & \ule \in \Ccal_r (0) \text{ s.t. } \lvert e \rvert > r , e \text{ open} \big\vert\, \lvert \Ccal_r (0) \rvert < k\big) \Ppc(\lvert \Ccal_r (0) \rvert < k)\\
    & \le \sum_{s=1}^{k-1} s \; \Ppc(\exists v \in Q_r^c \text{ s.t. } \{0,v\} \text{ open}) \Ppc(\lvert \Ccal_r (0) \rvert = s)\\
    &\le \frac{\zeta}{r^{\alpha}} \sum_{s=1}^{k-1} \Ppc( \lvert \Ccal_r(0) \rvert > s) \\
    &\le \frac{\zeta}{r^{\alpha}} \sum_{s=1}^{k-1} \Ppc( \lvert \Ccal(0) \rvert > s)\\
    &\le \frac{\zeta}{r^{\alpha}} \sum_{s=1}^{k-1} \frac{C_1}{\sqrt{s}} \le \frac{C_2 \sqrt{k}}{r^{\alpha}}
    \end{split}
\end{equation}
where we used \refeq{BAsize} in the second to last step.
Applying the above bound to \refeq{firststep} with $k = \vep^2 r^{\alpha}$ and some suitably small constant $\vep$ thus yields
\begin{equation}\label{e:firstbound}
    \Ppc(0 \conn \dqr) \ge \frac{\zeta  \vep^2 r^{\alpha} }{(2 r)^{\alpha}} \Ppc( \lvert \Ccal_r (0) \rvert \ge \vep^2 r^{\alpha}) \ge \frac{\zeta  \vep^2}{2^{\alpha}}\left(\frac{c_1}{\vep r^{\alpha/2}} - \frac{C_2 \vep}{r^{\alpha/2}}\right) \ge \frac{c'}{r^{\alpha/2}},
\end{equation}
completing the proof for $\alpha \in (0,4]$.

To prove the theorem for $\alpha > 4$, that is, to establish $\Ppc(0 \conn \dqr) \ge c / r^2$, we employ the second moment method. Fix $n$ large, and define
\begin{equation}
    N_{r,nr} = \# \{x : x \in Q_{nr} \setminus Q_r \text{ and }0 \conn x\}.
\end{equation}
Then,
\begin{equation}
    \Ppc (0 \conn \dqr) \ge \Ppc( N_{r,nr} \ge 1).
\end{equation}
By the second moment method, we have
\begin{equation}
    \Ppc(N_{r,nr} \ge 1) \ge \frac{\Epc[N_{r,nr}]^2}{\Epc[N_{r,nr}^2]}.
\end{equation}
We can write
\begin{equation}
    N_{r,nr} = \lvert Q_{nr} \cap \Ccal(0)\rvert - \lvert Q_r \cap \Ccal(0)\rvert.
\end{equation}
By Theorem \ref{ExpectationBounds}, when $n$ is large enough,
\begin{equation}
    \Epc[N_{r,nr}] = \Epc[\lvert Q_{nr} \cap \Ccal(0)\rvert]- \Epc[\lvert Q_r \cap \Ccal(0)\rvert]\ge c_3 (nr)^{\twa} - C_4 r^{\twa} \ge c_5 r^{\twa}.
\end{equation}

We can write $N_{r,nr}^2$ as
\begin{equation}
    N_{r,nr}^2 = \# \{\text{pairs }x,y : x,y \in Q_{nr}\setminus Q_r \text{ and } 0 \conn x, 0 \conn y\}.
\end{equation}
Obviously,
\begin{equation}
    N_{r,nr}^2 \le \#\{\text{triplets }x,y,z : x,y \in Q_{nr}\setminus Q_r, z \in \Zd \text{ and } \{0 \conn z\}\circ\{z \conn x\}\circ\{z \conn y\}\}
\end{equation}
where ``$\circ$'' denotes disjoint occurrence. Using the BK-inequality \cite{BerKes85} and techniques similar to those used in the proof of Theorem \ref{ExpectationBounds}, we can show
\begin{equation}
    \Epc[N_{r,nr}^2] \le \sum_{x,y \in Q_{nr}} \sum_{z \in \Zd} \taupc(z) \taupc(x-z) \taupc(y-z) \le C_{6} r^{3\twa}.
\end{equation}
Hence, it follows that
\begin{equation}\label{e:secondbound}
    \Ppc(0 \conn \dqr) \ge \frac{c_5 r^{2 \twa}}{C_{6} r^{3 \twa}} \ge \frac{c''}{r^{2}}.
\end{equation}

Finally, we combine the bounds \refeq{firstbound} and \refeq{secondbound}. This yields
\begin{equation}
    \Ppc(0 \conn \dqr) \ge \max \left\{\frac{c'}{r^{\alpha/2}}, \frac{c''}{r^2} \right\} \ge \frac{c}{r^{(4 \wedge \alpha)/2}},
\end{equation}
completing the proof. \qed 

\section{Bounds on lace expansion coefficients: proof of Proposition \ref{prop:TermBounds}}\label{app:Prop}
In this section and the next we prove Proposition \ref{prop:TermBounds}. We start by showing that the functions $\xi^{(n)}$ and $\gamma^{(n)}$ can both be bounded in terms of one-arm probabilities, $\pi^{(n)}$, and another function, $\phi^{(n)}$. Then we bound the complex expressions $\pi^{(n)}$ and $\phi^{(n)}$ in terms of simpler two-point functions. These bounds are known as \emph{diagrammatic estimates}.

Using the diagrammatic estimates we are able to obtain the bounds required to prove Proposition \ref{prop:TermBounds}, but it involves a lot of machinery to do so.

In the case of item Proposition \ref{prop:TermBounds}(i), this is mainly due to the fact that the function $\phi^{(n)}$ has not appeared in any other lace expansion (though a similar function is considered for oriented percolation in \cite{HofHolSla07b}), so there is little to fall back on.

In the case of Proposition \ref{prop:TermBounds}(ii), the reason for the difficulties is more fundamental. The bound that we require is quite strong while our knowledge of the two-point functions is relatively minimal and mainly consists of its properties in Fourier space. Significant effort is needed to evaluate these functions in Fourier space without sacrificing too much accuracy in the bounds.

In the course of the proof of Proposition \ref{prop:TermBounds}(ii) we introduce a method for obtaining lace expansion diagrams \emph{in Fourier space}. This construction makes use of ideas from graph theory, and in principle applies to any lace expansion whose terms can be bounded by `planar' diagrams (e.g. self-avoiding walk, lattice animals and lattice trees).
Moreover, the Fourier space diagrams have a simple combinatorial structure and are fairly easy to bound.

\subsection{Proof of Proposition \ref{prop:TermBounds} (i)}
In this section we prove Proposition \ref{prop:TermBounds}(i) subject to Proposition \ref{prop:TermBounds} (ii), which we prove in the next section, and subject to Lemma \ref{lem:Phibd}, which is stated further along in the section and proved in the final subsection. The techniques that we use are similar to those used in \cite{HofHolSla07b}, but much less refined, as we only need an upper bound.

\proof[Proof of Proposition \ref{prop:TermBounds}(i) subject to Proposition \ref{prop:TermBounds}(ii) and Lemma \ref{lem:Phibd}.]
Recall definitions \refeq{Expansion1} -- \refeq{Expansion4}.
We start with the following observation: $F$ is a cylinder event restricted to a finite box $Q_m$. Hence, there exists a finite positive constant $C_m$ that only depends on $Q_m$, such that
\begin{equation}
    \Epc[\indi_{F \cap \{0 \conn \dqr, Q_m \Conn \dqr\}} ] \le C_{m} \Epc [\indi_{E''(0,r;Q_m)}]\text{ and }\Epc[\indi_{F \cap \{0 \conn x, Q_m \Conn x\}} ] \le C_{m} \Epc [\indi_{E'(0,x;Q_m)}].
\end{equation}
In light of the above and the fact that we only need upper bounds for the proof of Proposition \ref{prop:TermBounds}, we can bound
\begin{equation}\label{e:xibd0}
    \xi^{(0)}(r;F) \le C_{m}\; \xi^{(0)}_{m}(r) \equiv C_{m}\; \Ppc(E''(0,r;Q_m)).
\end{equation}

Similarly, for $n \ge 1$,
\begin{equation}\label{e:xibd00}
    \begin{split}
    \xi^{(n)}(r;F) \le C_{m}^{(n)}\; \xi^{(n)}_{m}(r) &\equiv C_{m} \sum_{(u_0, v_0)\in\Be} p_{u_0, v_0} \dotsm \sum_{(u_{n-1}, v_{n-1})\in \Ee} p_{u_{n-1} v_{n-1}} \E_0 \bigg[\indi_{E'(0,u_0;Q_m)}\\
                        &\quad \times \E_1\Big[\indi_1 \E_2\big[ \indi_2 \dotsm \E_{n-1}[\indi_{n-1} \E_n[ \indi_{E''(v_{n-1}, r;\tCcal_{n-1})}]]\dotsm \big]\Big]\bigg].
    \end{split}
\end{equation}
By making similar replacements in \refeq{GammaZero}, \refeq{Expansion6}, \refeq{PiZero}, \refeq{Expansion3} and \refeq{Expansion4} we can also define for $n \ge 0$ the upper bounds
\begin{equation}\label{e:grbd0}
    \gamma^{(n)}(r;F) \le C_{m} \gamma^{(n)}_{m}(r), \qquad \pi^{(n)}(x,r;F) \le  C_{m} \pi^{(n)}_{m}(x,r), \quad \text{and}\quad \psi^{(n)}(x,r;F) \le C_{m}\psi^{(n)}_{m}(x,r).
\end{equation}
Since upper bounds on the functions on the right-hand sides of \refeq{xibd0}, \refeq{xibd00} and \refeq{grbd0} imply upper bounds on their respective left-hand sides, the influence of the event $F$ will only play a role through $Q_m$. For this reason we will from here on consider only the $Q_m$-dependent right-hand sides.

We start by showing that $\xi^{(n)}_{m}(r)$ and $\gamma^{(n)}_{m}(r)$ can be bounded as follows:
\begin{eqnarray}
    \label{e:xibd1} \xi^{(n)}_{m}(r) &\le& \sum_{x \in Q_r} \theta^{(n)}_{m}(x,r) \Ppc(x \conn \dqr) + \rho^{(n)}_{m}(r);\\
    \label{e:gammabd1} \gamma^{(n)}_{m}(r) &\le& \sum_{x \in Q_r} \theta^{(n)}_{m}(x,r) \Ppc(x \conn \dqr),
\end{eqnarray}
for the function $\theta^{(n)}_{m}$ defined below in \refeq{theta0} and \refeq{thetan} and $\rho^{(n)}_{m}$ defined below in \refeq{rho0} and \refeq{rhon}.

Then we show that $\rho^{(n)}_{m}$ can be bounded by
\begin{equation}\label{e:rhotopi}
    \rho^{(n)}_{m} (r)\le \sum_{x \in \dqr} \pi^{(n)}_{m} (x,r)
\end{equation}
and that $\theta^{(n)}_{m}$ can be bounded further by
\begin{equation}\label{e:thetatophi}
    \theta^{(n)}_{m}(x,r) \le \sum_{z \in Q_r} \phi^{(n)}_{m}(z,x,r)\Ppc(z \conn \dqr)
\end{equation}
for the function $\phi^{(n)}_{m}$ defined below in \refeq{phi0} and \refeq{phin}.

After showing that such bounds exist, we obtain diagrammatic bounds on $\theta^{(n)}_{m}$ and $\phi^{(n)}_{m}$ that suffice to prove Proposition \ref{prop:TermBounds}(i) (subject to Proposition \ref{prop:TermBounds}(ii)).

Define
\begin{equation}\label{e:theta0}
    \theta^{(0)}_{m}(x,r) = \sum_{y \in Q_r} p_{y,x} \Epc[\indi_{\{\{0 \conn y, Q_m \Conn y, 0 \conn \dqr\} \text{ on }\tCcal^{(y,x)} (0)\}}],
\end{equation}
then \refeq{gammabd1} for $n=0$ follows immediately from \refeq{GammaZero}, \refeq{grbd0} and the simple fact that
\begin{equation}
    \Ppc^A(x \conn \dqr) \le \Ppc (x \conn \dqr).
\end{equation}
For $n \ge 1$, define,
\begin{equation}\label{e:thetan}
    \begin{split}
        \theta^{(n)}_{m}(x, r) &= \sum_{y \in Q_r} p_{y,x} \sum_{(u_0, v_0)\in \Be} p_{u_0, v_0} \dotsm \sum_{(u_{n-1}, v_{n-1})\in \Be} p_{u_{n-1} v_{n-1}}\E_0 \bigg[\indi_{E'(0,u_0;Q_m)}\\
                        &\quad \times \E_1\Big[\indi_1 \E_2\big[ \indi_2 \dotsm \E_{n-1}[\indi_{\{E'(v_{n-1}, y; \tCcal_{n-1})\cap\{v_{n-1} \conn \dqr\} \text{ on }\tCcal^{(y,x)} (v_{n-1})\}}]\dotsm \big]\Big]\bigg].
    \end{split}
\end{equation}
Combined with \refeq{Expansion6}, \refeq{thetan} gives \refeq{gammabd1} for $n \ge 1$.

Although the purpose of the functions $\theta^{(n)}$ is to bound the probability of events $E'$ that are restricted to be connected to $\dqr$, it will come in handy later on to use that the bound
\begin{equation}\label{e:thetatopi}
    \theta^{(n)}_{m} (x, r) \le \pi^{(n)}_{m}(x,r)
\end{equation}
still holds, since
\begin{equation}
    E'(v_{n-1}, y; \tCcal_{n-1})\cap\{v_{n-1} \conn \dqr\} \subset E'(v_{n-1}, y; \tCcal_{n-1}).
\end{equation}

To show \refeq{xibd1}, we need to do a bit more work. In a similar fashion as in \cite{HofHolSla07b}, we define the set
\begin{equation}
    \Pcal_A = \bigl\{\text{edges }b \vert \text{ the event }E'(v, \ulb;A) \cap \{b \text{ open}\}\cap\{\olb \conn \dqr \text{ off }\tCcal^b (v)\} \text{ occurs}\bigr\}.
\end{equation}
In words, $\Pcal_A$ is the (unordered) set of \emph{cutting edges}, i.e., edges in $\Pcal_A$ have the property that they are open and they are the first pivotal edge after $A$ for at least one connection from $v$ to $\dqr$. (This means that these edges are not necessarily pivotal for all connections from $v$ to $\dqr$.)

We can decompose the event $E''(v,r;A)$ according to the size of $\Pcal_{A}$:
\begin{eqnarray}
       \nonumber \Ppc(E''(v,r;A))&=& \Ppc(E''(v,r;A)\cap\{\Pcal_{A} = \varnothing \}) + \sum_{l=1}^{\infty} \frac{1}{l} \sum_{b \in \Ee} \Ppc(E''(v,r;A) \cap \{b \in \Pcal_{A}\} \cap \{\lvert \Pcal_{A} \rvert = l\})\\
       \label{e:decomp1}                     &=& \frac{1}{2} \sum_{b \in \Ee} \Ppc(E''(v,r;A)\cap \{b \in \Pcal_{A}\}) + \rho^{(0)}(v,r;A)
\end{eqnarray}
where
\begin{equation}\label{e:rho0}
        \rho^{(0)} (v,r;A) = \Ppc(E''(v,r;A)\cap\{\Pcal_{A} = \varnothing \}) + \sum_{l=1}^{\infty}\left(\frac{1}{l}-\frac12 \right) \sum_{b \in \Ee} \Ppc(E''(v,r;A) \cap \{b \in \Pcal_{A}\} \cap \{\lvert \Pcal_{A} \rvert = l\}).
\end{equation}
Define $\rho^{(0)}_{m}(r) = \rho^{(0)}(0,r;Q_m)$, then it follows that
\begin{equation}\label{e:xibd2}
    \xi^{(0)}_{m}(r) = \frac{1}{2} \sum_{b \in \Ee} \Ppc(E''(0,r;Q_m)\cap \{b \in \Pcal_{Q_m}\}) + \rho^{(0)}_{m}(r).
\end{equation}
Similarly, by replacing the final expectation in \refeq{Expansion2} by \refeq{decomp1}, we can isolate a term $\rho^{(n)}_{m}(r)$ from $\xi^{(n)}_{m}(r)$:
\begin{equation}\label{e:rhon}
    \begin{split}
        \rho^{(n)}_{m}(r) &= \sum_{(u_0, v_0)\in \Be} p_{u_0, v_0} \dotsm \sum_{(u_{n}, v_{n})\in \Ee} p_{u_{n} v_{n}}\E_0 \bigg[\indi_{E'(0,u_0;Q_m)}\\
                        &\quad \times \E_1\Big[\indi_1 \E_2\big[ \indi_2 \dotsm \E_{n-1}[\rho^{(0)}(v_{n-1}, r; \tCcal^{(u_n,v_n)}(v_{n-1}))]\dotsm \big]\Big]\bigg],
    \end{split}
\end{equation}

From \cite[Proposition 4.3]{HofHolSla07b} we have the following useful identity: for $A \subseteq \Zd$, $v \in \Zd$, $r \ge 1$ and $b \in \Be$,
\begin{equation}\label{e:decomp2}
    E''(v,r;A) \cap \{b \in \Pcal_A\} = \{E'(v,\ulb;A) \cap \{v \Aconn \dqr\} \text{ on }\tCcal^b (v) \}\cap\{b \text{ open}\} \cap \{\ulb \conn \dqr \text{ off }\tCcal^b (v)\}.
\end{equation}
This equality is proved in \cite{HofHolSla07b} for oriented percolation, but the proof is easily adapted to the unoriented case.

Applying \refeq{decomp2} with $A =Q_m$ and $v =0$ to \refeq{xibd2} and using
\begin{equation}\label{e:subsetsimple}
    \{x \Aconn \dqr\} \subseteq \{x \conn \dqr\}
\end{equation}
and the Factorization Lemma yields \refeq{xibd1} for $n=0$. Applying \refeq{decomp2} with $A = \tCcal^b (v_{n-1})$ and $v = v_{n-1}$ to \refeq{xibd2} and again using \refeq{subsetsimple} and the Factorization Lemma yields \refeq{xibd1} for $n \ge 1$.

Now we show \refeq{rhotopi}. Observe that the sum in \refeq{rho0} only has positive contributions when $l=0,1$, so we do not have to consider the terms $l \ge 2$ for an upper bound.
Therefore,
\begin{equation}
    \rho^{(0)}(v,r;A) \le \Ppc(E''(v,r;A) \cap \{\lvert \Pcal_A \rvert \le 1\}).
\end{equation}
From \cite[Proposition 4.6]{HofHolSla07b} we have
\begin{equation}\label{e:decomp3}
    E''(v,r;A) \cap \{\lvert \Pcal_A \rvert \le 1\} \subseteq \bigcup_{x \in \dqr} E'(v,x;A).
\end{equation}
Again, \refeq{decomp3} is proved in \cite{HofHolSla07b} for oriented percolation, and again the proof is straightforwardly adapted to the unoriented case.

Applying \refeq{decomp3} with $A = Q_m$ and $v=0$ to \refeq{rho0} and applying \refeq{decomp3} with $A = \tCcal^b (v_{n-1})$ and $v=v_{n-1}$ to \refeq{rhon} yields \refeq{rhotopi} for $n\ge0$.

Next is the bound \refeq{thetatophi}.
From the tree-graph inequality \cite{AizNew84} and the definition of $E'$ it follows that
\begin{equation}
   E'(v,x;A) \cap \{v \conn \dqr\} \subseteq \bigcup_{z \in Q_r} E'(v, x; A) \cap \{\{v \conn z\} \circ \{z \conn x\} \circ \{z \conn \dqr\}\}.
\end{equation}
Hence, by the BK-inequality,
\begin{equation}\label{e:thetaBK}
    \Ppc(E'(v,x;A)\cap\{v \conn \dqr\}) \le \sum_{z \in Q_r} \Ppc(E'(v,x;A) \cap \{\{v \conn z\} \circ \{z \conn x\}\}) \Ppc(z \conn \dqr).
\end{equation}
Define
\begin{equation}\label{e:Hprime}
    H'(v,z,x;A) = E'(v,x;A) \cap \{\{v \conn z\} \circ \{z \conn x\}\},
\end{equation}
and define
\begin{equation}\label{e:phi0}
    \phi^{(0)}_{m}(z,x,r)  = \Ppc(H'(0,z,x;Q_m))
\end{equation}
(as with $\pi^{(0)}_{m}(x,r)$, this function is independent of $r$, but we write $r$ anyway for consistency). Also define, for $n\ge1$,
\begin{equation}\label{e:phin}
    \begin{split}
    \phi^{(n)}_{m}(z,x,r)  &= \sum_{(u_0, v_0)\in \Be} p_{u_0, v_0} \dotsm \sum_{(u_{n-1}, v_{n-1})\in \Ee} p_{u_{n-1} v_{n-1}}\E_0 \bigg[\indi_{E'(0,u_0;Q_m)}\\
                         &\quad \times \E_1\Big[\indi_1 \E_2\big[ \indi_2 \dotsm \E_{n-1}[\indi_{\{H'(v_{n-1},z, x; \tCcal_{n-1})\}}]\dotsm \big]\Big]\bigg].
    \end{split}
\end{equation}
Now it follows from \refeq{theta0}, \refeq{thetan}, \refeq{thetaBK} and \refeq{Hprime} that \refeq{thetatophi} holds.

For future use we define
\begin{equation} \label{e:PTPdef} \widetilde{\Pi}_{m}(x,r) = \sum_{n=0}^{\infty} \pi_{m}^{(n)} (x,r), \qquad
 \widetilde{\Theta}_{m}(x,r) = \sum_{n=0}^{\infty} \theta^{(n)}_{m} (x,r)\qquad \text{and}\qquad \widetilde{\Phi}_{m}(x,y,r) = \sum_{n=0}^{\infty} \phi^{(n)}_{m}(x,y,r).
\end{equation}

\pagebreak
Before we proceed we state the following lemma:
\bl\label{lem:Phibd}
    Under the same assumptions as Theorems \ref{QiicTh1} and \ref{QiicTh2} and for the same choice of $\delta>0$ as in Proposition \ref{prop:TermBounds}(ii), for $L \ge L_0$, all $r\in \N$, and $i=1,2,3$, there exists a constant $C_i = C_i(m,L,d,\alpha,\delta)$, such that
        \begin{eqnarray}
        \label{e:Phinodelta} \sum_{x,y \in \Zd} \widetilde{\Phi}_{m} (x,y,r) & \le & C_1;\\
        \label{e:Phidelta} \sum_{x,y \in \Zd} \lvert x - y \rvert^{\delta} \widetilde{\Phi}_{m}(x,y,r) &\le& C_2;\\
        \label{e:Pidelta} \sum_{x \in \Zd} \lvert x \rvert^{\twa+\delta} \widetilde{\Pi}_{m}(x,r) &\le& C_3.
        \end{eqnarray}
\el
We do not prove Lemma \ref{lem:Phibd} since it can be proved in a similar way as Proposition \ref{prop:TermBounds}(ii). At the end of Section \ref{sec:proofprop252} we do briefly discuss this proof for Lemma \ref{lem:Phibd}.

\medskip
From \refeq{grbd0} -- \refeq{thetatophi}, \refeq{thetatopi} and \refeq{PTPdef} it follows that
\begin{equation}
    \Xi(r;F) \le C_{m, \Xi} \sum_{x \in Q_r} \widetilde{\Theta}_{m}(x,r) \Ppc(x \conn \dqr) + C_{m, \Xi} \sum_{x \in \dqr} \widetilde{\Pi}_{m}(x,r)
\end{equation}
and
\begin{equation}
    \Gamma(r;F) \le C_{m, \Gamma} \sum_{x \in Q_r} \widetilde{\Theta}_{m}(x,r) \Ppc(x \conn \dqr),
\end{equation}
for constants $C_{m, \Xi}$ and $C_{m, \Gamma}$ that may depend on $F$.
Hence, Proposition \ref{prop:TermBounds}(i) is proved once we show
\begin{equation}\label{e:Thetabd}
    \sum_{x \in Q_r} \widetilde{\Theta}_{m}(x,r) \Ppc(x \conn \dqr) \le \frac{C_4}{r^{1/\rho+\delta}}
\end{equation}
for a constant $C_4$ that may depend on $L, d, \alpha$ and $\delta$, and
\begin{equation}\label{e:Pibdout}
    \sum_{x \in \dqr} \widetilde{\Pi}_{m} (x,r) \le \frac{C_3}{r^{\twa + \delta}}.
\end{equation}

That \refeq{Pibdout} holds follows immediately from \refeq{Pidelta}: $x \in \dqr$ implies $\lvert x \rvert / r >1$, so
\begin{equation}\label{e:Pibdout2}
    \sum_{x \in \dqr} \widetilde{\Pi}_{m}(x,r) \le \sum_{x \in \dqr} \frac{\lvert x \rvert^{\twa + \delta}}{r^{\twa+\delta}} \widetilde{\Pi}_{m}(x,r) \le \frac{C_3}{r^{\twa+\delta}}.
\end{equation}

To bound \refeq{Thetabd} we introduce the following notation: for $a,b \in \N$ and $a > b$,
\begin{equation}
    Q_{a,b} = Q_a \setminus Q_b.
\end{equation}

The sum on the left-hand side of \refeq{Thetabd} can be split into the contributions of $x \in Q_{r/4}$ and those of $x \in Q_{r,r/4}$:
\begin{equation}
    \sum_{x \in Q_r} \widetilde{\Theta}_{m}(x,r) \Ppc(x \conn \dqr)= \sum_{x \in Q_{r/4}} \widetilde{\Theta}_{m}(x,r) \Ppc(x \conn \dqr)+\sum_{x \in Q_{r,r/4}} \widetilde{\Theta}_{m}(x,r) \Ppc(x \conn \dqr).
\end{equation}
The second term can be bounded using \refeq{thetatopi} and \refeq{Pidelta}:
\begin{equation}\label{e:Thetabd1}
    \sum_{x \in Q_{r,r/4}} \widetilde{\Theta}_{m}(x,r) \Ppc(x \conn \dqr) \le \sum_{x \in Q_{r/4}^{c}} \widetilde{\Pi}_{m}(x,r) \le \frac{C_3}{r^{\twa + \delta}}.
\end{equation}
To bound the first term we use \refeq{thetatophi}:
\begin{equation}\label{e:ThetabdPhi}
    \sum_{x \in Q_{r/4}} \widetilde{\Theta}_{m}(x,r) \Ppc(x \conn \dqr) \le \qsum{r/4}{r} \widetilde{\Phi}_{m}(x,y,r) \Ppc(x \conn \dqr) \Ppc(y \conn \dqr).
\end{equation}
For convenience, denote
\begin{equation}
    S_{m}^{\Phi}(x,y,r) = \widetilde{\Phi}_{m}(x,y,r)\Ppc(x \conn \dqr) \Ppc(y \conn \dqr).
\end{equation}
The right-hand side of \refeq{ThetabdPhi} can again be split into two parts, now according to whether $y \in Q_{r/2}$ or $y \in Q_{r,r/2}$:
\begin{equation}\label{e:sphisplit}
    \qsum{r/4}{r} S_{m}^{\Phi}(x,y,r) = \qsum{r/4}{r/2} S_{m}^{\Phi}(x,y,r)+ \qsum{r/4}{r,r/2} S_{m}^{\Phi}(x,y,r).
\end{equation}
To bound the first term on the right-hand side, we note for $a,b \in \N$, $a > b$ and $x \in Q_{b}$, by \refeq{LRPArm} we have
\begin{equation}\label{e:onearmab}
    \Ppc(x \conn Q_{a}^{c}) \le \Ppc(0 \conn Q_{a-b}^{c}) \le \frac{C}{(a-b)^{1/\rho}}
\end{equation}
so that, by Lemma \ref{lem:Phibd},
\begin{equation}\label{e:Thetabd2}
    \qsum{r/4}{r/2} S_{m}^{\Phi}(x,y,r) \le \frac{C}{r^{2/\rho}} \qsum{r/4}{r/2} \widetilde{\Phi}_{m}(x,y,r) \le \frac{C_5}{r^{2 /\rho}}
\end{equation}
for some constant $C_5$ that may depend on $m, L, d, \alpha$ and $\delta$.

Finally, the second term in \refeq{sphisplit} can also be bounded using \refeq{onearmab} and Lemma \ref{lem:Phibd}: since $x \in Q_{r/4}$ and $\lvert x - y \rvert > r/4$, we can bound
\begin{equation}\label{e:Thetabd3}
    \begin{split}
     \qsum{r/4}{r,r/2} S_{m}^{\Phi}(x,y,r) &\le \frac{C}{r^{1/\rho}} \qsum{r/4}{r,r/2} \widetilde{\Phi}_{m}(x,y,r)\\
        &\le \frac{C}{r^{1/\rho}} \qsum{r/4}{r,r/2} \frac{\lvert x - y \rvert^{\delta}}{r^{\delta}} \widetilde{\Phi}_{m}(x,y,r) \le \frac{C_6}{r^{1/\rho + \delta}}
    \end{split}
\end{equation}
for some constant $C_6$ that may depend on $m, L, d, \alpha$ and $\delta$.

Combining \refeq{Thetabd1}, \refeq{Thetabd2} and \refeq{Thetabd3} gives the desired bound \refeq{Thetabd2} and completes the proof. \qed

\subsection{The proof of Proposition \ref{prop:TermBounds}(iii)}
From the definition of $R^{(N)}(r;F)$ in \refeq{Expansion5} and of $\pi^{(n)}(x,r;F)$ in \refeq{Expansion3} it is easy to see that
\begin{equation}
    R^{(N)}(r;F) \le p_c \sum_{x\in \Zd} \pi^{(N-1)}(x,r;F).
\end{equation}
It is a simple consequence of \refeq{PTPdef} and \refeq{Pidelta} that $\lim_{N \to \infty} \sum_x \pi^{(N-1)}(x,r;F) =0$. Furthermore, for all $N\ge 1$, $R^{(N)}(r;F)\ge 0$ and $\pi^{(N-1)}(x,r;F) \ge 0$, so by dominated convergence,
\begin{equation}
    \lim_{N \to \infty} R^{(N)}(r;F) = \lim_{N \to \infty} p_c \sum_{x \in \Zd} \pi^{(N-1)} (x,r;F) = 0.
\end{equation} \qed

\subsection{Diagrammatic estimates}\label{app:diagbd}
In this subsection we derive diagrammatic estimates on the functions $\pi_{m}^{(n)}$ and $\phi_{m}^{(n)}$. We need them to prove Proposition \ref{prop:TermBounds}(ii) and Lemma \ref{lem:Phibd}. Our derivation is based on the derivation given in \cite{BorChaHofSlaSpe05b}.

We start with $\pi^{(0)}_{m}$ and $\phi^{(0)}_{m}$. From the definition of $E'$ in \refeq{Eprime} it is easy to see that
\begin{equation}
    E'(0,x;Q_m) \subseteq \bigcup_{w \in Q_m}\left(\{0 \conn x\}\circ\{w \conn x\}\right).
\end{equation}
Hence, by the BK-inequality,
\begin{equation}\label{e:pi0d1}
    \pi^{(0)}_{m}(x,r) = \Ppc(E'(0,x;Q_m)) \le \Ppc\left(\bigcup_{w \in Q_m}\left(\{0 \conn x\}\circ\{w \conn x\}\right)\right) \le \sum_{w \in Q_m}\taupc(x)\taupc(x-w).
\end{equation}
Similarly, from the definition of $H'$ in \refeq{Hprime} it follows that
\begin{equation}
    H'(0,x,y;Q_m) \subseteq \bigcup_{w \in Q_m} \left(\{0 \conn x\}\circ\{w \conn y\}\circ\{y \conn x\}\right).
\end{equation}
Therefore, by \refeq{phi0} and the BK-inequality,
\begin{equation}\label{e:phi0d1}
        \phi^{(0)}_{m} (x,y,r)    \le \Ppc\left(\bigcup_{w \in Q_m} \left(\{0 \conn x\}\circ\{w \conn y\}\circ\{y \conn x\}\right)\right)
                                \le \sum_{w \in Q_m} \taupc(x)\taupc(y-x)\taupc(y-w).
\end{equation}
Furthermore, since on both right-hand sides of \refeq{pi0d1} and \refeq{phi0d1} we sum $w$ over the finite ball $Q_m$, we can bound both by $Q_m$-independent functions:
\begin{equation}\label{e:pibarzero}
    \pi^{(0)}_{m}(x,r) \le C'_m \bar{\pi}^{(0)}(x) := C'_m \taupc(x)^2
\end{equation}
and
\begin{equation}\label{e:phibarzero}
    \phi^{(0)}_{m} (x,y,r) \le C'_m \bar{\phi}^{(0)}(x,y) := C'_m \taupc(x)\taupc(x-y)\taupc(y)
\end{equation}
where $C'_m$ is a constant given by
\begin{equation}
    C'_m = \max \left\{\max_{x \in \Zd} \frac{\sum_{w\in Q_m} \taupc(x-w)}{\taupc(x)}, \max_{x \in \Zd} \frac{\sum_{w \in Q_m} (\taupc * \taupc)(x-w)}{\taupc(x)^2}\right\} < \infty.
\end{equation}

Let $\Ppc^{(n)}$ denote the product measure of $n+1$ copies of critical percolation on $\Zd$. We write $A_i$ to signify that the event $A$ occurs on the $i$th copy. By Fubini's theorem and \refeq{Expansion3} and \refeq{phin}, for $n \ge 1$,
\begin{equation}
    \begin{split}
        \pi^{(n)}_{m}(x,r) =    & \sum_{(u_0,v_0)\in \Ee}p_{u_0 v_0} \dotsm \sum_{(u_{n-1},v_{n-1})\in \Ee} p_{u_{n-1}v_{n-1}} \Ppc^{(n)}\Biggl(E'(0, u_0; Q_m)_0 \\
                            &  \qquad \cap \left(\bigcap_{i=1}^{n-1} E'(v_{i-1}, u_i; \tCcal_{i-1})_i \right) \cap E'(v_{n-1}, x; \tCcal_{n-1})_n \Biggr)
    \end{split}
\end{equation}
and
\begin{equation}
    \begin{split}
        \phi^{(n)}_{m}(x,y,r) =    & \sum_{(u_0,v_0)\in \Ee}p_{u_0 v_0} \dotsm \sum_{(u_{n-1},v_{n-1})\in \Ee} p_{u_{n-1}v_{n-1}} \Ppc^{(n)}\Biggl(E'(0, u_0; Q_m)_0  \\
                            &  \qquad \cap \left(\bigcap_{i=1}^{n-1} E'(v_{i-1}, u_i; \tCcal_{i-1})_i \right) \cap H'(v_{n-1}, x,y; \tCcal_{n-1})_n \Biggr).
    \end{split}
\end{equation}

\begin{figure}
  \includegraphics[width=490pt]{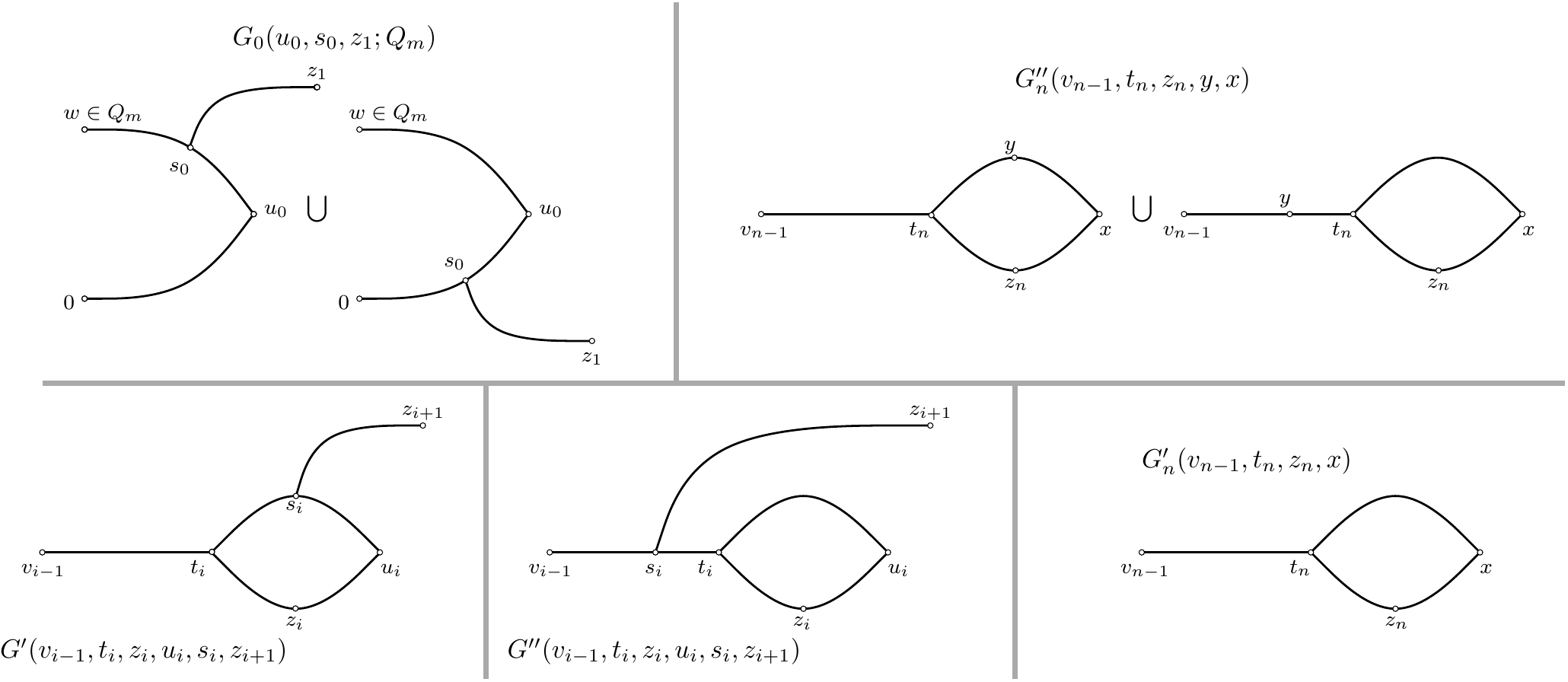}\\
  \caption{Depictions of the events $G_0$, $G'$, $G''$, $G'_n$ and $G''_n$.}\label{fig:GraphEvents}
\end{figure}

To estimate these functions, we define the events
\begin{eqnarray}
    G_0 (u_0, x_0, z_1 ; Q_m) &=& \left(\bigcup_{w \in Q_m} \{0 \conn u_0\}\circ\{w \conn s_0\}\circ\{s_0 \conn u_0\}\circ\{s_0 \conn z_1\}\right)\nonumber\\
    \label{e:Gzero}             & & \cup \left(\bigcup_{w \in Q_m} \{0 \conn s_0\}\circ\{s_0 \conn u_0\}\circ\{w \conn u_0\}\circ\{s_0 \conn z_1\}\right);\\
    G' (v_{i-1}, t_i, z_i, u_i, s_i, z_{i+1}) &=& \{v_{i-1} \conn t_i\}\circ\{t_i \conn z_i\}\circ\{t_i \conn s_i\} \circ \{z_i \conn u_i\} \circ \{s_i \conn u_i\} \circ \{s_i \conn z_{i+1}\};\quad\\
    G'' (v_{i-1}, t_i, z_i, u_i, s_i, z_{i+1}) &=& \{v_{i-1} \conn s_i\}\circ\{s_i \conn t_i\}\circ\{t_i \conn z_i\} \circ \{t_i \conn u_i\} \circ \{z_i \conn u_i\} \circ \{s_i \conn z_{i+1}\};\quad\\
    G (v_{i-1}, t_i, z_i, u_i, s_i, z_{i+1}) &=& G' (v_{i-1}, t_i, z_i, u_i, s_i, z_{i+1}) \cup G'' (v_{i-1}, t_i, z_i, u_i, s_i, z_{i+1});\\
    G'_{n} (v_{n-1}, t_n, z_n, x) &=& \{v_{n-1} \conn t_n\} \circ \{t_n \conn z_n\} \circ\{t_n \conn x\} \circ \{z_n \conn x\};\\
    \label{e:Gdprimen} G''_{n} (v_{n-1}, t_n, z_n, y,x) &=& \bigl(\{v_{n-1} \conn t_n\} \circ \{t_n \conn z_n\}\circ \{z_n \conn x\}\circ\{t_n \conn y \} \circ\{y \conn x\}\bigr)\\
                        & & \cup \bigl(\{v_{n-1} \conn y\}\circ\{y \conn t_n\} \circ \{t_n \conn z_n\}\circ \{z_n \conn x\} \circ\{t_n \conn x \}\bigr)\nonumber.
\end{eqnarray}
See Figure \ref{fig:GraphEvents} for depictions of these events.
All the events above are constructed of disjointly occurring, increasing events, and hence the BK-inequality can be used to factorize their probabilities.

The events inside $\pi^{(n)}_{m}$ and $\phi^{(n)}_{m}$ can be contained in constructions of the events \refeq{Gzero} -- \refeq{Gdprimen}: by definitions \refeq{Eprime} and \refeq{Hprime},
\begin{equation}\label{e:Eprimesetn}
    E' (v_{n-1}, x; \tCcal_{n-1})_n \subset \bigcup_{z_n \in \tCcal_{n-1}} \bigcup_{t_n \in \Zd} G'_{n} (v_{n-1}, t_n, z_n, x)_n
\end{equation}
and
\begin{equation}\label{e:Hprimesetn}
    H'(v_{n-1}, x,y; \tCcal_{n-1})_n \subset \bigcup_{z_n \in \tCcal_{n-1}} \bigcup_{t_n \in \Zd} G''_{n} (v_{n-1}, t_n, z_n, y,x)_n.
\end{equation}
For $n \ge 2$ and $i \in \{1,\dotso, n-1\}$,
\begin{equation}\label{e:Eprimeseti}
    E'(v_{i-1}, u_i; \tCcal_{i-1}) \cap \{z_{i+1} \in \tCcal_i\} \subset \bigcup_{z_i \in \tCcal_{i-1}} \bigcup_{t_i, s_i \in \Zd} G(v_{i-1}, t_i, z_i, u_i, s_i, z_{i+1})_i.
\end{equation}

\begin{figure}
  \includegraphics[width=490pt]{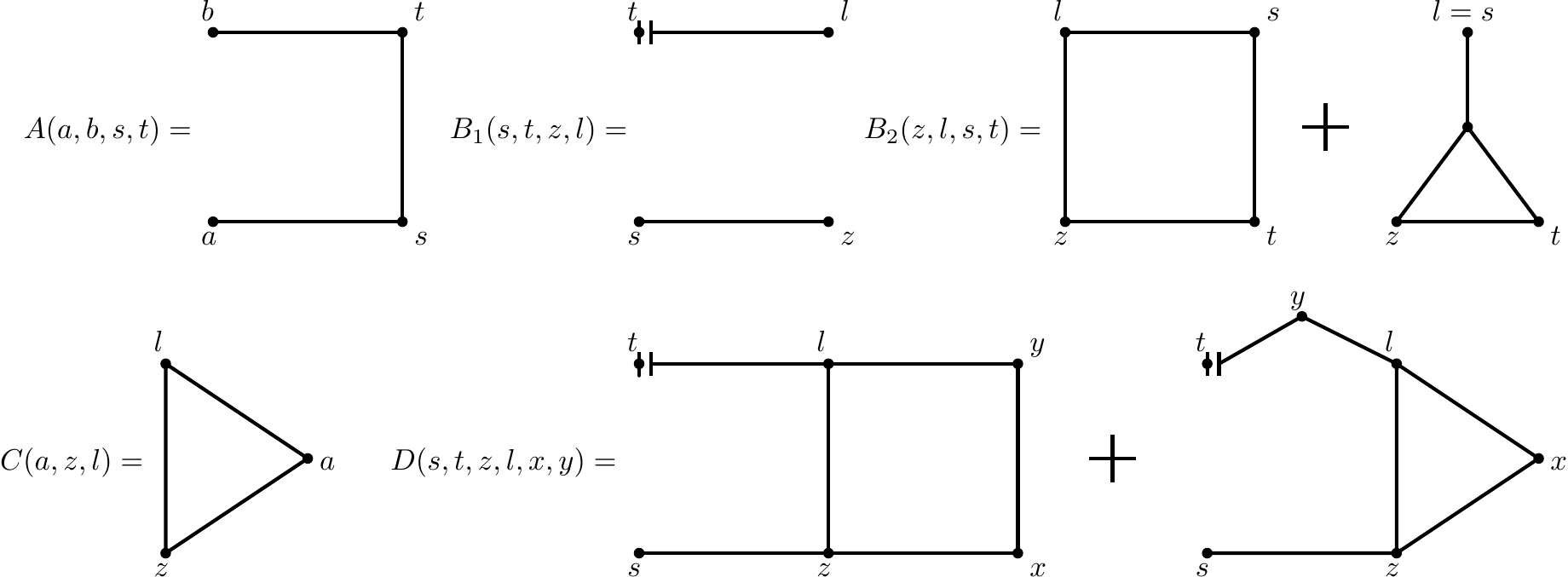}\\
  \caption{Diagrammatic representations of $A$, $B_1$, $B_2$, $C$ and $D$. Unbroken lines represent $\tau$'s, lines that start with a gap represent $\widetilde{\tau}$'s.}\label{fig:DiagExp}
\end{figure}

The relations \refeq{Eprimesetn} and \refeq{Eprimeseti} lead to
\begin{equation}
    \begin{split}
        E'(0,u_0 & ;Q_m)_0\cap \left(\bigcap_{i=1}^{n-1} E' (v_{i-1}, u_i; \tCcal_{i-1})_i\right)\cap E'(v_{n-1}, x ; \tCcal_{n-1})_n\\
        & \subset \bigcup_{\vec{t}, \vec{s}, \vec{z}} \Biggl(G_0(u_0, s_0, z_1 ;Q_m)_0 \cap \left(\bigcap_{i=1}^{n-1} G(v_{i-1}, t_i, z_i, u_i, s_i, z_{i+1})_i \right) \cap G'_n (v_{n-1}, t_n, z_n, x)_n \Biggr),
    \end{split}
\end{equation}
where $\vec{t}=(t_1, \dotso, t_n)$, $\vec{s} = (s_0, \dotso, s_{n-1})$ and $\vec{z}=(z_1, \dotso, z_n)$, and all elements are allowed to take values in $\Zd$.
The relations \refeq{Hprimesetn} and \refeq{Eprimeseti} lead to
\begin{equation}
    \begin{split}
        E'(0,& u_0  ;Q_m)_0 \cap  \left(\bigcap_{i=1}^{n-1} E' (v_{i-1}, u_i; \tCcal_{i-1})_i\right)\cap H'(v_{n-1}, x,y ; \tCcal_{n-1})_n\\
        & \subset \bigcup_{\vec{t}, \vec{s}, \vec{z}} \Biggl(G_0(u_0, s_0, z_1 ;Q_m)_0 \cap \left(\bigcap_{i=1}^{n-1} G(v_{i-1}, t_i, z_i, u_i, s_i, z_{i+1})_i \right) \cap G''_n (v_{n-1}, t_n, z_n, y,x)_n \Biggr).
    \end{split}
\end{equation}
Therefore, we can get an upper bound on $\pi^{(n)}_{m}$ and $\xi^{(n)}_{m}$:
\begin{equation}\label{e:piBK}
    \begin{split}
        \pi^{(n)}_{m}(x,r) \le  & \sum_{\vec{z},\vec{t}, \vec{s},\vec{u},\vec{v}}\left[\prod_{i=0}^{n-1} p_{u_i v_i}\right] \Ppc(G_0(u_0, s_0, z_1;Q_m)) \\
                            & \times \prod_{i=1}^{n-1} \Ppc(G(v_{i-1}, t_i, u_i, s_i, z_{i+1})) \Ppc(G'_n (v_{n-1}, t_n, z_n, x)),
    \end{split}
\end{equation}
where $\vec{u} = (u_0, \dotso, u_{n-1})$ and $\vec{v}=(v_0, \dotso, v_{n-1})$ with all elements are restricted to $\Zd$, and
\begin{equation}\label{e:phiBK}
    \begin{split}
        \phi_{m}^{(n)}(x,y,r) \le  & \sum_{\vec{z},\vec{t}, \vec{s},\vec{u},\vec{v}}\left[\prod_{i=0}^{n-1} p_{u_i v_i}\right] \Ppc(G_0(u_0, s_0, z_1;Q_m)) \\
                            & \times \prod_{i=1}^{n-1} \Ppc(G(v_{i-1}, t_i, u_i, s_i, z_{i+1})) \Ppc(G''_n (v_{n-1}, t_n, z_n, y,x)).
    \end{split}
\end{equation}
The probabilities in \refeq{piBK} and \refeq{phiBK} factorize because $G_0, \dotso, G'_n$ and $G_0, \dotso, G''_n$ are events on different percolation models. The separate probabilities can all be estimated using the BK-inequality. To organize the resulting sum, define
\begin{equation}
    \ttaupc(x) = p_c (D * \taupc)(x)
\end{equation}
and
\begin{eqnarray}
    \label{e:Adiagdef} A (a,b,s,t) &=& \taupc(a-s)\taupc(s-t)\taupc(t-b);\\
    \label{e:Bdiagdef} B_1(s,t,z,l) &=& \ttaupc(l-t)\taupc(z-s);\\
    B^{(0)}_2 (z,l,s,t) &=& \taupc(l-z)\taupc(t-z)\taupc(s-l)\taupc(t-s)\\
    \label{e:B21diagdef} B^{(1)}_2 (z,l,s,t) &=& \sum_{a \in \Zd} \delta_{l,s} \taupc(a-s)\taupc(z-a)\taupc(t-a)\taupc(t-z);\\
    B_2(z,l,s,t) &=& B^{(0)}_2(z,l,s,t) + B^{(1)}_2(z,l,s,t);\\
    \label{e:Cdiagdef} C (a,z,l) &=& A(a,a,z,l) = \taupc(a-z)\taupc(l-a)\taupc(z-l);\\
    D^{(0)}(s,t,z,l,x,y) &=& B_1 (s,t,z,l)A(z,l,x,y);\\
    D^{(1)}(s,t,z,l,x,y) &=& \ttaupc(y-t)\taupc(l-y)\taupc(z-s)C(x,z,l);\\
    \label{e:Ddiagdef} D(s,t,z,l,x,y) &=& D^{(0)}(s,t,z,l,x,y)+ D^{(1)}(s,t,z,l,x,y).
\end{eqnarray}
See Figure \ref{fig:DiagExp} for diagrammatic representations of these functions.

\begin{figure}
  \includegraphics[width=16cm]{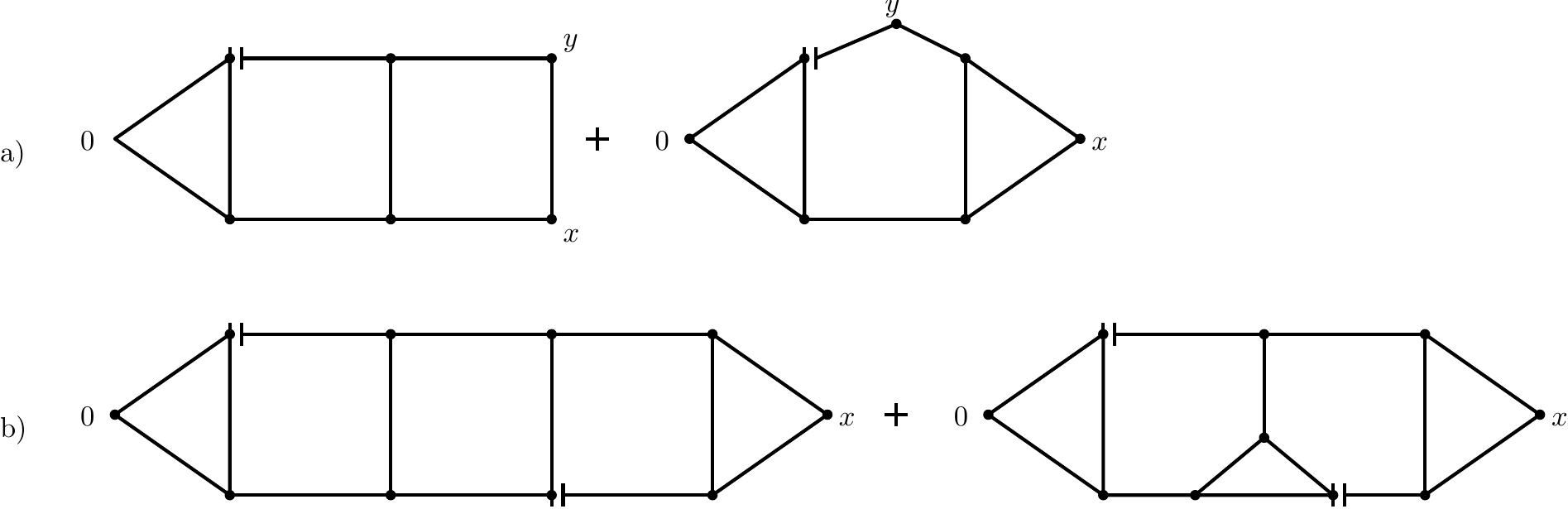}\\
  \caption{Diagrams bounding a) $\bar{\phi}^{(1)}(x,y)$ and b) $\bar{\pi}^{(2)}(x)$ }\label{fig:DiagPhiPi}
\end{figure}

Application of the BK-inequality yields
\begin{eqnarray}
    \label{e:firstdiag} \Ppc(G_0(s_0,t_0,z_1;Q_m)) &\le& \sum_{w \in Q_m} A(0,w,s_0,t_0)\taupc(s_0,z_1);\\
    \label{e:nthdiagpi} \sum_{v_{n-1}\in \Zd} p_{t_{n-1}v_{n-1}} \Ppc(G'_n(v_{n-1},t_n,z_n,x)) &\le& \frac{B_1 (s_{n-1}, t_{n-1}, z_n, l_n)}{\taupc(z_{n}-s_{n-1})} C (x,z_n,l_n);\\
    \label{e:nthdiagxi} \sum_{v_{n-1}\in \Zd} p_{t_{n-1}v_{n-1}} \Ppc(G''_n(v_{n-1},t_n,z_n,y,x)) &\le& \frac{D (s_{n-1}, t_{n-1}, z_n, l_n,x,y)}{\taupc(z_{n}-s_{n-1})}.
\end{eqnarray}
For $G'$ and $G''$ we obtain
\begin{equation} \label{e:ithdiagGp}
        \sum_{v_{i-1} \in \Zd} p_{t_{i-1}v_{i-1}} \Ppc(G'(v_{i-1}, l_i,z_i,s_i,t_i,z_{i+1}))\le \frac{B_1 (s_{i-1}, t_{i-1}, z_i, l_i)}{\taupc(z_{i}-s_{i-1})} B^{(0)}_2 (z_i, l_i,s_i,t_i)\taupc(z_{i+1}-s_i);
\end{equation}
\begin{equation} \label{e:ithdiagGdp}
        \sum_{v_{i-1}, l_i \in \Zd} p_{t_{i-1}v_{i-1}}  \Ppc(G''(v_{i-1},l_i,z_i,s_i,t_i,z_{i+1}))\le \frac{B_1 (s_{i-1}, t_{i-1}, z_i, l_i)}{\taupc(z_{i}-s_{i-1})} B^{(1)}_2 (z_i, l_i,l_i,t_i)\taupc(z_{i+1}-s_i).
\end{equation}
The Kronecker delta in $B^{(1)}_2$ guarantees that it can only be nonzero when its second and third argument are equal, so we can replace the third argument of $B^{(1)}_2$ by $s_i$ and combine \refeq{ithdiagGp} and \refeq{ithdiagGdp} to obtain
\begin{equation}\label{e:ithdiagG}
        \sum_{v_{i-1}, l_i \in \Zd} p_{t_{i-1}v_{i-1}} \Ppc(G(v_{i-1},l_i,z_i,s_i,t_i,z_{i+1}))\le \frac{B_1 (s_{i-1}, t_{i-1}, z_i, l_i)}{\taupc(z_{i}-s_{i-1})} B_2 (z_i, l_i,s_i,t_i)\taupc(z_{i+1}-s_i).
\end{equation}
Substituting \refeq{firstdiag}, \refeq{nthdiagpi} and \refeq{ithdiagG} into \refeq{piBK}, and \refeq{firstdiag}, \refeq{nthdiagxi} and \refeq{ithdiagG} into \refeq{phiBK}, respectively, we obtain, for $n \ge 1$
\begin{equation}\label{e:pidiag}
    \begin{split}
        \pi^{(n)}_{m} (x,r) \le & \sum_{\vec{s},\vec{t}, \vec{z},\vec{l}} \sum_{w \in Q_m} A(0,w,s_0,t_0) \prod_{i=1}^{n-1} \bigl[B_1(s_{i-1}, t_{i-1}, z_{i}, l_i) B_2 (z_i, l_i, s_i, t_i) \bigr]\\
        & \qquad \times B_1(s_{n-1}, t_{n-1}, z_n, l_n) C(x,z_n, l_n)
    \end{split}
\end{equation}
and
\begin{equation}\label{e:phidiag}
    \begin{split}
        \phi^{(n)}_{m} (x,y,r) \le & \sum_{\vec{s},\vec{t}, \vec{z},\vec{l}} \sum_{w \in Q_m} A(0,w,s_0,t_0) \prod_{i=1}^{n-1} \bigl[B_1(s_{i-1}, t_{i-1}, z_{i}, l_i)  B_2 (z_i, l_i, s_i, t_i) \bigr]\\
        & \qquad \times D(s_{n-1}, t_{n-1}, z_n, l_n,x,y).
    \end{split}
\end{equation}
The summation over the vectors $\vec{s}=(s_0, \dotso, s_{n-1})$, $\vec{t}=(t_0, \dotso t_{n-1})$, $\vec{z}=(z_1, \dotso, z_{n})$ and $\vec{l}=(l_1, \dotso, l_n)$ on the right-hand sides of \refeq{pidiag} and \refeq{phidiag} is over all of $\Zd$ for each element, so in both cases the dependence of $r$ has been removed.
Also observe that the sum over $w$ is again restricted to $Q_m$, so that once again we may replace $A(0,w,s_0,t_0)$ by $C(0,s_0, t_0)$ in both instances, to bound, for $n \ge 1$,
\begin{equation}\label{e:pidiag1}
    \begin{split}
        \pi^{(n)}_{m} (x,r) \le C'_m \bar{\pi}^{(n)}(x) := & C'_{m}\sum_{\vec{s},\vec{t}, \vec{z},\vec{l}}  C(0,s_0,t_0) \prod_{i=1}^{n-1} \bigl[B_1(s_{i-1}, t_{i-1}, z_{i}, l_i) B_2 (z_i, l_i, s_i, t_i) \bigr]\\
        & \qquad \times B_1(s_{n-1}, t_{n-1}, z_n, l_n) C(x,z_n, l_n)
    \end{split}
\end{equation}
and
\begin{equation}\label{e:phidiag1}
    \begin{split}
        \phi^{(n)}_{m} (x,y,r) \le C'_m \bar{\phi}^{(n)}(x,y):= & C'_{m}\sum_{\vec{s},\vec{t}, \vec{z},\vec{l}} C(0,s_0,t_0) \prod_{i=1}^{n-1} \bigl[B_1(s_{i-1}, t_{i-1}, z_{i}, l_i)  B_2 (z_i, l_i, s_i, t_i) \bigr]\\
        & \qquad \times D(s_{n-1}, t_{n-1}, z_n, l_n,x,y).
    \end{split}
\end{equation}
The two bounds above are commonly referred to as diagrammatic estimates.
In Figure \ref{fig:DiagPhiPi} we show two examples of diagrams.

Finally, for ease of notation in the coming sections, we define
\begin{equation}\label{e:PiPhidef}
    \bar{\Pi} (x) = \sum_{n=0}^{\infty} \bar{\pi}^{(n)}(x) \qquad \text{ and }\qquad \bar{\Phi} (x,y) = \sum_{n=0}^{\infty} \bar{\phi}^{(n)}(x,y).
\end{equation}

\section{Finite moments of $\bar{\Pi}(x)$: proof of Proposition \ref{prop:TermBounds}(ii)}\label{sec:proofprop252}

In this section we prove Proposition \ref{prop:TermBounds}(ii), which states that the $(\twa+\delta)${\rm 'th} moment of $\lvert\Pi(x,r;F) \rvert$ is finite for some $\delta >0$. We do this by showing
\begin{equation}
    \sum_{x\in\Zd} \lvert x\rvert^{\twa+\delta} \bar{\Pi} (x) \le K'
\end{equation}
for $\bar{\Pi}(x)$ as defined in \refeq{PiPhidef}.
In showing this, we bound certain quantities that are similar to quantities bounded by Chen and Sakai \cite{CheSak09}, and the proof of this bound is based in part on their proofs.

We assume $p=p_c$ throughout and suppress all subscripts $p_c$. We also omit the area of integration $(-\pi,\pi]^d$ below the integral signs, whenever it occurs.
\\
\proof[Proof of Proposition \ref{prop:TermBounds}(ii)]
The proof is split up into three sections. In the first section we describe a way of distributing the weight $\lvert x \rvert^{\twa+\delta}$ over the path elements of the diagram. The second section deals with taking the Fourier transform of lace expansion diagrams. In the third section we bound the elements of these Fourier space diagrams.

\subsection{Distributing the weight}
For $\alpha >0$ and $d > 3 \twa$ we choose $\delta$ such that
\begin{equation}
    \delta \in (0, \twa \wedge (d - 3 \twa)\wedge 1).
\end{equation}
By the definition of $\bar{\Pi}(x)$,
\begin{equation}\label{e:Pitopi}
    \sum_{x \in \Zd} |x|^{\twa + \delta} \bar{\Pi} (x) = \sum_{x \in \Zd} \sum_{n=0}^{\infty} |x|^{\twa + \delta} \bar{\pi}^{(n)} (x).
\end{equation}

For $x \in \Zd$ we write $x = (x_1, x_2, \dotso, x_d)$. Because the functions $\bar{\pi}^{(n)}(x)$ are invariant under the symmetries of $\Zd$, we can bound \refeq{Pitopi} as follows:
\begin{equation}\label{e:pionedim}
    \sum_{x \in \Zd} \sum_{n=0}^{\infty} |x|^{\twa + \delta} \bar{\pi}^{(n)} (x) \le d^{(\twa + \delta)/2 +1} \sum_{x \in \Zd} \sum_{n=0}^{\infty} |x_1|^{\twa + \delta}\bar{\pi}^{(n)} (x).
\end{equation}
To make the sum tractable, we `distribute the weight $\lvert x_1 \rvert^{\twa + \delta}$ along the top and bottom paths of the diagram': For $t>0$ and $\zeta \in (0,2)$, let
\begin{equation}
    K'_\zeta \equiv \int\limits_{0}^{\infty} \frac{1- \cos(x)}{x^{1+\zeta}} \dnd x \in (0, \infty).
\end{equation}
This gives the identity
\begin{equation}\label{e:tzIdentity}
    t^{\zeta} = \frac{1}{K'_\zeta} \int\limits_{0}^{\infty} \frac{1-\cos(s t)}{s^{1+\zeta}} \dnd s.
\end{equation}

For $u,v \in (0,\infty)$, define the $d$-dimensional vectors $\vec{u} = (u,0,\dotso,0)$ and $\vec{v} = (v,0,\dotso,0)$.
Let $\delta_1$ and $\delta_2$ be constants, such that
\begin{equation}\label{e:deltadef}
    \delta_1 \in \left(\delta, \twa \wedge \left(\frac{(1+\delta) d}{3+\delta} - \twa\right)\right), \quad \text{and} \quad \delta_2 = \twa + \delta - \delta_1,
\end{equation}
so that $\delta_1 + \delta_2 = \twa+\delta$. Applying \refeq{tzIdentity} twice to \refeq{pionedim} with $\zeta= \delta_1, \delta_2$, we obtain
\begin{equation}\label{e:picos}
    \sum_{x \in \Zd}  \sum_{n=0}^{\infty} |x|^{\twa + \delta} \bar{\pi}^{(n)} (x) \le C \int\limits_{0}^{\infty} \frac{\dnd u}{u^{1+\delta_1}} \int\limits_{0}^{\infty} \frac{\dnd v}{v^{1+\delta_2}} \sum_{x \in \Zd} \sum_{n=0}^{\infty} [1-\cos(\vec{u}\cdot x)] [1-\cos(\vec{v}\cdot x)] \bar{\pi}^{(n)}(x).
\end{equation}
\begin{rk} The exponent $\delta_1$ can be viewed as being only slightly larger than $\delta$, making $\delta_2$ only slightly smaller than $\twa$. Were we to consider the case where $\delta_1 = 0$, this would reduce the problem to that of Proposition 4.1 in \cite[(4.32)]{BorChaHofSlaSpe05b}.
\end{rk}

The double integral can be split into four parts: $I_1 + I_2 + I_3 + I_4$, where
\begin{equation}
    I_1 = O(1) \int\limits_{0}^{1} \frac{\dnd v}{v^{1+\delta_1}} \int\limits_{0}^{1} \frac{\dnd u}{u^{1+\delta_2}}\sum_{x \in \Zd} \sum_{n=0}^{\infty} [1-\cos(\vec{u}\cdot x)] [1-\cos(\vec{v}\cdot x)]\bar{\pi}^{(n)}(x)
\end{equation}
and $I_2, I_3$ and $I_4$ are similarly defined but with different areas of integration $A_i$, $i=2,3,4$, where
\begin{equation}
    A_2 = [0,1] \times (1, \infty], \quad A_3 = (1,\infty] \times [0,1], \quad \text{and} \quad A_4 = (1, \infty] \times (1, \infty].
\end{equation}

It remains to show that $I_1, \dotso, I_4$ are finite.
To prove that this is so, we need an upper bound on
\begin{equation}\label{e:cosxysum}
    \sum_{x \in \Zd} \sum_{n=0}^{\infty} [1 - \cos(\vec{u}\cdot x)][1 -\cos(\vec{v}\cdot x)] \bar{\pi}^{(n)}(x).
\end{equation}
Let
\begin{equation}
        \gamma_1 = \left(\frac{(1+\delta) d}{3+\delta}-  \twa\right)\wedge \twa, \quad \text{ and } \quad \gamma_2 = \left(\frac{2 d}{3+\delta} - \twa\right)\wedge \twa,
\end{equation}
then $\gamma_1 > \delta_1$ and $\gamma_2 > \delta_2$. Proposition \ref{prop:TermBounds}(ii) follows once we show
\begin{equation}\label{e:Wuvbound}
        \sum_{x \in \Zd} \sum_{n=0}^{\infty} [1 - \cos(\vec{u}\cdot x)][1 -\cos(\vec{v}\cdot x)] \bar{\pi}^{(n)}(x) = O\left((u\wedge 1)^{\gamma_1} (v\wedge 1)^{\gamma_2}\right).
\end{equation}

The bounds are easy for $u$ or  $v$ in $(1,\infty]$. In particular, $I_4 < \infty$ follows from the fact that $\sum_{x \in \Zd} \sum_{n=0}^{\infty}  \bar{\pi}^{(n)}(x) \le C < \infty$ and $1-\cos(t) \le 2$.

The remainder of this section is devoted to proving \refeq{Wuvbound} when both $u,v \in [0,1]$, that is, the bound needed for the finiteness of $I_1$. The bounds on $I_2$ and $I_3$ can be obtained in a similar, but much easier, way.

\begin{figure}
  \includegraphics[width=490pt]{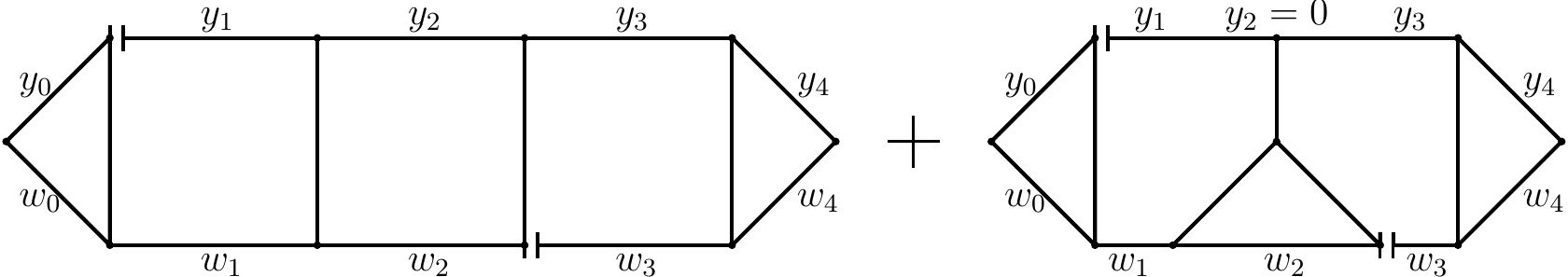}\\
  \caption{The path-elements of $\bar{\pi}^{(2)}(x)$ labeled according to the proposed scheme.}\label{pic:defyiwi}
\end{figure}

We start by only considering $n\ge1$. The case $n=0$ is much simpler, and we will comment on the right bound for $n=0$ when it is appropriate (around equation \refeq{cospi0}). Using \refeq{pidiag} we can rewrite the right-hand side of \refeq{picos} (with the term for $n=0$ omitted) as
\begin{equation}\label{e:pionedimdiag}
    \begin{split}
        C & \int\limits_{0}^{\infty} \frac{\dnd u}{u^{1+\delta_1}} \int\limits_{0}^{\infty} \frac{\dnd v}{v^{1+\delta_2}} \sum_{x \in \Zd}\sum_{n=1}^{\infty} [1-\cos(\vec{u} \cdot x)] [1-\cos(\vec{v}\cdot x)] \sum_{\vec{s}, \vec{t}, \vec{z}, \vec{l}} C (0, s_0, t_0) \\
        & \times \prod_{m=1}^{n-1} \left[ B_1 (s_{m-1}, t_{m-1}, z_m, l_m)B_2 (z_m, l_m, s_m, t_m) \right] B_1 (s_{n-1}, t_{n-1}, z_n, l_n ) C (x , z_n, l_n).
    \end{split}
\end{equation}
Define for $i = 0,1, \dotso, n$:
\begin{eqnarray}
    y_{2i} = \left\{  \begin{array}{ll}   t_0 & \text{if }i=0; \\
                                            t_i - z_i & \text{if $i$ is odd;}\\
                                            s_i - l_i & \text{if $i$ is even;}
                        \end{array}\right.
    \qquad
        \begin{array}{l}    y_{2i+1} = \left\{    \begin{array}{ll}   l_i - t_{i-1} & \text{if $i<n$ is odd;}\\
                                                                    z_i - s_{i-1} & \text{if $i<n$ is even;}
                                                \end{array}\right. \\
                            y_{2n} = \left\{  \begin{array}{ll}   x - z_n & \text{if $n$ is odd;}\\
                                                                    x - l_n & \text{if $n$ is even;}
                                                \end{array}\right.
        \end{array} \\
    w_{2i} = \left\{  \begin{array}{ll}   s_0 & \text{if }i=0; \\
                                            s_i - l_i & \text{if $i$ is odd;}\\
                                            t_i - z_i & \text{if $i$ is even;}
                        \end{array}\right.
    \qquad
        \begin{array}{l}    w_{2i+1} = \left\{    \begin{array}{ll}   t_i - z_{i-1} & \text{if $i<n$ is odd;}\\
                                                                    s_i - l_{i-1} & \text{if $i<n$ is even;}
                                                \end{array}\right. \\
                            w_{2n} = \left\{  \begin{array}{ll}   x - l_n & \text{if $n$ is odd;}\\
                                                                    x - z_n & \text{if $n$ is even.}
                                                \end{array}\right.
        \end{array}
\end{eqnarray}
The $y$'s and $w$'s can be viewed as the path elements along the top and bottom of the diagram $\bar{\pi}^{(n)}$, respectively. An example is given in Figure \ref{pic:defyiwi}.

The result is that we obtain two telescoping sums:
\begin{equation}
    \sum_{i=0}^{2n} y_i = \sum_{i=0}^{2n} w_i = x.
\end{equation}

By \cite[(4.51)]{BorChaHofSlaSpe05b}, for $a = \sum_{j=1}^J a_j$,
\begin{equation}\label{e:costbd}
    1- \cos a \le (2 J + 1) \sum_{j=1}^{J}[ 1 - \cos a_j].
\end{equation}

Applying this with $a_i = \vec{u} \cdot y_i$ and $\vec{v} \cdot w_i$ gives that \refeq{pionedimdiag} is bounded from above by
\begin{equation}\label{e:Rnij}
    \begin{split}
        C \int\limits_{0}^{\infty} \frac{\dnd u}{u^{1+\delta_1}}& \int\limits_{0}^{\infty} \frac{\dnd v}{v^{1+\delta_2}} \sum_{x \in \Zd}\sum_{n=1}^{\infty} (4n+3)^2  \sum_{i,j = 0}^{2n} \sum_{\vec{s}, \vec{t}, \vec{z}, \vec{l}} [1-\cos(\vec{u}\cdot y_i)] [1-\cos(\vec{v} \cdot w_j)]\\
        &\times C (0, s_0, t_0) \prod_{m=1}^{n-1} \left[ B_1 (s_{m-1}, t_{m-1}, z_m, l_m)B_2 (z_m, l_m, s_m, t_m) \right] B_1 (s_{n-1}, t_{n-1}, z_n, l_n ) C (x , z_n, l_n) \\
        &\qquad \qquad \equiv C  \int\limits_{0}^{\infty} \frac{\dnd u}{u^{1+\delta_1}} \int\limits_{0}^{\infty} \frac{\dnd v}{v^{1+\delta_2}} \sum_{n=1}^{\infty} (4n+3)^2 \sum_{i,j = 0}^{2n} \Anij.
    \end{split}
\end{equation}

Each of the $\Anij$ is  the sum of $2^{n-1}$ terms: one for each sequence of $B_{2}^{(0)}$ and $B_{2}^{(1)}$ diagrams possible. The possible sequences of $B_{2}^{(0)}$ and $B_{2}^{(1)}$ diagrams from left to right (say), corresponds one-to-one to the binary expansion of an integer between $0$ and $2^{n-1} -1$, so we can write
\begin{equation}\label{e:Anmij}
    \Anij \equiv \sum_{m=0}^{2^{n-1}-1} \Anmij
\end{equation}
where each of the $\Anmij$ corresponds to exactly one realization of a diagram. Furthermore, the diagrams are products of functions of two variables, the (possibly weighted) connectivity functions. Hence, we can associate a graph to each of the $\Anmij$ in such a way that the edges of the graph correspond to the two-variable functions of $\Anmij$ and the vertices correspond to the variables in $\Zd$ that are being summed over.
The graph structure implies certain properties of the Fourier transform of the diagrams that are useful in obtaining upper bounds.

We use these properties to bound the diagrams in Fourier space. Our strategy is the following: The first step is to use graph properties to write $\Anmij$ as the integral over a function of $2n+1$ Fourier variables, rather than the $6n+2$ variables that would be obtained from taking the Fourier transform for each of the $6n+2$ connectivity functions that are contained in $\Anmij$ separately. Then, using a duality argument on the graph structure, we determine the order in which to integrate over these $2n+1$ variables (similar approaches exist for bounding Feynman diagrams in the quantum field theory literature, cf.\ \cite{FelKnoSalTru99}, \cite{FelMagRivSen85}). We show that there exists an order such that we can integrate over the product of at most three functions of the same variable. Roughly speaking, this corresponds to integration over the triangle diagram in Fourier space, which we assumed to be bounded by a small constant in the statement of Proposition \ref{prop:TermBounds}. This way we are able to show that all the $\Anmij$ have an upper bound of the order of $\beta^{n-3} u^{\gamma_1} v^{\gamma_2}$, and this suffices to show \refeq{Wuvbound} and hence Proposition \ref{prop:TermBounds}(ii) holds.

\subsection{Fourier space diagrams}
We start by carrying out the above program with some general considerations. Let $\Mcal$ be an inner product space. Let $\Vcal$ be a set of vertices with $\lvert \Vcal \rvert = V$  and let $\Ecal \subseteq \Vcal \times \Vcal$ with $\lvert \Ecal \rvert = E$ be a set of unoriented edges (i.e., $\{i,j\} = \{j,i\}$). The graph $\Gcal=(\Vcal, \Ecal)$ plays the role of an index set for a diagram. We call a function $F: \Mcal^V \mapsto \R$ an \emph{edge diagram} if it can be written as a product of functions on `edges' as indexed by $\Ecal$, i.e.,
\begin{equation}\label{e:edgediag}
    F(x_1, \dotso, x_{V}) =  \prod_{\{i,j\} \in \Ecal} f_{i,j} (x_i, x_j),
\end{equation}
where $f_{i,j}: \Mcal^2 \mapsto \R^{+}$.
We call the edge diagram simple and connected, respectively, if the associated graph $\Gcal$ is simple and connected.

We say that the functions $f_{i,j}$ are translation invariant if, for every $a \in \Mcal$,
\begin{equation}
    f_{i,j}(x,y) = f_{i,j} (x+a, y+a).
\end{equation}

The upcoming lemma and its proof make use of certain elementary graph theoretic notions.
It is a basic fact from graph theory that associated to every graph $\Gcal$ there is a vector space $\Ccal(\Gcal)$ whose elements represent formal combinations of cycles in $\Gcal$. We call this vector space the \emph{cycle space of }$\Gcal$. Given a spanning tree $\Tcal = (\Vcal, \Ecal')$ of $\Gcal$, we can define a \emph{fundamental cycle of }$\Tcal$ as the single cycle in the graph $\Scal=(\Vcal, \Ecal' \cup e)$ for an edge $e \in \Ecal \setminus \Ecal'$. The set of all fundamental cycles of $\Tcal$ is a basis for $\Ccal(\Gcal)$.
For further definitions and a proof of the above statement we refer the reader to the literature of this field (e.g. \cite{Dies06}).

\bl[Characterizing the independent Fourier variables]\label{lem:ConMom} Let $F(x_1, \dotso, x_V)$ be a translation invariant simple and connected edge diagram indexed by a graph $\Gcal = (\Vcal, \Ecal)$ with $V$ vertices and $E$ edges. Then the Fourier transform of $F$ can be expressed in terms of $E-V+1$ linearly independent Fourier variables. These variables can be chosen to correspond to a basis of the cycle space $\Ccal(\Gcal)$ of $\Gcal$.
\el

\proof We start by expressing $f_{i,j}(x_i, x_j)$ in terms of its Fourier transform:
\begin{equation}\label{e:fij}
    f_{i,j}(x_i, x_j) = \int\limits_B \frac{\dd p_i}{\lvert B \rvert} \int\limits_B \frac{\dd p_j}{\lvert B \rvert} e^{i p_i \cdot x_i} e^{i p_j \cdot x_j} \hat{f}_{i,j} (p_i, p_j),
\end{equation}
where $B$ is the fundamental domain of the reciprocal space of $\Mcal$ (e.g. when $\Mcal = \Zd$, then $B=(-\pi, \pi]^d$).

Shifting $x_i$ and $x_j$ by a vector $a$ we obtain
\begin{equation}\label{e:fijtrans}
    \begin{split}
    f_{i,j}(x_i +a, x_j +a) &= \int\limits_B \frac{\dd p_i}{\lvert B \rvert} \int\limits_B \frac{\dd p_j}{\lvert B \rvert} e^{i p_i \cdot (x_i +a)} e^{i p_j \cdot (x_j + a)} \hat{f}_{i,j} (p_i, p_j)\\
                            &=\int\limits_B \frac{\dd p_i}{\lvert B \rvert} \int\limits_B \frac{\dd p_j}{\lvert B \rvert} e^{i p_i \cdot x_i} e^{i p_j \cdot x_j} e^{i(p_i + p_j) \cdot a} \hat{f}_{i,j} (p_i, p_j).
    \end{split}
\end{equation}
However, by translation invariance, the left-hand sides of \refeq{fij} and \refeq{fijtrans} are equal, and so the right-hand sides must also be equal. This is only the case for every $a \in \Mcal$ when
\begin{equation}
    e^{i(p_i + p_j) \cdot a}=1 \quad \text{ or, equivalently, }\quad p_i + p_j = 0 \mod B.
\end{equation}
We have such a constraint for every pair $p_i$, $p_j$ for which $\{i,j\} \in \Ecal$. Therefore, we can write these constraints as a system of linear equations, so that in matrix notation we have
\begin{equation}
    A \cdot \vec{p} = \vec{0}.
\end{equation}
where $\vec{p}$ is a vector of length $V$ with entries $p_i \in B$, $i=1,\dotso,V$, and $A$ is a $V \times E$ matrix. In fact, the transpose of $A$ is the incidence matrix of $\Gcal$. As such, it is an elemental result from graph theory (\cite[Proposition 1.9.7]{Dies06}), that the rank of $A$ is equal to the dimension of $\Ccal (\Gcal)$, the cycle space of $\Gcal$. Another elemental result is that the dimension of $\Ccal(\Gcal)$ is $E-V+1$ (\cite[Theorem 1.9.6]{Dies06}), so the rank of $A$ is $E-V+1$, and hence there are $E-V+1$ linearly independent Fourier variables associated to the Fourier transform of $F$.

Furthermore, since the kernel of $A^T$ is $\Ccal(\Gcal)$ (\cite[Proposition 1.9.7]{Dies06}), we can express these linearly independent Fourier variables in terms of a basis of $\Ccal(\Gcal)$, such as the fundamental cycles of a spanning tree $\Tcal$ of $\Gcal$. \qed

\begin{figure}
  \includegraphics[width=300pt]{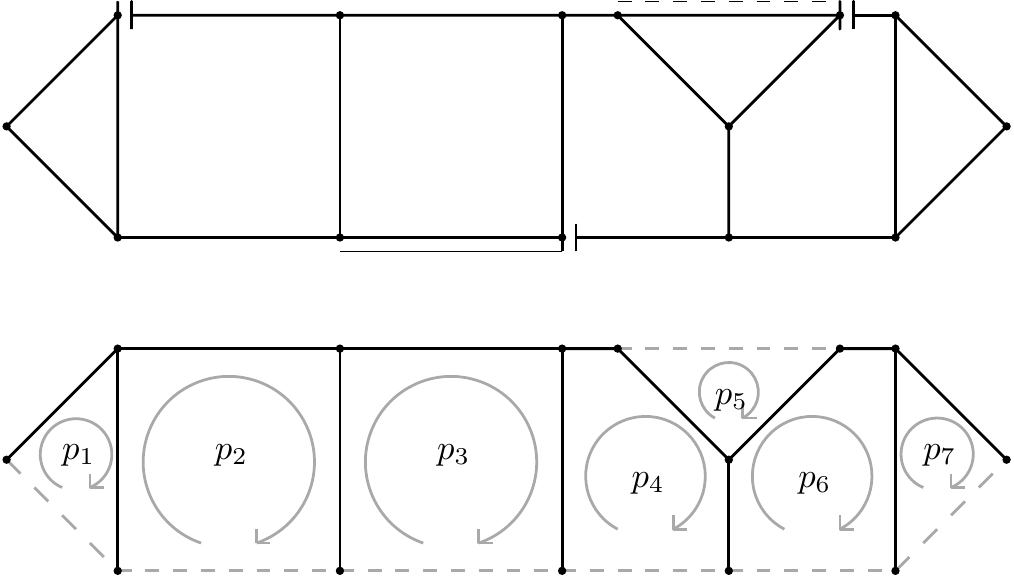}
  \includegraphics[width=170pt]{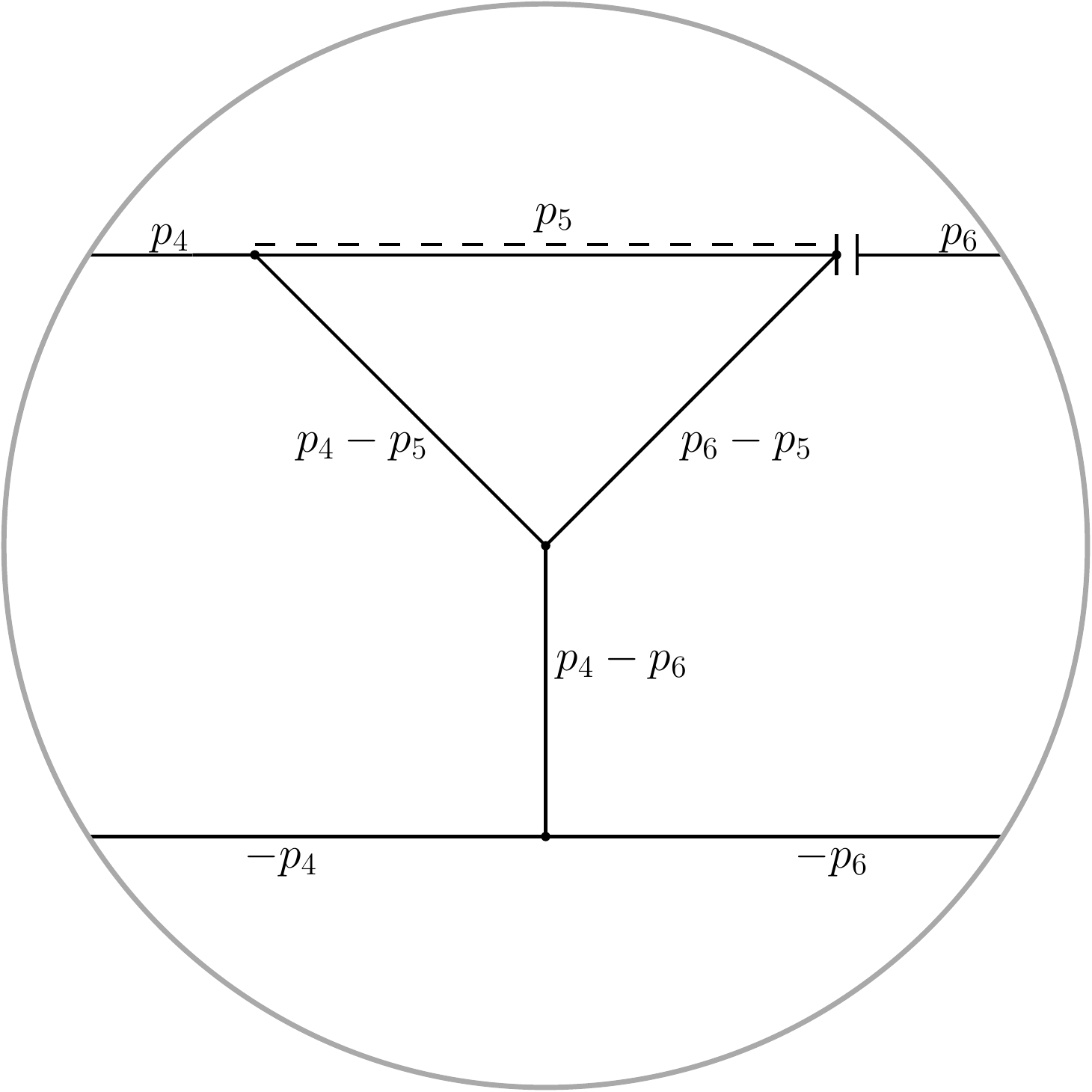}\\
  \caption{On the left the diagram $\Rcal^{\sss (2,1)}_{\sss(5,3)}$ and its spanning tree, with associated loop momenta. The solid and the dashed line in the upper diagram represent the weighted paths. On the right a portion of the diagram with loop momenta associated to the lines.}\label{pic:diagtree}
\end{figure}

Recall definition \refeq{Anmij}. Every $\Anmij$ is an edge diagram with all its variables summed over, so we can write
\begin{equation}
    \Anmij = \sum_x \sum_{\vec{z},\vec{t},\vec{s},\vec{l}} \Fiox
\end{equation}
where $\iota$ is a shorthand for the quartet of indices $n,m,i,j$ and the dependence of $\Fio$ on $\vec{u}$ and $\vec{v}$ is implicit.
Let $\Gio$ be the graph associated to $\Anmij$. From the construction of the diagrams it follows that all $\Gio$ are planar graphs. Furthermore, $\Gio$ has $6n+2$ edges and $4n+2$ vertices, so by Lemma \ref{lem:ConMom}, the Fourier transform of $\Fio$ has $2n+1$ independent variables. By choosing the right spanning tree of $\Gio$, these variables can be associated to loops along the $2n+1$ internal faces of the graph $\Gio$ (cf. Figure \ref{pic:diagtree}). There are $2n+2$ faces (also counting the external face) and $2n+1$ linearly independent variables, so Fourier variable associated to the external face can be set to zero. Furthermore, we are free to choose the direction of the variables. We always take the variables to run clockwise along a face. In the physics literature, such variables are commonly known as \emph{loop momenta} and hence we use the same term.

A property of planar graphs is that each edge lies between exactly two faces (where the area on the `outside' of the graph is also considered a face). Furthermore, it is a well-known fact from graph theory that each planar graph $\Gcal$ has a unique (up to isomorphism) dual (multi-)graph $\Gcal^{\star}$ such that each vertex of $\Gcal^{\star}$ can be associated to a face of $\Gcal$, and each edge of $\Gcal$ is crossed by exactly one dual edge of $\Gcal^{\star}$ and vice versa. It follows that the degree of vertices in $\Gcal^{\star}$ corresponds to the number of sides of the associated face in $\Gcal$, and therefore, it corresponds to the number of separate occurrences of the associated loop momentum in the Fourier transform of the edge diagram that $\Gcal$ indexes.
\begin{figure}
  \includegraphics[width=235pt]{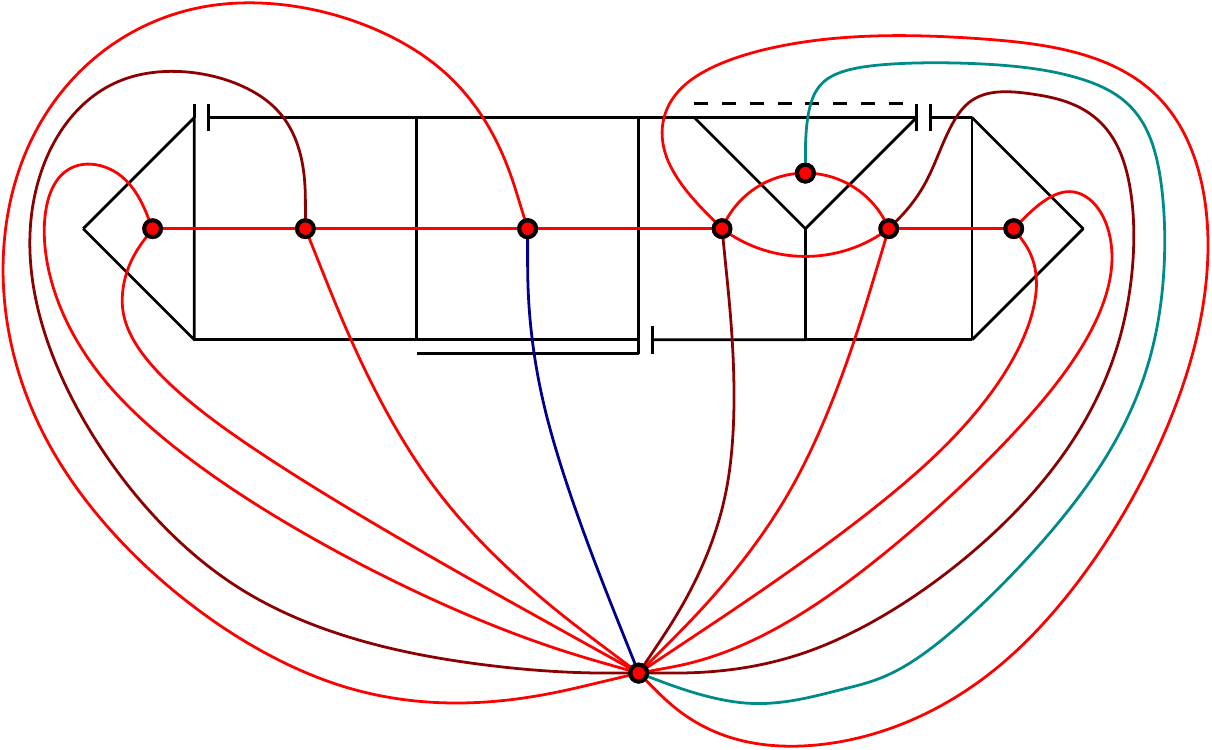}
  \includegraphics[width=235pt]{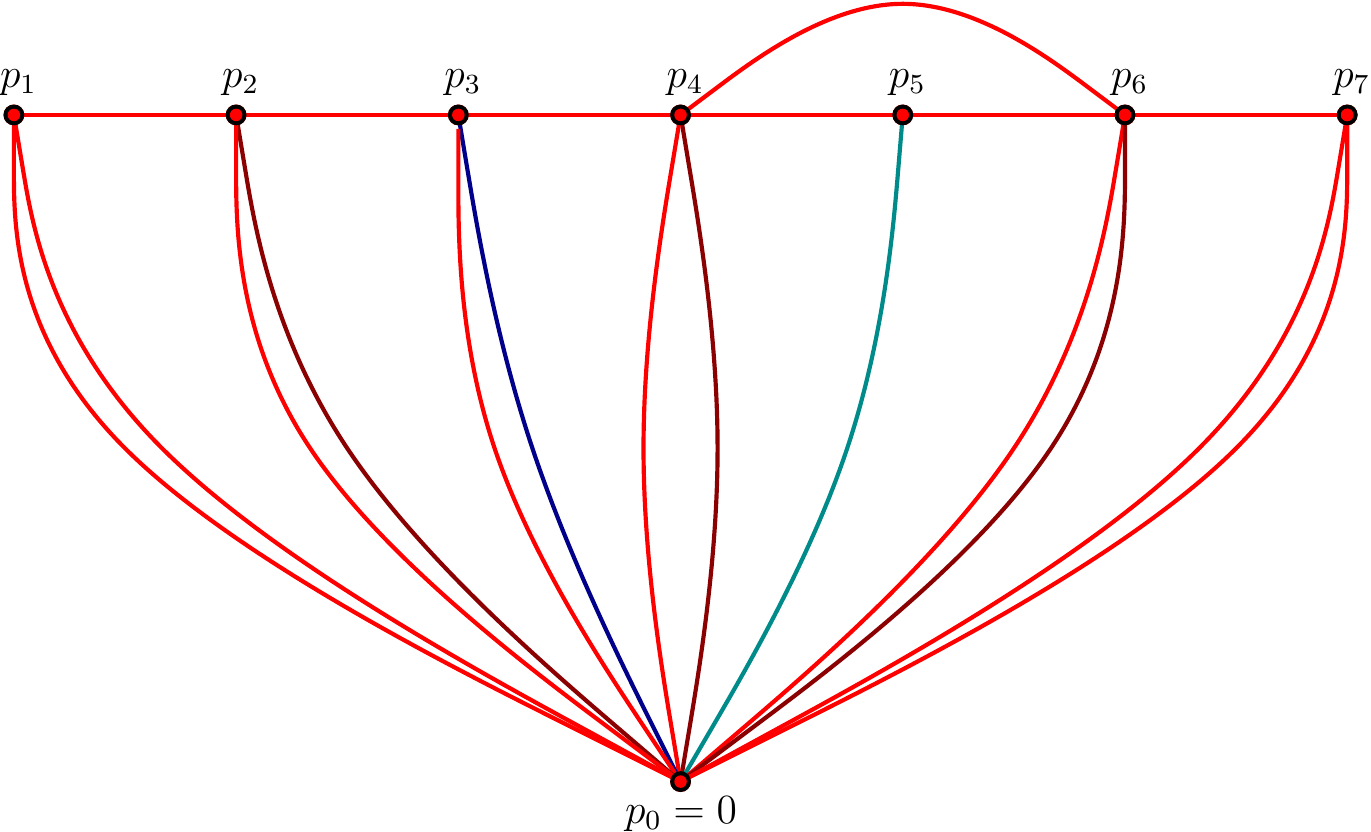}\\
  \caption{On the left the diagram $\Rcal^{\sss (2,1)}_{\sss (5,3)}$ and its dual. On the right an isomorphism of the dual.}\label{pic:diagdual}
\end{figure}

Hence, the dual graph $(\Gio)^{\star}$ indexes $\hFio$, the Fourier transform of $\Fio$. The dual graph $(\Gio)^{\star}$ again has a very simple structure that allows us to write $\hFio$ as the product of $2n+2$ simple elements. In Figure \ref{pic:diagdual} we show an example of a diagram and its dual diagram.

The construction that follows does not work for one particular subset of $\Anmij$, namely those where a weight is associated to a path element that is forced to be zero by the Kronecker delta in the definition of $B^{(1)}_2$, \refeq{B21diagdef}. These weights are an artifact of our notation: they are trivially zero. From here on, we assume that the weights lie on path elements that have a non-zero displacement.

Define
\begin{eqnarray}
    \label{e:Bcaldef} \Bcal(p_a, p_b) &=& \htau(p_a) \htau(p_a - p_b) \htau(p_b);\\
    \widetilde{\Bcal}(p_a, p_b) &=& \hat{D}(p_a) \htau(p_a) \htau(p_a-p_b) \htau(p_b);\\
    \Ccal(p_a, p_b, p_c) &=& \htau(p_a-p_b) \htau(p_b) \htau(p_b-p_c).
   \end{eqnarray}
Also define the functions
\begin{eqnarray}
    \label{e:taul} \tau_q (x) &=& [1- \cos(q \cdot x)]\tau (x);\\
    \label{e:ttaul} \widetilde{\tau}_q (x) &=& [1-\cos(q \cdot x)](D * \tau) (x),
\end{eqnarray}
\begin{equation}
    \bDc_q (p_i) = \left\{ \begin{array}{ll}   \httau_q (p_i)/\httau(p_i) & \text{ if }i/2 \text{ is an odd integer;}\\
                                            \htau_q (p_i)/\htau(p_i) & \text{ otherwise}
                        \end{array} \right.
\end{equation}
and
\begin{equation}
    \uDc_q (p_i) = \left\{ \begin{array}{ll}   \httau_q (p_i)/\httau(p_i) & \text{ if }i/2 \text{ is an even integer;}\\
                                            \htau_q (p_i)/\htau(p_i) & \text{ otherwise.}
                        \end{array} \right.
\end{equation}

Write $m$ as a binary expansion, i.e., $m = m_{n-1} \dotsm m_2 m_1$. Taking the Fourier transform of $\Anmij$, using the definitions \refeq{Bdiagdef} -- \refeq{Cdiagdef}, and rewriting the Fourier variables in terms of the loop momenta as described above, we obtain
\begin{eqnarray}
       \nonumber \Anmij      &=&\int \frac{\dd p_{1}}{(2 \pi)^d} \dotsm \int \frac{\dd p_{2n+1}}{(2 \pi)^d} \htau(p_1) \Bcal(p_1, p_2)\\
       \nonumber             & & \times \left[\prod_{\ell=1}^{n-1} \delta_{0, m_{\ell}} \widetilde{\Bcal}(p_{2 \ell}, p_{2 \ell+1}) \Bcal(p_{2 \ell+1}, p_{2 \ell + 2}) + \delta_{1,m_\ell}\widetilde{\Bcal}(p_{2 \ell}, p_{2 \ell + 2}) \Ccal(p_{2 \ell}, p_{2 \ell+1}, p_{2 \ell+2})\right]\\
       \label{e:AnmijF}      & &\times \widetilde{\Bcal}(p_{2n}, p_{2n+1}) \htau(p_{2n+1}) \bDc_{\vec{u}} (p_i) \uDc_{\vec{v}} (p_j)
\end{eqnarray}
where $\delta_{0,m_\ell}$ and $\delta_{1, m_\ell}$ are Kronecker deltas.

\subsection{A recursive scheme for bounding $\Anmij$}
The simple structure that $\Anmij$ has in Fourier space allows us to recursively integrate over all the variables in such a way that all integrals converge. This is not necessarily obvious if we perform the integrals in some arbitrary order. There may be as many as six functions of the same loop momentum, while we know that the integrals converge when there are at most three connectivity functions present. Furthermore, we need to take special care of the integral over the variables $p_i$ and $p_j$, as the weights makes the integral more divergent (though the weight on $p_i$ more so than the one on $p_j$).

One of the main tools we need for bounding $\Anmij$ is the following iterative version of H\"olders inequality:
\bl[An application of H\"older's inequality]\label{lem:Holder}
    For any $n\ge 2$, let $\alpha_1, \dotso, \alpha_n \in \mathbb{R}^{+}$. Let $S_n = \sum_{i=1}^{n} \alpha_i$. Let $f_1, \dotso, f_n$ be $L^{S_n}$-integrable functions. Then
    \begin{equation}
        \int \prod_{i=1}^n f_{i}(x)^{\alpha_i} \dnd x \le \prod_{i=1}^n \left( \int f_i (x)^{S_n} \dnd x \right)^{\alpha_i / S_n}.
    \end{equation}
\el
\proof The proof is by induction over $n$. The case $n=2$ follows directly from H\"older's inequality with conjugates $S_2/\alpha_1$ and $S_2/\alpha_2$. The inductive step is performed by applying H\"older's inequality with conjugates $S_n/\alpha_n$ and $S_n/S_{n-1}$ to establish that the hypothesis holds for $n$ if it holds for $n-1$.
\qed
\medskip
Note that for any function $f\colon\Zd \mapsto \R$, its Fourier transform $\hat{f}(k)$ will be periodic with period $2 \pi$ in all dimensions, and therefore, we have for any vector $\vec{q}$ and any $s \in \R$,
\begin{equation}\label{e:intshift}
    \int \limits_{(-\pi, \pi]^d} \dd k \hat{f}(k + \vec{q})^s = \int \limits_{(-\pi, \pi]^d + \vec{q}} \dd k' \hat{f}(k')^s = \int \limits_{(-\pi, \pi]^d} \dd k' \hat{f}(k')^s.
\end{equation}

One of the bounds that the recursion is based on is
\begin{equation}\label{e:bdAB}
    \begin{split}
        \int \frac{\dd p_a}{(2 \pi)^d}  \htau(p_a) \Bcal(p_a, p_b) &=  \htau(p_b) \int \frac{\dd p_a}{(2 \pi)^d} \htau(p_a)^2 \htau (p_a - p_b) \\
        &\le  \htau(p_b) \left(\int \frac{\dd p_a}{(2 \pi)^d} \htau(p_a)^3\right)^{2/3}\left(\int \frac{\dd p_a}{(2 \pi)^d} \htau (p_a - p_b)^3 \right)^{1/3} \\
        &= \htau(p_b)\int \frac{\dd p_a}{(2 \pi)^d} \htau(p_a)^3 \le \htau(p_b)\tbar.
    \end{split}
\end{equation}
where $\tbar$ is given in \refeq{triangle}. The first inequality follows from Lemma \ref{lem:Holder}, the second equality follows from \refeq{intshift}, and the second inequality is a consequence of the triangle condition.
In a similar vein, but with a slightly longer calculation, it can be shown that
\begin{equation}\label{e:bdAtB}
        \int \frac{\dd p_a}{(2 \pi)^d} \htau(p_a) \widetilde{\Bcal}(p_a, p_b) \le T\, \htau(p_b)
\end{equation}
where $T$ is given in \refeq{triangleT}.
Furthermore, it also follows from Lemma \ref{lem:Holder} and \refeq{intshift} that
\begin{equation}\label{e:bdC}
    \int \frac{\dd p_b}{(2 \pi)^d} \Ccal(p_a, p_b, p_c) \le \tbar.
\end{equation}

From the bounds \refeq{bdAB}, \refeq{bdAtB} and \refeq{bdC} it is easy to see that we can perform the integrals over the Fourier variables that are not associated with a term $\bDc_{\vec{u}}$ or $\uDc_{\vec{v}}$ in \refeq{AnmijF} in such a way that we can bound every integral by either a factor $T$ or a factor $\tbar$.

Associate the terms $\bDc_{\vec{u}}(p_i)$ and $\uDc_{\vec{v}} (p_j)$ with the first term $\htau, \Bcal, \widetilde{\Bcal}$ or $\Ccal$ of the same variable, as seen when viewed from left to right in the Fourier diagram's construction in \refeq{AnmijF}.

Assume for the moment that $i \neq j$. Sequentially integrate over all other Fourier variables from the left, the right, and any $\Ccal$ that may be in between using the bounds \refeq{bdAB}, \refeq{bdAtB} and \refeq{bdC}.  Once all these variables are integrated over, the resulting expression either contains an integral of the form
\begin{eqnarray}
    \Xcal(\vec{v}) &\equiv& \int \frac{\dd p_j}{(2 \pi)^d} \htau(p_j) \uDc_{\vec{v}} (p_j) \Bcal(p_j, p_{a});\\
    \underline{\Xcal}(\vec{v}) &\equiv& \int \frac{\dd p_j}{(2 \pi)^d} \httau(p_j) \uDc_{\vec{v}} (p_j) \Bcal(p_j, p_{a})
\end{eqnarray}
(where the value of the second index depends on the structure of the diagram) or an integral of the form
\begin{equation}
    \Xcal'(\vec{v}) \equiv  \int \frac{\dd p_i}{(2 \pi)^d}  \Ccal (p_{i-1}, p_i, p_{i+1})\uDc_{\vec{v}} (p_i).
\end{equation}
It can be shown that
\begin{equation}\label{e:bdXs}
    \Xcal(\vec{v}) = O(v^{\gamma_2}) \htau(p_a), \qquad \underline{\Xcal}(\vec{v}) = O(v^{\gamma_2})\htau(p_a) \qquad\text{ and }\qquad \Xcal'(\vec{v}) = O(v^{\gamma_2}).
\end{equation}
Assume that these bounds hold. We continue integrating over the Fourier variables that lie between $p_i$ and $p_j$ until, for $i\neq j$, we end up with either of the following integrals:
\begin{equation}
    \Ycal(\vec{u})\equiv \int \frac{\dd p_i}{(2 \pi)^d} \int \frac{\dd p_a}{(2 \pi)^d} \htau (p_i) \bDc_{\vec{u}} (p_i) \Bcal(p_i, p_{a}) \htau(p_{a}),
\end{equation}
the integral $\overline{\Ycal}(\vec{u})$, which we define to be the same integral but with $\htau(p_i)$ replaced by $\httau (p_i)$, or
\begin{equation}
    \Ycal'(\vec{u}) \equiv \int \frac{\dd p_{i-1}}{(2 \pi)^d} \int \frac{\dd p_i}{(2 \pi)^d} \int \frac{\dd p_{i+1}}{(2 \pi)^d} \htau(p_{i-1}) \Bcal(p_{i-1}, p_{i+1}) \Ccal (p_{i-1}, p_i, p_{i+1}) \htau(p_{i+1})\bDc_{\vec{u}} (p_i).
\end{equation}
Indeed,
\begin{equation}\label{e:bdYs}
    \Ycal(\vec{u}) = O(u^{\gamma_1}), \qquad \overline{\Ycal}(\vec{u}) = O(u^{\gamma_1}) \qquad\text{ and }\qquad \Ycal' (\vec{u}) = T \tbar O(u^{\gamma_1}),
\end{equation}
as we discuss below.
When $i=j$ we integrate over variables from the left and the right until we obtain
\begin{equation}
    \Zcal(\vec{u}, \vec{v}) \equiv \int \frac{\dd p_i}{(2 \pi)^d} \htau(p_i)^2 \bDc_{\vec{u}} (p_i) \uDc_{\vec{v}} (p_i) = \int \frac{\dd p_i}{(2 \pi)^d} \htau_{\vec{u}} (p_i) \htau_{\vec{v}} (p_i),
    \end{equation}
the integral $\underline{\Zcal}(\vec{u},\vec{v})$, or the integral $\overline{\Zcal}(\vec{u},\vec{v})$, depending on the value of $i$, where $\underline{\Zcal}$ and $\overline{\Zcal}$ follow the same definition as $\Zcal$, but with $\htau_{\vec{v}}$ and $\htau_{\vec{u}}$ replaced by $\httau_{\vec{v}}$ and $\httau_{\vec{u}}$, respectively.

This is the right time to mention the case $n=0$, because then, by \refeq{pibarzero} and the Fourier techniques described above we can write
\begin{equation}\label{e:cospi0}
    \sum_{x \in \Zd}[1-\cos(\vec{u} \cdot x)][1-\cos(\vec{v}\cdot x)] \bar{\pi}^{(0)}(x) = \Zcal(\vec{u}, \vec{v}).
\end{equation}

We show below that
\begin{equation}\label{e:bdZ}
    \Zcal (\vec{u}, \vec{v}) = O(u^{\gamma_1} v^{\gamma_2}).
\end{equation}
Very similar proofs can be given for the following bounds:
\begin{equation}\label{e:bdZt}
    \underline{\Zcal}(\vec{u}, \vec{v}) = O(u^{\gamma_1} v^{\gamma_2}) \qquad \text{ and } \qquad \overline{\Zcal}(\vec{u}, \vec{v}) = O(u^{\gamma_1} v^{\gamma_2}).
\end{equation}

When the bounds \refeq{bdXs}, \refeq{bdYs}, \refeq{bdZ} and \refeq{bdZt} hold, it follows that
\begin{equation}
    \Anmij \le  T^{n-3} \tbar^{n+1} O(u^{\gamma_1} v^{\gamma_2}),
\end{equation}
and therefore
\begin{equation}
    \Anij = \sum_{m=0}^{2^{n-1}} \Anmij \le 2^{n-1} T^{n-3} \tbar^{n+1} O(u^{\gamma_1} v^{\gamma_2})
\end{equation}
and finally, by \refeq{Rnij} and Lemma \ref{orderbounds},
\begin{equation}
    \sum_{x \in \Zd} \sum_{n=0}^{\infty} [1 - \cos(\vec{u}\cdot x)][1 -\cos(\vec{v}\cdot x)] \bar{\pi}^{(n)}(x) \le \sum_{n=0}^{\infty} (4n +3)^2 \Anij =  O(u^{\gamma_1} v^{\gamma_2})
\end{equation}
when $\beta$ is sufficiently small, as we set out to prove.

We complete the proof by establishing \refeq{bdZ} and the third bound in \refeq{bdYs}. The two other bounds in \refeq{bdYs} and those in \refeq{bdXs} can be obtained similarly.

Before we start with the proof of \refeq{bdZ}, we briefly indicate how to deal with factors $\httau_{q}(k)$ when they appear.
Define
\begin{equation}
    D_q (x) = [1-\cos(q \cdot x)]D(x).
\end{equation}
Recall the definition of $\widetilde{\tau}$, \refeq{ttaul}. We begin by distributing the weight once more, now over $D$ and $\tau$:
\begin{equation}\label{e:httaulsplit}
    \begin{split}
        [1-\cos(q \cdot x)](D * \tau)(x)    &= \sum_{y \in \Zd}[1-\cos(q \cdot x)]D(y)\tau(x-y)\\
                                            &\le 5 \sum_{y \in \Zd}([1-\cos(q \cdot y)] + [1-\cos(q \cdot (x-y))])D(y)\tau(x-y)\\
                                            &= 5 (D_q * \tau)(x) + 5 (D * \tau_q)(x)
    \end{split}
\end{equation}
where we used \refeq{costbd} for the inequality.
The Fourier transform of $(D_q * \tau)(x)$ can be bounded as follows:
\begin{equation}
    \begin{split}
        \widehat{(D_q * \tau)}(k)   &= \hat{D}_q(k) \htau (k) = \left(\sum_{x \in \Zd} \cos(k \cdot x)[1-\cos(q \cdot x)]D(x) \right)\htau(k)\\
                                    &\le \left(\sum_{x \in \Zd} [1 - \cos(q \cdot x)]D(x)\right) \htau(k) = [1-\hat{D}(q)]\htau(k)= O(q^{\twa}) \htau(k).
    \end{split}
\end{equation}

For the second term of \refeq{httaulsplit}, we observe that $\hat{D}(k) \le 1$, uniformly in $k$, so
\begin{equation}
    \widehat{(D * \tau_q)}(k) = \hat{D}(k) \htau_q (k) \le \htau_q (k).
\end{equation}
Hence, we can bound
\begin{equation}
    \httau_q (k) \le O(q^{\twa}) \htau(k) + 5 \htau_q (k).
\end{equation}
Applying this bound whenever a weighted factor $\httau$ occurs, we can make use of the bounds on weighted and unweighted factors $\htau$ for an upper bound.

Now we give the full proof of \refeq{bdZ}. Recall definition \refeq{taul}. From symmetry of the cosine it follows that that
\begin{equation}\label{e:htaul}
    \htau_q (k) = \tfrac12 \htau(k-q) + \tfrac12 \htau (k+q) - \htau (k) = -\tfrac12 \Delta_{q} \htau (k),
\end{equation}
where $\Delta_q$ is the discrete Laplacian operator with shift $q$.
Therefore, using \refeq{htaul} we obtain
\begin{equation}\label{e:Zuvbd}
        \Zcal (\vec{u},\vec{v})= \int \frac{\dd k}{(2 \pi)^d} \htau_{\vec{u}} (k) \htau_{\vec{v}} (k)  \le \intF{\left\lvert \tfrac12 \Delta_{\vec{u}} \htau(k)\right\rvert \left\lvert \tfrac12 \Delta_{\vec{v}} \htau(k) \right\rvert}.
\end{equation}
Define
\begin{equation}
    \Ck = \frac{1}{1-\hat{D}(k)}.
\end{equation}
Recall \refeq{htaubd}. It follows that
\begin{equation}\label{e:taukbd}
    \htau(k) = O(1) \Ck.
\end{equation}
Hence,
\begin{equation}\label{e:laptautriv}
    \left\lvert \tfrac12 \Delta_q \htau(k) \right\rvert \le O(1) (\Cml + \Ck + \Cpl).
\end{equation}
Define
\begin{equation}
    \Ucal(q, k) = \frac{1}{\hat{C}(q)} \left\{\Cml \Ck + \Ck \Cpl + \Cml \Cpl\right\}.
\end{equation}
From \cite[(5.17)]{HeyHofSak08} we also have the following bound:
\begin{equation}\label{e:laptaunontriv}
    \left\lvert \tfrac12 \Delta_q \htau(k)\right\rvert \le O(1)\, \Ucal (q, k).
\end{equation}
Combining \refeq{laptautriv} and \refeq{laptaunontriv} we obtain an interpolating bound for $\theta \in (0,1)$:
\begin{equation}\label{e:laptauhybrid}
    \left\lvert \tfrac12 \Delta_q \htau (k) \right\rvert \le [1-\hat{D}(q)]^\theta \Ucal (q, k)^\theta [\Cml + \Ck + \Cpl]^{1-\theta}.
\end{equation}

Now we apply \refeq{laptaunontriv} to $\lvert \frac12 \Delta_{\vec{u}} \htau(k) \rvert$ and \refeq{laptauhybrid} with $\theta =\delta$ to $\lvert \frac12 \Delta_{\vec{v}} \htau (k)\rvert$ in \refeq{Zuvbd}. This gives
\begin{eqnarray}
       \nonumber \Zcal (\vec{u}, \vec{v}) &\le&  O(1) [1-\hat{D}(\vec{u})][1-\hat{D}(\vec{v})]^\delta  \int \dd k [\Cmu \Ck + \Ck \Cpu + \Cmu \Cpu] \\
       \nonumber & &\qquad\qquad\times[\Cmv + \Ck + \Cpv]^{1-\delta}\\
       \nonumber  & &\qquad\qquad \times[\Cmv \Ck + \Ck \Cpv + \Cmv \Cpv]^\delta  \\
       \nonumber   &\le& O(1) [1-\hat{D}(\vec{u})] [1-\hat{D}(\vec{v})]^\delta \int  \dd k [\Cmu \Ck + \Ck \Cpu + \Cmu \Cpu] \\
       \nonumber  & &\qquad\qquad\times [\Cmv^{1-\delta} + \Ck^{1-\delta} + \Cpv^{1-\delta}]\\
       \label{e:Clong} & &\qquad\qquad\times [\Cmv^{\delta} \Ck^{\delta} + \Ck^{\delta} \Cpv^{\delta} + \Cmv^{\delta} \Cpv^{\delta}]
\end{eqnarray}
where we used for the second inequality that $(x+y)^\delta \le x^\delta + y^\delta$ for $x,y \ge 0$ and $\delta \in (0,1)$.
The integral contains $27$ distinct product-terms of the function $\hat{C}$ with different shifts and different powers. One term, for instance, is $\Cmu \Ck^{2-\delta} \Cmv^\delta \Cpv^\delta$. To generalize the structure of these terms, we write
\begin{equation}\label{e:Cgeneral}
    \int \Cmu^{a_1} \Cpu^{a_2} \Cmv^{b_1} \Cpv^{b_2} \Ck^{c_1 + c_2} \dd k.
\end{equation}
Here, $c_1$ is the exponent due to the bound on $\tau_{\vec{u}}$, whereas $c_2$ is due to the bound on $\tau_{\vec{v}}$. Note that for every term the following relations hold for the exponents:
\begin{eqnarray}
    \nonumber a_1 + a_2 + b_1 + b_2 + c_1 + c_2 &=& 3 + \delta;\\
    \label{e:exprules} a_1 + a_2 + c_1 &=& 2;\\
    \nonumber b_1 + b_2 + c_2 &=& 1 + \delta.
\end{eqnarray}

Applying Lemma \ref{lem:Holder} to \refeq{Cgeneral} we obtain the upper bound
\begin{multline}\label{e:Csplit}
     \left(\int \Cmu^{3+\delta} \dd k \right)^{\frac{a_1}{3+\delta}}\left(\int \Cpu^{3+\delta} \dd k \right)^{\frac{a_2}{3+\delta}}\left(\int \Ck^{3+\delta} \dd k \right)^{\frac{c_1}{3+\delta}}\\
    \times \left(\int  \Cmv^{3+\delta} \dd k \right)^{\frac{b_1}{3+\delta}}\left(\int  \Cpv^{3+\delta} \dd k \right)^{\frac{b_2}{3+\delta}}\left(\int \Ck^{3+\delta} \dd k \right)^{\frac{c_2}{3+\delta}}.
\end{multline}

Using \refeq{intshift} and \refeq{exprules} we can bound \refeq{Csplit} from above by
\begin{equation}\label{e:Cdouble}
    \left(\int \Ck^{3+\delta} \dd k\right)^{\frac{2}{3+\delta}}\left(\int \Ck^{3+\delta} \dd k\right)^{\frac{1+\delta}{3+\delta}}.
\end{equation}
Now we bound the first factor by a function of $u$, and the second factor by a function of $v$: use \refeq{lbDhat}, from which it follows that $\Ck = O(\lvert k \rvert^{-\twa})$ for any $k \in [-\pi, \pi]^d$.
Whenever $d > (3+\delta)\twa$ and for any $a \in [0,1]$,
\begin{equation}
    \begin{split}
        \int\limits_{[-\pi,\pi]^d} \Ck^{(3+\delta)} \dd k &\le O(1) \int\limits_{|k| \le a} \frac{1}{\lvert k \rvert^{(3+\delta) \twa}} \dd k+ O\left(a^{-(3+\delta)\twa}\right) \int\limits_{|k| \ge a}1 \dd k \\
        & = O\left(a^{(d-(3+\delta) \twa)\wedge 0}\right) + O\left(a^{d-(3+\delta)\twa}\right).
    \end{split}
\end{equation}
Hence we can bound \refeq{Cdouble} by
\begin{equation}\label{e:Ouvsimple}
    O\left(u^{(\frac{2 d}{3+\delta}- 2 \twa)\wedge 0}\, v^{(\frac{(1+\delta)d}{3+\delta} - (1+\delta) \twa)\wedge 0}\right) = O\left(u^{\gamma_1 -\twa} v^{\gamma_2 - \delta \twa}\right).
\end{equation}
Plugging this bound into \refeq{Clong} and using \refeq{lbDhat} again, we obtain
\begin{equation}\label{e:Zcalfinish}
    \Zcal(\vec{u},\vec{v}) = [1-\hat{D}(\vec{u})][1-\hat{D}(\vec{v})]^{\delta} O\left(u^{\gamma_1 - \twa} v^{\gamma_2 - \delta \twa}\right)= O\left(u^{\gamma_1} v^{\gamma_2}\right),
\end{equation}
establishing \refeq{bdZ}.

The Fourier space diagram corresponding to the integrated function in $\Ycal'(\vec{u})$ has two vertices of degree four and only one vertex of degree three, which, unfortunately, is the weighted vertex. As we saw while bounding $\Zcal$, the integral associated to the weighted vertex is only just convergent for $d$ near the critical dimension when it is of degree two. The other two vertices correspond to integrals that are divergent near the critical dimension.

\begin{figure}
  \includegraphics[width=490pt]{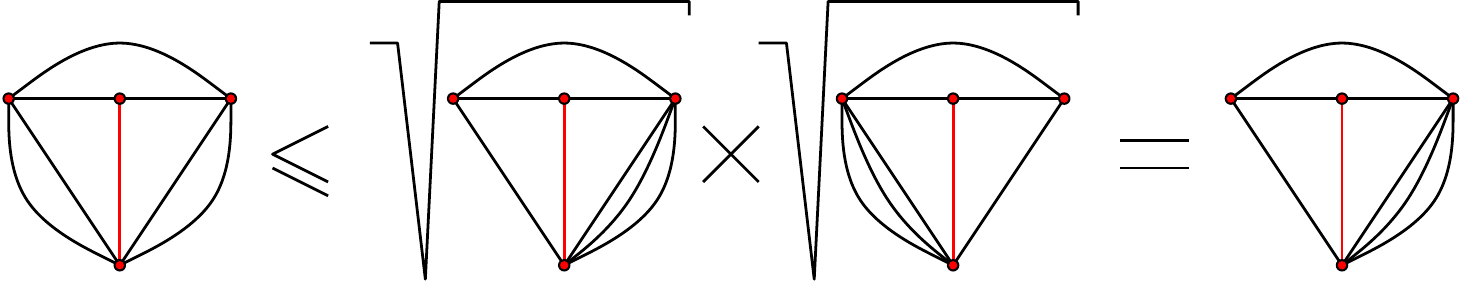}\\
  \caption{A graphic representation of the bound on $\Ycal'(\vec{u})$. The red (vertical) line corresponds to the weighted edge.}\label{pic:Ycalbd}
\end{figure}

However, the diagram has three integrated variables and eight functions, so we should be able to bound it by two triangles and a weighted bubble. To see this, we need to bound the integral by something simpler before we evaluate it. We use the Cauchy-Schwarz inequality for this. Roughly speaking, using the Cauchy-Schwarz inequality and the symmetry of the integral under relabeling allows us to bound the diagram by the same diagram with one factor $\htau(p_{i-1})$ replaced by a factor $\htau(p_{i+1})$. See Figure \ref{pic:Ycalbd} for an illustration of this.

Applying the bound described above, and by positivity of the $\htau$-functions, we obtain
\begin{eqnarray}
  \nonumber  \Ycal'(\vec{u}) & \le &\int \frac{\dd p_{i-1}}{(2 \pi)^d} \int \frac{\dd p_i}{(2 \pi)^d} \int \frac{\dd p_{i+1}}{(2 \pi)^d} \widetilde{\Bcal}(p_{i-1}, p_{i+1}) \Ccal (p_{i-1}, p_i, p_{i+1}) \htau(p_{i+1})^2 \bDc_{\vec{u}} (p_i)\\
  \nonumber  &=& \int \frac{\dd p_i}{(2 \pi)^d} \int \frac{\dd p_{i+1}}{(2 \pi)^d} \htau_{\vec{u}} (p_i) \htau(p_i -p_{i+1}) \htau(p_{i+1})^3 \int \frac{\dd p_{i-1}}{(2 \pi)^d} \hat{D}(p_{i-1}) \htau(p_{i-1}) \htau(p_{i-1} -p_i) \htau(p_{i-1}-p_{i+1}) \\
    & \le& T \int \frac{\dd p_{i+1}}{(2 \pi)^d} \htau(p_{i+1})^3 \int \frac{\dd p_i}{(2 \pi)^d} \htau_{\vec{u}} (p_i) \htau(p_i -p_{i+1})\\
  \nonumber  & \le&  T O(u^{\twa}) \int \frac{\dd p_{i+1}}{(2 \pi)^d} \htau(p_{i+1})^3 \le T \tbar O(u^{\twa}) \le T \tbar O(u^{\gamma_1}).
\end{eqnarray}
The second to fourth inequality follow from a calculation similar to \refeq{bdAB} and $\Zcal(\vec{u}, \vec{v})$. The final bound is just there to fit the statement of \refeq{Wuvbound}. This completes the proof of Proposition \ref{prop:TermBounds}(ii).
\qed

\subsection{About the proof of Lemma \ref{lem:Phibd}}
The equality \refeq{Pidelta} in Lemma \ref{lem:Phibd} follows immediately from \refeq{pi0d1}, \refeq{pidiag1}, \refeq{PiPhidef} and the proof of Proposition \ref{prop:TermBounds}(ii).

The equality \refeq{Phidelta} in Lemma \ref{lem:Phibd} can be proved in the same way as Proposition \ref{prop:TermBounds}(ii), but with much less bookkeeping, so we do not give it. Heuristically, the validity of the claim can be understood by noting that the diagrams $\bar{\phi}^{(n)}(x,y)$ are like $\bar{\pi}^{(n)}(x)$ diagrams with an extra point $y$ placed on either the last or second-to-last upper path element (cf. Figure \ref{fig:DiagPhiPi}). From \refeq{laptaunontriv} it can be seen that in Fourier space, adding a point to a path element has more or less the same effect as having a `heavy' weight on that path element. Therefore, the diagrams $\bar{\phi}^{(n)}(x,y)$ with the small weight $\lvert x - y\rvert^{\delta}$ can be bounded in a similar way as the diagrams $\bar{\pi}^{(n)}(x)$ with the weight $\lvert x \rvert^{\twa + \delta}$, and hence the bounds should also be similar. Following the proof of Proposition \ref{prop:TermBounds}(ii) confirms that this is the case. In the course of this proof, \refeq{Phinodelta} in Lemma \ref{lem:Phibd} also follows naturally.

\vspace{0.5cm}
{\bf Acknowledgements.}
This work is supported by the Netherlands Organization for Scientific Research (NWO). MH and RvdH thank the Institut Henri Poincar\'e Paris for kind hospitality during their visit in October 2009.
\begin{small}
\bibliographystyle{abbrv}
\bibliography{TimsBib}
\end{small}
\end{document}